\documentclass[preprint,12pt]{elsarticle}

\usepackage[utf8]{inputenc}
\usepackage{amsmath}
\usepackage{amstext} 
\usepackage{array}
\usepackage{accents}
\usepackage{hyperref}
\usepackage{cleveref}
\usepackage{commath}

\usepackage{graphicx}
\usepackage{float}
\usepackage[colorinlistoftodos,prependcaption,textsize=tiny]{todonotes}
\usepackage{subcaption}
\usepackage{tikz}
\usepackage{bm}
\usetikzlibrary{arrows}
\usepackage{comment}
\usepackage[numbers]{natbib}

\newcommand{\nfa}{\textrm{NFA}}

\newcommand{\bx}{\bm{x}}
\newcommand{\bn}{\bm{n}}
\newcommand{\bv}{\bm{v}}

\newcommand{\bA}{\bm{A}}

\newcommand{\bu}{\bm{u}}
\newcommand{\sett}[2]{\left\{{#1}\,:~{#2}\right\}}

\newcommand{\sym}{\nabla_{sym}}

\newcommand{\gd}{\nabla}

\DeclareMathOperator{\dv}{div}

\newcolumntype{L}{>{$}l<{$}} 

\usepackage{fancyhdr}
\usepackage{amsmath}
\usepackage{amssymb}
\usepackage{tikz}
\usepackage{tikz-cd}
\usetikzlibrary{cd,arrows}
\tikzcdset{scale cd/.style={every label/.append style={scale=#1},
cells={nodes={scale=#1}}}}
\usetikzlibrary{decorations.markings}
\tikzset{
invisible/.style={opacity=0},
visible on/.style={alt={#1{}{invisible}}},
alt/.code args={<#1>#2#3}{%
\alt<#1>{\pgfkeysalso{#2}}{\pgfkeysalso{#3}}%
}
}
\usepackage{ragged2e}
\usepackage{tikz}
\usetikzlibrary{positioning}

\usepackage{alltt}
\usepackage{url}
\usepackage[absolute,overlay]{textpos}
\usepackage[utf8]{inputenc}
\usepackage[T1]{fontenc}
\usepackage[english]{babel}
\usepackage{mathtools,amsmath,amsthm,amssymb}
\usepackage{extarrows}
\usepackage{transparent}
\usepackage{eucal}
\usepackage{color}
\usepackage[inkscapeformat=png]{svg}
\usepackage{graphicx}
\usepackage{placeins}       
\usepackage{tikz}
\usepackage[export]{adjustbox}

\usepackage{ifthen}
\usetikzlibrary{backgrounds}
\usetikzlibrary{shapes}
\usetikzlibrary{snakes}
\usetikzlibrary{patterns}        
\usepackage{pgfplots}
\usetikzlibrary{intersections, pgfplots.fillbetween}
\usepackage{pgfplotstable}
\usepgfplotslibrary{fillbetween}
\usepgfplotslibrary{groupplots}
\usepackage{layouts}
\usepackage{cancel}
\usepackage{multimedia}
\usepackage{qrcode}

\usepackage{amssymb}
\usepackage{cleveref}

\newcommand{\R}{\ensuremath{\mathbb{R}}} 

\journal{Journal}

\begin{document}

\begin{frontmatter}

\title{Preconditioning for a Cahn-Hilliard-Navier-Stokes model for morphology formation in organic solar cells}
\author[TUchem]{Pelin \c{C}\.{i}lo\u{g}lu}
	\ead{pelin.ciloglu@mathematik.tu-chemnitz.de}
 
\author[AugsU]{Carmen Tretmans}
	\ead{carmen.tretmans@uni-a.de}

\author[Heide]{Roland Herzog}
	\ead{roland.herzog@iwr.uni-heidelberg.de}

\author[AugsU]{Jan-F. Pietschmann}
	\ead{jan-f.pietschmann@uni-a.de}
 
\author[TUchem]{Martin Stoll}
	\ead{martin.stoll@mathematik.tu-chemnitz.de}

\cortext[cor1]{Corresponding author}    
\address[TUchem]{Faculty of Mathematics, University of Technology Chemnitz, 09107, Chemnitz, Germany}
\address[AugsU]{Faculty of Mathematics, University of Augsburg, 86159, Augsburg, Germany}
\address[Heide]{Interdisciplinary Center for Scientific Computing (IWR), Heidelberg University, 69120, Heidelberg, Germany}

\begin{abstract}
We present a  model for the morphology evolution of printed organic solar cells which occurs during the drying of a mixture of polymer, the non-fullerene acceptor and the solvent. Our model uses a phase field approach coupled to a Navier-Stokes equation describing the macroscopic movement of the fluid. Additionally, we incorporate the evaporation process of the solvent using an Allen-Cahn equation. 
The model is discretized using a finite-element approach with a semi-implicit discretization in time. The resulting (non)linear systems are coupled and of large dimensionality. We present a preconditioned iterative scheme to solve them robustly with respect to changes in the discretization parameters. We illustrate that the preconditioned solver shows parameter-robust iteration numbers and that the model qualitatively captures the behavior of the film morphology during drying.
\end{abstract}



\begin{keyword}
Preconditioning \sep  phase–field models \sep  organic solar cells \sep Navier-Stokes \sep  Cahn–Hilliard \sep  finite element analysis


\end{keyword}

\end{frontmatter}


\section{Introduction}
Research on the fabrication of organic photovoltaics (OPV) has gained much momentum in recent years due to its high potential in terms of a low-cost and energy efficient production process, combined with favourable properties like flexibility and lightweightness. Compared to silicon-based alternatives, OPVs still lack in various aspects like efficiency, stability, and lifespan, but have gained further traction recently. In particular, the use of non-fullerene acceptors (NFA) has the potential to make organic solar cells a competitive alternative to silicon-based cells with an environmental advantage \cite{weitz_revealing_2023, wopke_traps_2022}. NFAs are especially beneficial compared to the previously used fullerene acceptors due to their highly tunable photoelectric properties \cite{zhang2022renewed}. Key to improving upon existing issues is the active layer of OPVs, as it is the most important component effecting the performance of the cell. The layer is a thin film consisting of pure regions of an organic electron donor and acceptor. These thin films are formed by a solvent-based fabrication process.

An organic blend consisting of a polymer acting as an electron donor, an NFA acting as an electron acceptor, and a solvent are deposited onto a substrate, using, for example, spin coating techniques \cite{janssen2007optimization}. The polymer and the NFA undergo a spontaneous phase separation as the solvent evaporates, resulting in a thin layer containing pure donor and acceptor regions that act as active layer. The resulting morphology within the thin film is key to the device performance. 
However, the morphology formation is a complex process being highly sensitive to the processing conditions as well as specifics of the polymer, NFA, and solvent blends. 
Finding the optimal production conditions experimentally would be very time-consuming and is therefore virtually impossible. Therefore, accurate modelling of the morphology formation and subsequent in-silico numerical simulations are an indispensable tool towards a highly optimized OPV production process. 

The mathematical modelling of the aforementioned morphology formations has previously been studied in \cite{Bergermann2023,Wodo2012_height}.
In what follows, we will present a novel model that links fluid dynamical influences, phase separation of acceptor and donor regions, and solvent evaporation. This continuous problem is then discretized using a finite element discretization \cite{brenner2008mathematical,elman2014finite}, which results in large-scale linear and nonlinear systems of equations. These systems typically require the careful choice of a preconditioner to avoid a large number of iterations. More precisely, the linear systems we encounter are in saddle-point form
$$
\begin{bmatrix}
     \mathbf{A} &\mathbf{B}^\top\\
    \mathbf{B}&0
\end{bmatrix}
\begin{bmatrix}
\mathbf{x}_1\\
\mathbf{x}_2
\end{bmatrix}
=\begin{bmatrix}
\mathbf{b}_1\\
\mathbf{b}_2
\end{bmatrix},
$$
where $\mathbf{A}\in \R^{n,n}, \mathbf{B}\in\R^{m,n}$ are matrices representing discretized continuous operators. Writing this system as $\mathcal{A}\mathbf{x}=\mathbf{b}$ we will  construct a preconditioning matrix $\mathcal{P}$ and instead solve the preconditioned system $\mathcal{P}^{-1}\mathcal{A}\mathbf{x}=\mathcal{P}^{-1}\mathbf{b}.$ 
%
In the context of Cahn-Hilliard equations preconditioners have already been suggested in \cite{bosch2015preconditioning,bosch2014fast, boyanova2012efficient,boyanova2014efficient} where these are typically block preconditioners exploiting the saddle-point structure of the discretized equations and then relying on efficient approximations of the individual blocks or components. The problem becomes more challenging when the phase field equations are coupled to the Navier-Stokes equations, for which the Cahn-Hilliard-Navier-Stokes equations \cite{boyer2010cahn,abels2013incompressible} are the main representative. Yet also in this situation, preconditioners have been proposed in \cite{kay2007efficient,JBosch_CKahle_MStoll_2018}. 
One of the main contribution of this work is the construction of an appropriate preconditioning technique for the model of OPVs which includes evaporation effects. 

The remainder of this paper is structured as follows. In \Cref{sec::model} we derive the model describing the morphology formation by coupling a phase separation process to an evaporation equation and Navier-Stokes equations. We then derive the spatial and temporal discretization in \Cref{sec::discretization} followed by the introduction of several structured preconditioners for the various coupled PDE models in \Cref{sec::precon}. In \Cref{sec::results} we illustrate the performance of the suggested preconditioners for various setups of the OPV model.

\section{Model equations} \label{sec::model}
\subsection{Governing equations}
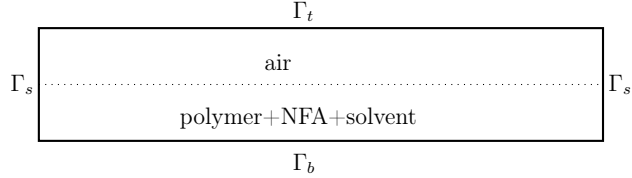
\begin{figure}[h!]
    \centering
    \hspace{1.5cm}
    \begin{tikzpicture}[scale=0.75]
    \draw[thick] (1,0) rectangle (11,2);
    \draw[dotted] (1,1) -- (11,1);
    \node[scale=0.75, text width=3cm]  at (6.5,1.4) {air};
    \node[scale=0.75, text width=3cm]  at (5,0.4) {polymer+NFA+solvent};
    \node[scale=0.75, text width=3cm]  at (7,2.3) {$\Gamma_t$};
    \node[scale=0.75, text width=3cm]  at (7,-.4) {$\Gamma_b$};
    \node[scale=0.75, text width=3cm]  at (2,1) {$\Gamma_s$};
    \node[scale=0.75, text width=3cm]  at (12.6,1) {$\Gamma_s$};
    \end{tikzpicture}
\caption{Sketch of the domain.}
\label{fig_domain_sketch}
\end{figure}

To simulate the morphology evolution during the drying of the thin film, we will derive a model capturing two main features. First, the behavior of the liquid organic blend needs to be described. Secondly, the evaporation of the solvent into the surrounding air needs to be included. We will consider the open domain  $\Omega  \subset \mathbb{R}^n, \; n\in {1, 2,3} $ with Lipschitz boundary $\partial \Omega$ as shown in Figure~\ref{fig_domain_sketch}. Initially, the blend of polymer, NFA, and solvent is deposited on a substrate, indicated by the bottom boundary $\Gamma_b$. The blend is surrounded by the ambient air, taking up the top part of the domain, which is enclosed at the top by the upper boundary $\Gamma_t$. The side boundaries are denoted by $\Gamma_s$.

The morphology is described by a phase field model consisting of four different material phases, namely polymer, NFA, solvent, and the surrounding air. They are described by the scalar fields 
	\begin{align*}
		\phi_p, \phi_{\nfa}, \phi_s, \phi_a \colon \Omega \times [0,t_{\mathrm{max}}] \rightarrow [0,1], 
		\end{align*}
    which are volume fractions and thus satisfy the pointwise constraint 
    \begin{align*}
        \phi_p + \phi_{\nfa} + \phi_s + \phi_a = 1,\quad \bx \in \Omega, \text{ and every }t \in [0,t_{\mathrm{max}}].
    \end{align*}  
To include evaporative effects of solvent into the ambient air, we follow the approach presented by Ronsin et al.~\cite{brabec2021_evaporation, brabec2020_evaporation}, and extend the system to allow species to change from the liquid state to a gaseous state. For this purpose, an additional order parameter
\begin{align*}
    \Phi_v \colon \Omega \times [0,t_{\mathrm{max}}] \rightarrow [0,1]
\end{align*}
is introduced denoting the fraction of the mixture in vapor phase. Here, $\Phi_v = 1$ describes the gas phase, $\Phi_v = 0$ the liquid phase. The motion of the four different species can be described by the separate velocities $\bv_i$, $i \in \{ p, \nfa, s, a \}$. Following the approach of \mbox{Abels et al.~\cite{Garcke2011_CHNS},} one can divide the motion of each species into the macroscopic motion of the fluid, describing the motion of the (volume averaged) fluid mixture $\bv = \sum_i \phi_i \bv_i$ and the microscopic motion of the individual components $(\bv - \bv_i)$ describing the phase separation of the components due to thermodynamic driving forces. The macroscopic motion of the fluid is governed by the laws of fluid dynamics, the microscopic motion of the fluid can be described by thermodynamic reasoning and is based in the minimisation of a free energy functional $f(\boldsymbol{\phi}, \gd\boldsymbol{\phi} , \Phi_v, \gd \Phi_v)$. Here, $\boldsymbol{\phi}$ and $\gd \boldsymbol{\phi}$ denote the set of $\phi_i$ and $\gd \phi_i$, $i \in \{ p, \nfa, s, a \}$, respectively. The resulting system consists of Navier-Stokes equations describing the macroscopic movement of the fluid, coupled to the Cahn-Hilliard equations describing the microscopic movement. An elaborate derivation of this system can be found in \cite{Garcke2011_CHNS}, yet not in the context of OPVs and without evaporation effects.

The evaporation of the solvent is described by the evolution of the order parameter $\Phi_v$, combined with the solvent volume fraction $\phi_s$. Following the approach by Ronsin et al.~\cite{brabec2021_evaporation}, evaporation of solvent is driven by two processes. First, solvent changes its state of matter from liquid to gaseous at the top surface of the liquid due to thermodynamic forces, forming a gas layer above the liquid. This process continues until the partial pressure of the solvent, measured in the gas layer above the liquid, approaches the equilibrium vapor pressure. This equilibrium is reached once the thermodynamic forces, transforming solvent from the gaseous state to the liquid state and vice versa, equate, and is a material specific property. The second process is the diffusion of gaseous solvent into the surrounding air. This process steers the partial pressure at the fluid-gas interface away from the equilibrium vapor pressure, which in turn results in transformation of more solvent from the liquid to the gaseous state. As mentioned in \cite{brabec2020_evaporation}, the second process is in fact the one limiting evaporation. It comes to a halt when the partial solvent pressure equals the constant equilibrium vapor pressure, i.e. the vapor is saturated. 

To model the diffusion of vapor into the surrounding environment, not specifically included in the domain, we define an outflux $j_{out}$ of solvent at the top of the domain $\Gamma_t$. The outflow of solvent is coupled to a respective inflow of air. The evaporation process can now continue until no more solvent is present in the simulation domain. The thermodynamic forces causing the change of state of matter are, as for the phase separation, based on the minimization of a free energy functional $F = \int_\Omega f \, dx $. The resulting equation governing the vapor order parameter is not mass conserving, and thus yields an Allen-Cahn equation. The Allen-Cahn equation can be coupled to the Navier-Stokes equations in a thermodynamically consistent way using a microforce balance, see Abels et al.~\cite{Garcke2011_CHNS}, with microforces derived by Gurtin \cite{gurtin_generalized_1996}. 

Various assumptions and simplifications are used in the model for practical purposes. While we allow the species to have different pure densities $\Tilde{\rho}_i$, ${ i \in \{ p, \nfa, s, a \} }$, we assume these densities to be constant and independent of the species' liquid or gaseous state. This is a strong assumption that will considerably affect the representation of the vapor phase. However, this does not directly affect the model behavior in the liquid phase, which is our main interest. Secondly, solvent is replaced by air in the simulation domain in the description of the evaporation process. This is achieved by creating an outward flux of (vapor) solvent and a respective inward flux of (vapor) air at the boundary of the domain. Since this flux should be zero for liquid solvent, we restrict this outward flux to the top boundary of the domain, where the vapor phase can be assumed. Furthermore, we consider the process to be isothermal with constant temperature $T$, allowing for further simplifications of the model. Taking these aspects into account yields the coupled equations
\begin{align}\label{eqn:govern}
    \span \rho \frac{\partial \bv}{\partial t} + ((\rho \bv - \sum_{i}\frac{V_0 \Tilde{\rho}_i}{R T} M_i\nabla \mu_i) \cdot \nabla ) \bv  = \dv (2 \eta \sym \bv) \nonumber \\
    \span - \sum_i \dv \left( \gd \phi_i \otimes \frac{\partial f}{\partial \gd \phi_i} \right)- \dv \left( \gd \Phi_v \otimes \frac{\partial f}{\partial \gd \Phi_v} \right) - \gd p \quad \text{ in } \Omega \times [0,t_{\mathrm{max}}],\nonumber\\ 
    \dv \bv & = 0 \quad \text{ in } \Omega \times [0,t_{\mathrm{max}}], \nonumber\\  
    \frac{\partial \phi_i}{\partial t} + \bv \cdot \gd \phi_i  & = \frac{V_0}{RT} \dv (M_i \nabla \mu_i)  \quad \text{ in } \Omega \times [0,t_{\mathrm{max}}] \text{ for } i \in \{p, \nfa, s, a \}, \nonumber\\
	\mu_i &	= \frac{\partial f}{\partial \phi_i} -  \dv \frac{\partial f}{\partial \gd \phi_i} \quad \text{ in } \Omega \times [0,t_{\mathrm{max}}] \text{ for } i \in \{p, \nfa, s, a \}, \nonumber\\ 
    \frac{\partial \Phi_v}{\partial t} + \bv \cdot \gd \Phi_v & = - \frac{V_0}{RT}M_v \left( \frac{\partial f}{\partial \Phi_v} -  \dv \frac{\partial f}{\partial \gd \Phi_v} \right) \quad \text{ in } \Omega \times [0,t_{\mathrm{max}}]. 
\end{align}
Here, all summations are over the set $\{ p, \nfa, s, a \}$. The symmetric gradient is denoted by $\sym ( \bv ) = \frac{1}{2} \left( (\gd \bv ) + ( \gd \bv )^T \right) $. The divergence of a matrix is defined row-wise as $(\dv \bA )_{i} = \sum_j \frac{\partial A_{ij}}{\partial x_j}$.  We introduce the chemical potentials $\mu_i(\bx, t)$, the density of the fluid mixture $\rho(\bx, t) = \sum_i \phi_i(\bx, t) \Tilde{\rho}_i$, and pressure $p(\bx, t)$. The constant parameters $V_0$, $R$, and $T$ represent the reference molar volume, gas constant, and temperature, respectively. The molecular mobility coefficients $M_i$, the vapor mobility coefficient $M_v$, and viscosity $\eta$ generally depend on the composition of the fluid blend $\boldsymbol{\phi}$ as well as the vapor fraction $\Phi_v$. To avoid additional nonlinearities within the model, all three parameters are chosen constant in this setting.
We apply no-slip boundary conditions to the velocity and no flux boundary conditions for all volume fractions. An exception to the latter is the outflux of solvent and inflow of air $j_{out}$ at the top boundary $\Gamma_t$ due to the evaporation effects. For the chemical potential as well as the vapor order parameter we impose no flux conditions,
\begin{align}\label{eqn:bc}
        \bv & = 0  \quad \text{on} \quad \partial \Omega \times [0,t_{\mathrm{max}}], &   \nonumber \\
        \frac{\partial \phi_i}{\partial \bn } & = 0 \quad \text{on} \quad \partial \Omega \times [0,t_{\mathrm{max}}] \text{ for } i \in \{ p, \nfa \}, \nonumber \\
        \frac{\partial \phi_s}{\partial \bn } & = \begin{cases}
        j_{out}  \quad \text{on} \quad \Gamma_{top}  \times [0,t_{\mathrm{max}}], \\
        0  \quad \text{on} \quad \partial \Omega \setminus \Gamma_{top}  \times [0,t_{\mathrm{max}}],
        \end{cases} \nonumber  \\
          \frac{\partial \phi_a}{\partial \bn } &=  \begin{cases}
        -j_{out}  \quad \text{on} \quad \Gamma_{top}  \times [0,t_{\mathrm{max}}] ,\\
        0  \quad \text{on} \quad \partial \Omega \setminus \Gamma_{top}  \times [0,t_{\mathrm{max}}],
         \end{cases}  \nonumber  \\
         \frac{\partial \Phi_v}{\partial \bn } & = 0  \quad \text{on} \quad \partial \Omega \times [0,t_{\mathrm{max}}] , &  \nonumber \\
         \frac{\partial \mu_i}{\partial \bn } &= 0 \quad \text{on} \quad \partial \Omega \times [0,t_{\mathrm{max}}] \text{ for }  i \in \{p, \nfa, s, a \} . 
    \end{align}
Following the Hertz-Knudsen theory described by Ronsin et al.~\cite{brabec2021_evaporation}, the flux at the top of the boundary $\Gamma_{top}$ is defined as 
\begin{align*}
    j_{out} = \alpha \sqrt{\frac{V_0}{2 \pi RT} \frac{N_s}{\Tilde{\rho}_s}} P_0 (\phi^{vap}_s - \phi_s^\infty)
\end{align*}
for a given flow rate $\alpha$, solvent volume fraction in the ambient air $\phi_s^\infty$, and reference pressure $P_0$. Here, $\phi_s^{vap}$ is equal to the mean volume fraction of solvent in the gas phase. Since one can assume high diffusivity in the gas phase, one can assume approximately constant solvent fraction within the gas phase. We simplify $\phi_s^{vap}$ by setting its value equal to the solvent volume fraction at a reference point at the top of the domain.

Imposing the initial conditions 
    \begin{align} \label{eqn:ic}
        \bv(\cdot, 0) &= \bv_0(\cdot) \quad \text{in} \quad \Omega ,  \nonumber \\
        \phi_i (\cdot, 0) &= \phi_{i,0}(\cdot) \quad \text{in} \quad \Omega \text{ for } i \in \{ p, \nfa, s, a \} ,\nonumber \\
        \Phi_i (\cdot, 0) &= \Phi_{v,0}(\cdot) \quad \text{in} \quad \Omega 
\end{align} 
completes the governing equations. 

\subsection{Free energy functional}
Phase separation and the evaporation of solvent are driven by thermodynamic forces, minimizing a free energy functional $F = \int_\Omega f \, dx $. Key to an accurate representation of the morphology evolution is hence a detailed representation of the functional. Locally, the free energy consists of entropy and enthalpy contributions. While the entropy favours mixing of different components, enthalpy contributions might induce phase separation due to stronger intermolecular bonds. We capture the change of the free energy density compared to the unmixed, liquid ground state of all components. Following the Flory-Huggins theory, see e.g \cite{MDijk_AWakker_1997}, the change of energy density due to mixing of polymers with molar volume $V_i = N_i V_0$ is given by
\begin{align*}
    f^{loc}_{liq}(\boldsymbol{\phi}) & = \frac{RT}{V_0}\left( \sum_i \frac{\phi_i}{N_i} \ln \phi_i + \sum_i \sum_{j<i} \chi_{i,j} \phi_i \phi_j \right).  
\end{align*}
Here, the first term sums over all species $i \in \{ p, \nfa, s, a\}$ and represents the entropy of mixing. The second term sums over all pairs $i \neq j$ with $i, j \in \{ p, \nfa, s, a \}$ and captures the change of enthalpy with the interaction parameters $\chi_{ij}$. In the gas phase, all gases are assumed to be ideal with the same molecular size, canceling the enthalpy contribution. Compared to the unmixed liquid ground state, the energy density due to vaporization is captured by $- \phi_i \ln \phi_{sat(i)}$, where $\phi_{sat(i)}$ are the saturated vapor volume fractions, equal to the saturated partial pressure of fluid $i$ in the surrounding air. The energy density in the gas phase is given by 
\begin{align*}
    f^{loc}_{gas}(\boldsymbol{\phi}) & = \frac{RT}{V_0}\left( \sum_i \phi_i \ln \frac{\phi_i}{\phi_{sat(i)}} \right).
\end{align*}
To obtain the total change in energy density, we interpolate between the liquid and gas phase using the order parameter $\Phi_v$, 
\begin{align*}
    f^{loc}(\boldsymbol{\phi}, \Phi_v) &= \left( 1- p(\Phi_v)\right) f^{loc}_{liq}(\boldsymbol{\phi}) + p(\Phi_v) f^{loc}_{gas} (\boldsymbol{\phi})  .
\end{align*}
Here, we use the smooth interpolation function 
\begin{align*}
    p(\Phi_v) = \Phi_v^2(3-2\Phi_v) . 
\end{align*}
In addition to the local part of the free energy density, a nonlocal term is introduced penalizing the field gradients,
\begin{align*}
        f^{nonloc}(  \gd \phi_i , \gd \Phi_v) & = \sum_i \frac{\epsilon_i^2}{2} (\gd \phi_i)^2 + \frac{\lambda^2}{2}(\gd \Phi_v)^2, 
\end{align*}
where $\epsilon_i$ and $\lambda$ are interface parameters related to the width of the transition layers. Combining all contributions, the change of free energy density reads 
\begin{align*}
        f(\boldsymbol{\phi},  \gd  \boldsymbol{\phi}, \Phi_v, \gd \Phi_v) =& \left( 1- p(\Phi_v)\right) f^{loc}_{liq}(\boldsymbol{\phi}) + p(\Phi_v)  f^{loc}_{gas} (\boldsymbol{\phi}) \\
        & + f^{nonloc}( \gd \boldsymbol{\phi}, \gd \Phi_v) .
\end{align*}

\subsection{Nondimensionalized governing equations}
In order to obtain a nondimensionalized version of the governing equations \eqref{eqn:govern}, we fix $x_c$, $t_c$, and $m_c = \frac{RT}{V_0}x_c t_c^2$. Setting $\phi_i  = \phi_i^*$, $\Phi_v = \Phi_v^*$, $\rho = \rho^* \frac{m_c}{x_c^3}$, $\bv = \bv^* \frac{x_c}{t_c}$, $\mu_i = \mu_i^* \frac{m_c}{x_c t_c^2} $, and $p = p^*\frac{m_c}{x_c t_c^2}$, substituting the nonlocal part of the free energy density, and omitting the asterisk, the dimensionless form of the governing equations read 
\begin{align}
\label{eqn:model}
       & \rho \frac{\partial \bv}{\partial t} + ((\rho \bv - \sum_{i} \alpha_i \nabla \mu_i) \cdot \nabla ) \bv = \dv (2 \hat{\eta} \sym \bv)   - \sum_i \beta_i \dv \left( \gd \phi_i \otimes \gd \phi_i \right) \span \span \nonumber \\
       & \qquad - \beta_v \dv \left( \gd \Phi_v \otimes \gd \Phi_v \right)- \gd p   && \text{ in } \Omega \times [0,t_{\mathrm{max}}], \nonumber  \\ 
        &\dv \bv = 0   && \text{ in } \Omega \times [0,t_{\mathrm{max}}], \nonumber  \\  
        &\frac{\partial \phi_i}{\partial t} + \bv \cdot \gd \phi_i  =\gamma_i \dv ( \nabla \mu_i)   && \text{ in } \Omega \times [0,t_{\mathrm{max}}] \text{ for } i \in \{ p, \nfa, s, a \} , \nonumber\\  
       & \frac{\partial \Phi_v}{\partial t} + \bv \cdot \gd \Phi_v = - \gamma_v \frac{\partial \hat{f}}{\partial \Phi_v} + \delta_v \nabla^2 \Phi_v    && \text{ in } \Omega \times [0,t_{\mathrm{max}}] , \nonumber\\  
	&	\mu_i 	= 	\frac{\partial \hat{f}}{\partial \phi_i} - \beta_i \nabla^2 \phi_i  && \text{ in } \Omega \times [0,t_{\mathrm{max}}] \text{ for } i \in \{ p, \nfa, s, a \}, 
\end{align}
with the rescaled free energy density change $\hat{f} = \frac{V_0}{RT} f$, rescaled viscosity $\hat{\eta} = \frac{V_0}{RT} \frac{\eta}{t_c}$, and dimensionless parameters $\alpha_i$, $\beta_i$, $\beta_v$, $\gamma_i$, and $\gamma_{v}$. For readability, we use the notation $f$ and $\eta$ instead of $\hat{f}$ and $\hat{\eta}$ throughout the remainder of this paper. The nondimensionalized governing equations are completed with the nondimensionalized form of the initial~\eqref{eqn:ic} and boundary conditions~\eqref{eqn:bc}.

\section{Spatial \& temporal discretization}
\label{sec::discretization}
We use the finite element method \cite{brenner2008mathematical,strang2008analysis} to discretize the governing equations~\eqref{eqn:model}.
When using the finite element method with a primal variational formulation, we can avoid the limitations related to the continuity of basis functions by solving equations using the split form, see, e.g. \cite{Wodo2012_height}. The weak formulation of the Navier-Stokes part of equations~\eqref{eqn:model} reads 
\begin{subequations}\label{eqn:weak}
\begin{align}\label{eqn:weakns}
   &  \int \limits_{\Omega} \rho \frac{\partial \bv}{\partial t} \bu \,dx + \int \limits_{\Omega} 2\eta \sym \bv : \sym \bu \, dx +  \int \limits_{\Omega} ((\rho \bv - \sum_{i}\alpha_i \nabla \mu_i) \cdot \nabla ) \bv   \bu \, dx     \nonumber \\
   & \hspace{4.2cm} - \int \limits_{\Omega}  p \, \dv \bu \,dx  =  \int \limits_{\Omega}  \mathcal{G} :  \nabla \bu \, dx,  \\
   & \int \limits_{\Omega}  \dv \bv \, q  \,dx = 0, 
\end{align}
for all $\bu \in H^1_0(\Omega)^n$ and all $q \in L^2_{\diamond}(\Omega)$, where we define the space of square integrable functions with mean value zero,
\begin{align*}
    L^2_{\diamond}(\Omega) := \sett{ q \in L^2(\Omega) } {\int_{\Omega} q \,dx = 0 }.
\end{align*}
We furthermore denote  $\mathcal{G} = \sum_i \beta_i \left( \nabla \phi_i \otimes  \nabla\phi_i \right) +\beta_v  \left( \nabla \Phi_v \otimes \nabla \Phi_v \right) $. For the Cahn-Hilliard equations we obtain  
\begin{align}\label{eqn:weakch}
    & \int \limits_{\Omega} \left(\frac{\partial \phi_i}{\partial t } +\bv \cdot \nabla \phi_i \right) q_j \,dx +  \int \limits_{\Omega} \gamma_i \nabla \mu_i \cdot \nabla q_j \, dx =  0, \\ 
     & \int \limits_{\Omega} \mu_i \, w_j \,dx  - \int \limits_{\Omega} \beta_i \nabla \phi_i \cdot \nabla w_j \, dx  = \int \limits_{\Omega} \frac{\partial f}{\partial \phi_i} w_j \, dx  - \int \limits_{\Gamma_t} \beta_i\nabla \phi_i \cdot \mathbf{n} \, w_j \, ds ,
\end{align}
for all $q_j \in H^1(\Omega)$, $w_j \in H^1(\Omega)$, $i \in \{\hbox{p, \nfa, s, a} \}$ and $j = \{\hbox{p, \nfa, s, a}\}$. We point out that the Cahn-Hilliard equations are always treated in coupled form relating $\phi_i$ and $\mu_i.$ This is important when discussing the solution of the discretized systems in matrix form. Lastly, we obtain for the Allen-Cahn equation
\begin{align} \label{eqn:weakac}
     \int \limits_{\Omega} \left(\frac{\partial \Phi_{v}}{\partial t } + \bv \cdot \nabla \Phi_{v} \right) q_v \,dx + &  \int \limits_{\Omega}  \delta_v \nabla \Phi_{v} \cdot \nabla q_v \, dx 
      & = - \int \limits_{\Omega}  \gamma_v \frac{\partial f }{ \partial \Phi_v } q_v \, ds, 
\end{align}
\end{subequations}
for all $q_v \in H^1(\Omega)$. 

Let $\{\mathcal{T}_h\}_h$ be a family of shape-regular simplicial triangulations of the spatial domain $\Omega$. Each mesh $\mathcal{T}_h$ consists of closed triangles such that ${\overline{\Omega} = \bigcup_{K \in \mathcal{T}_h } \overline{K}}$  holds. We denote the finite dimensional approximations of the polymer, non-fullerene acceptor, solvent, and air volume fraction fields by ${\phi_i^h\in \mathcal{S}^h}$, the finite dimensional approximations for their corresponding chemical potential fields by ${\mu_i^h\in \mathcal{M}^h}$. Additionally, we define the approximation $\Phi_v^h \in \mathcal{S}^h$ for the vapor fraction field, as well as the approximations $\bv^h \in  \mathcal{V}^h$ and $p^h \in \mathcal{P}^h$ for the velocity and pressure fields, respectively. The discretized function spaces are chosen as:
\begin{subequations}
    \begin{align}
        \mathcal{S}^h &= \sett{ \phi^h \in C(\overline{\Omega})}{ \phi^h \mid_{K}\in \mathbb{P}^\ell(K) \quad \forall K \in \mathcal{T}_h}, \label{spa:vol}\\
        \mathcal{M}^h &= \sett{ \mu^h \in C(\overline{\Omega})}{ \mu^h \mid_{K}\in \mathbb{P}^\ell(K) \quad \forall K \in \mathcal{T}_h}, \label{spa:che}\\
        \mathcal{V}^h &= \sett{ \bv^h \in C(\overline{\Omega})^n}{\bv^h \mid_{K}\in (\mathbb{P}^{\ell +1}(K))^n \quad \forall K \in \mathcal{T}_h, \quad \bv \mid_{\partial \Omega} = 0}, \\
        \mathcal{P}^h &= \sett{ p^h \in C(\overline{\Omega})}{ p^h \mid_{K}\in \mathbb{P}^\ell(K) \quad \forall K \in \mathcal{T}_h}, 
    \end{align}
\end{subequations}
with  $\mathbb{P}^\ell(K)$ being the set of all polynomials on $K$ of degree at most $\ell$. The spaces \eqref{spa:vol} and \eqref{spa:che} are identical, but they are provided separately. As mentioned earlier, we are using the splitting form from \cite{Wodo2012_height} to circumvent the limitations associated with the continuity of the basis functions. It is noted that we use Taylor-Hood elements \cite{CTaylor_PHood_1973} for the discretization of the Navier-Stokes equation as for the stability of the numerical discretization velocity ($\ell +1=2$) and pressure ($\ell=1$) fields are chosen from different spaces.

To obtain the fully discrete form of the governing equations~\eqref{eqn:model}, we decide on a semi-implicit scheme \cite{Bergermann2023, bosch2015preconditioning} where the linear part of the right-hand side of \eqref{eqn:weak} is treated implicitly and the nonlinear terms coming from the potential are handled explicitly. On the left-hand side of the Navier-Stokes equation \eqref{eqn:weakns}, the nonlinear terms are taken from the previous iteration. By $ { t_0 < t_1 < \ldots < t_m =t_{\mathrm{max}}}$, we denote an equidistant subdivision of $ (0,t_{\mathrm{max}}]$ with fixed time step size $\tau$. The existence and uniqueness of the solution of the fully discrete system obtained using the semi-implicit time scheme for the coupled Navier-Stokes and single Cahn-Hilliard equations is discussed in \cite{HGarke_MHinze_CKahle_2016, JBosch_CKahle_MStoll_2018}. 

\section{Iterative solvers and preconditioning} 
\label{sec::precon}
In the previous section, the discretization scheme for the system of PDEs was derived and as a result we obtained several large-scale linear systems. In general, without the choice of a preconditioner or non-Eucledian inner product the conditioning of the stiffness matrix scales with $h^d$ where $h$ is the mesh parameter \cite{elman2014finite} and as such we require an iterative scheme that performs robust with changes to $h$ and in turn with changes to the degrees of freedom of the finite element space. For this, we use Krylov subspace methods \cite{elman2014finite,saad2003iterative} and as a general rule we will aim to sufficiently cluster the eigenvalues to achieve a robust convergence behavior with respect to changes in the system parameters such as the mesh size, the time step or the model parameters of Cahn-Hilliard, Allen-Cahn, or Navier-Stokes equations. 

We start by discussing the preconditioning of the Cahn-Hilliard equations discretized by the finite element scheme discussed in \Cref{sec::discretization} that lead to the following systems
\begin{subequations}
\begin{align}
    \frac{1}{\tau} \mathbf{M} \left(\boldsymbol{\phi}_i^{k+1} - \boldsymbol{\phi}_i^{k} \right) + \mathbf{C}\boldsymbol{\phi}_i^{k+1} &= - \gamma_i \mathbf{K} \boldsymbol{\mu}_i^{k+1}, \nonumber \\ 
    \mathbf{M} \boldsymbol{\mu}_i^{k+1} &=  \boldsymbol{f}_i^{k} + \beta_i \mathbf{K} \boldsymbol{\phi}_i^{k+1}, \label{eqn:airdisc}
      \end{align}
\end{subequations}
where $\boldsymbol{\phi}_i^k$ 
and $ \boldsymbol{\mu}_{i}^k$ 
are the coefficients vectors representing the discretized order parameters and chemical potentials at time $t_k$ for each of the individual materials, respectively. The matrices $\mathbf{M}$, $\mathbf{C}$, and $\mathbf{K}$ denote the standard mass, convection and stiffness matrices arising  from the spatial discretization.
The right-hand side vectors $\boldsymbol{f}_i^{k}$ represent the discretized representations of terms containing $\frac{\partial f}{\partial \phi_i}$, along with the added boundary contributions for the different phases indicated by $i \in \{p, \nfa, s, a \}$. Note that they are the same for polymer, the non-fullerene acceptor, the solvent, and the air equations. We now write the above equations in the notation of a saddle-point system \cite{elman2014finite,benzi2005numerical} as this will help us in deriving suitable preconditioners. After some reordering and rescaling, we obtain the following system
\begin{equation}\label{eqn:CH_linear}
            \begin{bmatrix}
            \mathbf{M} + \tau\mathbf{C} &  \tau\gamma_s \mathbf{K}  \\
            \tau\gamma_s\mathbf{K} & - \frac{\tau\gamma_s}{\beta_s}\mathbf{M} 
        \end{bmatrix}  
        \begin{bmatrix}
          \boldsymbol{\phi}_s^{k+1}  \\
          \boldsymbol{\mu}_s^{k+1}
        \end{bmatrix} = 
        \begin{bmatrix}
           \mathbf{M} \boldsymbol{\phi}_s^{k} \\ 
           -\frac{\tau\gamma_s}{\beta_s}\boldsymbol{f}_s^{k}
        \end{bmatrix},  
\end{equation}
which represents the equation for the solvent. Analogously, we can write the system for NFA, polymer and air. We arrive at a saddle-point system, which is nonsymmetric only because of the convection term represented by the matrix $\mathbf{C}$. The design of our preconditioner is guided by the idealized block preconditioner
    \begin{align}
       P= \begin{bmatrix}
            \mathbf{M} + \tau\mathbf{C} & 0 \\
             0 & \mathbf{S}
        \end{bmatrix} ,        
    \end{align} 
where $ \mathbf{S} =\frac{\tau\gamma_s}{\beta_s} \mathbf{M} + \left( \tau\gamma_s \right)^2 \mathbf{K} \left( \mathbf{M} + \tau \mathbf{C}  \right)^{-1} \mathbf{K} $  is the (negative) Schur complement. Block-diagonal preconditioners have the advantage of being cheap to apply once good approximations to the inverse of the diagonal blocks are available  
\cite{bosch2014fast,bosch2015preconditioning,murphy2000note}. For this, we now discuss a practical version of the idealized preconditioner as
\begin{align}
       \widetilde{P} = \begin{bmatrix}
            \mathbf{L} & 0 \\
             0 & \widetilde{\mathbf{S}}
        \end{bmatrix} ,        
\end{align} 
where $\mathbf{L}\approx \mathbf{M}+\tau \mathbf{C}$ is an approximation of the $(1,1)$-block of the saddle-point system. For this approximation, we rely on the use of an algebraic multigrid scheme \cite{ruge1987algebraic,xu2017algebraic,falgout2006introduction}. 
The more difficult term to approximate is the Schur complement $S$ for which we proceed using the following idea. We first assume that $\mathbf{C}=0$ to get the modified Schur complement as
\begin{equation}
\label{eqn::schurCH}
\frac{\tau\gamma_s}{\beta_s}\mathbf{M} + \left( \tau\gamma_s \right)^2 \mathbf{K} \mathbf{M}^{-1} \mathbf{K}.
\end{equation}
The direct solution to this system is still difficult given the sum of large-scale matrices. For this, we follow a matching approach \cite{pearson2012regularization} where we approximate the sum by a product of matrices whose inverses are cheap to approximate. The resulting approximation is given by 
$$ \left(\tau \gamma_s \mathbf{K} + \sqrt{\frac{\tau\gamma_s}{\beta_s}}\mathbf{M}\right) \mathbf{M}^{-1} \left(\tau \gamma_s \mathbf{K} + \sqrt{\frac{\tau\gamma_s}{\beta_s}}\mathbf{M}\right),$$
where we point out that this product matches the two terms in the modified Schur complement. For the symmetric case $\mathbf{C}=0$, it has been shown that such approximations lead to a provably optimal performance \cite{Bergermann2023}. The use of block preconditioners for the Cahn-Hilliard equation has also been used in \cite{boyanova2012efficient} and extensions to vector-valued Cahn-Hilliard equations are applied in \cite{bosch2015preconditioning,boyanova2014efficient}. For practical purposes, we again emphasize that we are interested in approximating the inverse of \Cref{eqn::schurCH} for which we need approximations for the inverse of a weighted sum of a mass and stiffness matrix. As all coefficients in this sum are positive, we can use the algebraic multigrid approximation  of $\tau \gamma_s \mathbf{K} + \sqrt{\frac{\tau\gamma_s}{\beta_s^2}}\mathbf{M}$ using \texttt{PYAMG}. Hence, the approximation of the Schur complement requires applying the algebraic multigrid scheme twice and one additionally matrix vector product with the sparse mass matrix $\mathbf{M}.$ 

We now discuss the approximation of the Allen-Cahn part written in discretized form as
\begin{align}
    \frac{1}{\tau} \mathbf{M} \left(\boldsymbol{\Phi}_v^{k+1} - \boldsymbol{\Phi}_v^{k} \right) + \mathbf{C}\boldsymbol{\Phi}_v^{k+1}   + \delta_v \mathbf{K} \boldsymbol{\Phi}_v^{k+1}  &= - \gamma_v\boldsymbol{f}_v^{k},
\end{align}
and as the solution of one linear system we obtain
\begin{align}
 \left( \mathbf{M} + \tau \delta_v \mathbf{K} +\tau \mathbf{C}\right) \boldsymbol{\Phi}_v^{k+1}  &= \mathbf{M}\boldsymbol{\Phi}_v^{k} - \tau \gamma_v\boldsymbol{f}_v^{k}.
\end{align}
As this equation is only of second order the system is not coupled and a straightforward approximation of the matrix $\mathbf{M} + \tau \delta_v \mathbf{K} +\tau \mathbf{C}$ using \texttt{PYAMG} is our approximation of choice. Again, this is a nonsymmetric system due to the convection term $\mathbf{C}$ and we apply the preconditioner within a nonsymmetric Krylov subspace solver.

It remains to discuss the efficient solution of the discretized Navier-Stokes equations written in the following form 
\begin{align}\label{eqn:matrixns}
    \frac{1}{\tau} \rho^k  \mathbf{M} \left( \boldsymbol{v}^{k+1} - \boldsymbol{v}^{k} \right) + \mathbf{J} \boldsymbol{v}^{k+1} + \mathbf{K_{sym}}\boldsymbol{v}^{k+1} + \mathbf{B}^\top\boldsymbol{p}^{k+1}
      &= \mathbf{G}^{k} ,\\
    -\mathbf{B} \boldsymbol{v}^{k+1}  &= 0,
\end{align}
where $\mathbf{J} = (\mathbf{j}_{mn}) $, $\mathbf{K_{sym}} = (\mathbf{k}_{mn}) $, and $\mathbf{B} = (\mathbf{b}_{ij})$ with the individual terms given by $ \mathbf{j}_{mn} = \left( (\rho^k \bv^k - \sum_{i}\alpha_i\nabla \mu_i^k )\cdot \nabla \zeta^n, \zeta^m \right)$,   $\mathbf{k}_{mn} = \left(2\eta \sym \zeta^n ,\sym \zeta^m  \right)$, and $\mathbf{b}_{mn} = -\left( \dv \, \zeta^m , \hat{\zeta}^n\right)$  for $\{\zeta^m\} \subset \mathcal{V}^h$ and $\{\hat{\zeta}^m\} \subset \mathcal{P}^h$. The right-hand side vector $\mathbf{G}^{k}$ is the discretized representation of the terms involving $\mathcal{G}$. After reordering and rescaling the above equations, during each time step we need to solve a linear system of the following structure
\begin{align}\label{eqn:linearnav}
        \begin{bmatrix}
            \mathbf{A}  &   \mathbf{B}^\top  \\
             \mathbf{B}  & 0
        \end{bmatrix}  
        \begin{bmatrix}
          \boldsymbol{v}^{k+1}  \\
          \boldsymbol{p}^{k+1}
        \end{bmatrix} = 
        \begin{bmatrix}
         \frac{\rho^k}{\tau} \mathbf{M}\boldsymbol{v}^{k} +  \mathbf{F}_v^{k} \\
           0
        \end{bmatrix},       
\end{align}  
with $ \mathbf{A} = \frac{\rho^k}{\tau} \mathbf{M +  \mathbf{J}} +  \mathbf{K}_{\epsilon}$. 
The idealized preconditioner is again of diagonal form 
$$
\mathbf{P}=
\begin{bmatrix}
    \hat{\mathbf{A}}&\\
    &\hat{\mathbf{S}}
\end{bmatrix},
$$
where $\hat{\mathbf{A}}\approx \mathbf{A}.$ The matrix $\mathbf{A}$ corresponds to a vector-valued convection diffusion problem for which we can use an algebraic multigrid approximation. Note that if a Newton scheme is employed the matrix $\mathbf{A}$ would change in every Newton iteration. The more intricate part comes from the Schur complement of the system $\mathbf{B}\mathbf{A}^{-1}\mathbf{B}^\top$. Following the derivations in \cite{elman2014finite}, we suggest the use of a pressure-convection-diffusion preconditioner that is based on the commutator of the Navier-Stokes system. The main idea follows from assuming that the commutator here only illustrated for the stationary equations
$$
\mathcal{E}=(-\Delta +w\cdot\nabla)\nabla-\nabla (-\Delta +w\cdot\nabla)_p,
$$
where the index $p$ means the operator defined on the pressure space. The assumption that the commutator is close to zero and the discretization of the corresponding operators leads to the approximation
$$
\mathbf{B}\mathbf{A}^{-1}\mathbf{B}^\top\approx \mathbf{K}_p\mathbf{A}_p^{-1}\mathbf{M}_p,
$$
where $\mathbf{K}_p$, $\mathbf{A}_p,$ and $\mathbf{M}_p$ are the pressure stiffness matrix, the pressure diffusion convection matrix, and the pressure mass matrix, respectively. Note that $\mathbf{A}_p$ represents the discretization of the bilinear form $\mathbf{A}$ defined on the pressure space, which include the contribution from diffusion, convection and the time-discretization via the Euler method. To further illustrate the efficiency of this approach, we point out that for the application of the preconditioner we require the application of 
$$
\mathbf{M}_p^{-1}\mathbf{A}_p\mathbf{K}_p^{-1},
$$
where the two-inverse for stiffness and mass matrix can be cheaply approximated using a multigrid scheme or other approximation as both matrices are symmetric and positive (semi-)definite.

\section{Numerical results}
\label{sec::results}
In this section, we present numerical results on several setups of the discussed models to examine the efficiency of the proposed numerical approaches. The model discussed in Section~\ref{sec::precon} is solved with a \texttt{PYTHON} implementation using the finite element libraries \texttt{DOLFINX} \cite{dolfinx}
, \texttt{Basix} \cite{basix}
, and \texttt{UFL}  \cite{ufl} 
from the \texttt{FENICS} project  \cite{fenics, autofenics} 
with the version 0.8.0. The codes in \cite{Bergermann2023}, which are publicly available\footnote[1]{\url{https://github.com/KBergermann/Precond-Cahn-Hilliard-OSC}}, form the basis for the treatment of the Cahn-Hilliard model. All our experiments are performed on an Ubuntu Linux machine having 13th Gen Intel i7-1355U 12-Core processor with 64 GB RAM. 

For better numerical stability, the free energy density change was substituted by a polynomial approximation. All terms of the form $\phi \ln \phi$ are approximated by their respective Taylor series expansion around $\phi=1/2$ up to 4th order. Furthermore, to guarantee the requirement $\phi_p + \phi_{\nfa} + \phi_s + \phi_a = 1$, we substitute the air volume fraction by $\phi_a = 1 - \phi_p - \phi_{\nfa} - \phi_s$, and reduce the linear system by one variable. However, note that the chemical potential field $\mu_a$, appearing in the Navier-Stokes equation \eqref{eqn:matrixns}, still needs to be computed. This is simply done by setting $\phi_a^{k+1} = 1 - \phi_p^{k+1} - \phi_{\nfa}^{k+1} - \phi_s^{k+1}$ after finding the iterates $\phi_p^{k+1}, \phi_{\nfa}^{k+1}$, and $ \phi_s^{k+1}$, and solve the chemical potential field equation \eqref{eqn:airdisc}. 

\begin{figure}[h!]
    \centering
    \begin{tikzpicture}[scale=1]
    \draw[thick, name path = r] (0,0) rectangle (10,2);
    \path[name path = h] (0,0) -- (10,0);
    \draw[snake=snake, segment length = 100, segment amplitude = 1, name path = s] (0,1) -- (10,1);                          
    \node[scale=0.75, text width=9cm]  at (5.5,1.4) {$ \phi_p = \phi_{\nfa} = \phi_s = 0 $, \qquad  \;\qquad  \qquad $ \phi_a = \Phi_v = 1$};
    \node[scale=0.75, text width=11cm]  at (5,0.4) {$ \phi_p = \phi_{\nfa} = a \pm 0.01 $, \; $ \phi_s = b \pm 0.01$  \qquad  \qquad  $ \phi_a = \Phi_v = 0$};
    \tikzfillbetween[of=h and s]{cyan, opacity = 0.1};
    \end{tikzpicture}
    \caption{Initial condition in 2D}
    \label{fig_domain}
\end{figure}
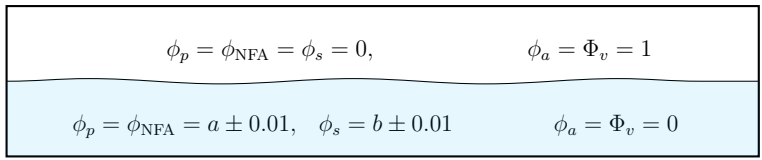

Figure~\ref{fig_domain} shows the setup of the initial conditions in the 2D-setting with $2a + b = 1$. Here, $\pm 0.01 $ denotes uniformly distributed random fluctuations. At time $t = 0$, these cause nonzero concentration gradients, which are crucial for phase separation to begin. 

In this paper, the generalized minimal residual (\texttt{GMRES}) method \cite{saad2003iterative} is applied as a nonsymmetric Krylov subspace solver and its implementation within \texttt{SCIPY}  \cite{SciPy} 
is used. As stopping criterion, the \texttt{GMRES} method uses \linebreak ${ \lVert \mathfrak{b} - \mathcal{A}z^n \rVert \leq \max \{ \texttt{rtol}  \lVert \mathfrak{b} \rVert, \texttt{atol} \} },$
where $\mathfrak{b}$ and $\mathcal{A}$ denote the corresponding right-hand side and the coefficient matrix, and  $z^n$ is the calculated solution at \texttt{GMRES} step $n$. In solving the Navier-Stokes equations, the \texttt{AMG} preconditioning for symmetric matrices and  non-symmetric matrices are implemented using the Ruge-Stüben AMG method and the approximate ideal restriction (AIR) method from the \texttt{Python} package \texttt{PYAMG} \cite{pyamg2023}
, respectively. For the Cahn-Hilliard equations, the Ruge-Stüben method is also used for the symmetric matrices, whereas the Smoothed Aggregation (SA) method is employed for non-symmetric matrices. Unless stated otherwise, we set the given tolerance values to $\texttt{rtol} = 10^{-6}$ and $\texttt{atol} = 10^{-8}$ for the solutions of Cahn-Hilliard and Allen-Cahn equations and $\texttt{rtol} = 10^{-4}$ and $\texttt{atol} = 10^{-8}$ for the solution of Navier-Stokes equation. Furthermore, the \texttt{AMG} preconditioner tolerance is set to $10^{-4}$. 

With an adjustable number of grid points $n_y$, $n_x \times n_y$, and  $n_x \times n_y \times n_z $, the uniform triangulations of $1$-, $2$- and $3$-dimensional domains of size $10$, $10 \times 10$, $10 \times 10 \times 10$ are created. Linear triangular Lagrange elements are chosen. The saturated vapor volume fractions are selected as $\phi_{sat(p)} = \phi_{sat(\nfa)} = 10^{-8}$, $\phi_{sat(s)}  = 10^{-6}$, and $\phi_{sat(a)}  = 1$.  

\subsection{Evaporation of solvent} \label{subsec:sa}
In this section, we investigate the evaporative effects of the model by simulating solely solvent and air. The Cahn-Hilliard equations involving solvent and air are considered, and the free energy functional is modified by eliminating all terms involving $\phi_p$ or $\phi_\nfa$. One and two-dimensional results are presented with initial conditions chosen as $a=0$ and $b=1$, as illustrated in Figure~\ref{fig_domain}. Given that the solution of the discretized Cahn-Hilliard equation does not satisfy the maximum principle, we use the truncated phase field function \cite{NAdam_FFranke_SAland_2020, JShen_XYang_2010}
\begin{align} \label{eqn:truncateCH}
   \hat{\phi }= \begin{cases}
         0 & \text{ for } \phi< 0, \\
         \phi & \text{ for } 0 \leq \phi \leq 1, \\
         1 & \text{ for } \phi >1
       \end{cases} 
\end{align}
to find the next iterates. All parametric details are given separately in the caption of each figure.

\begin{figure}[htp!]

\begin{minipage}[t]{0.32\textwidth}
\includegraphics[width=1.05\textwidth]{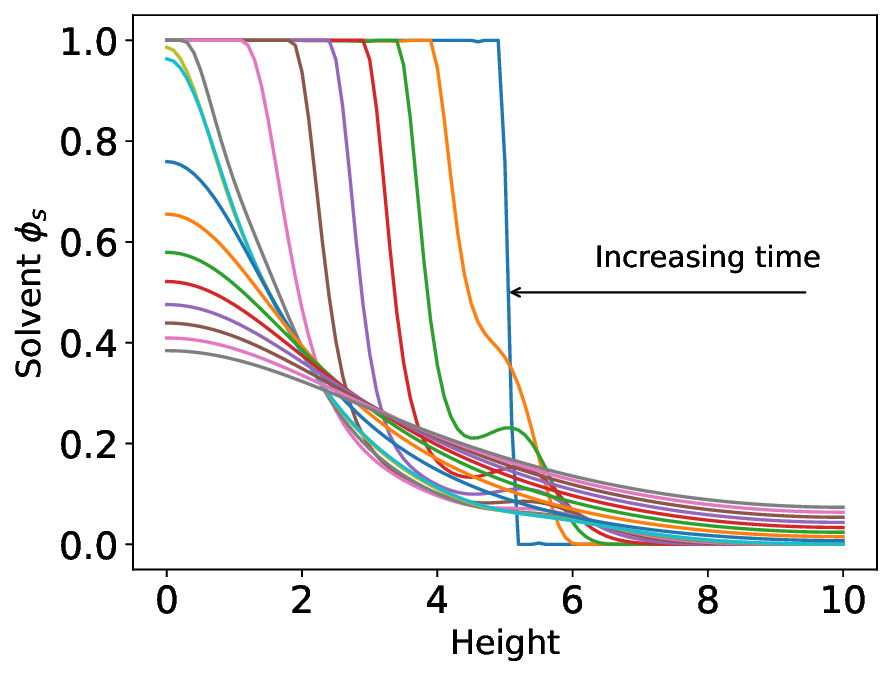}
\end{minipage}
\hfill
\begin{minipage}[t]{0.32\textwidth}
\includegraphics[width=1.05\textwidth]{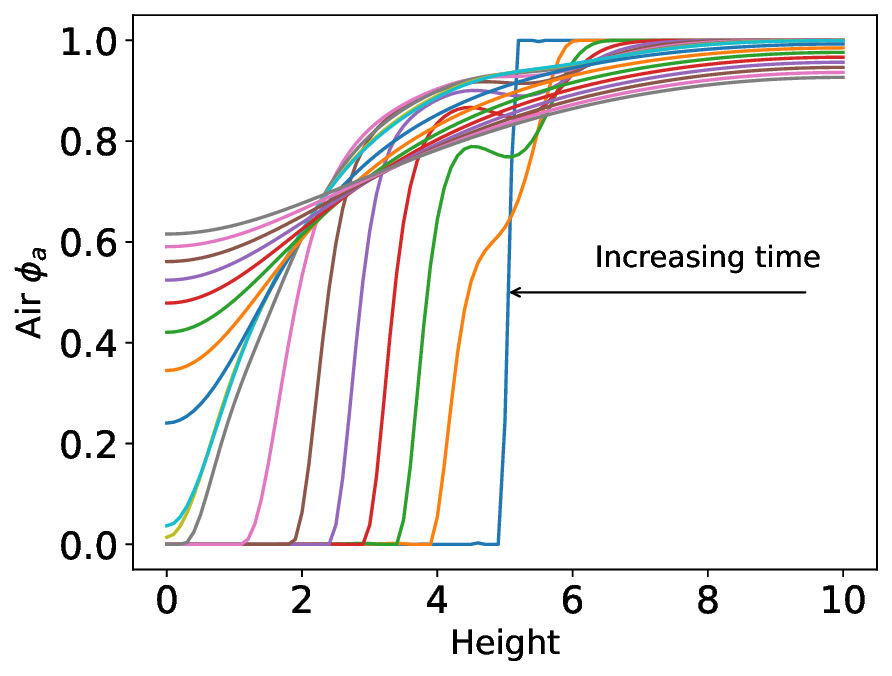}
\end{minipage}
\hfill
\begin{minipage}[t]{0.32\textwidth}
\includegraphics[width=1.05\textwidth]{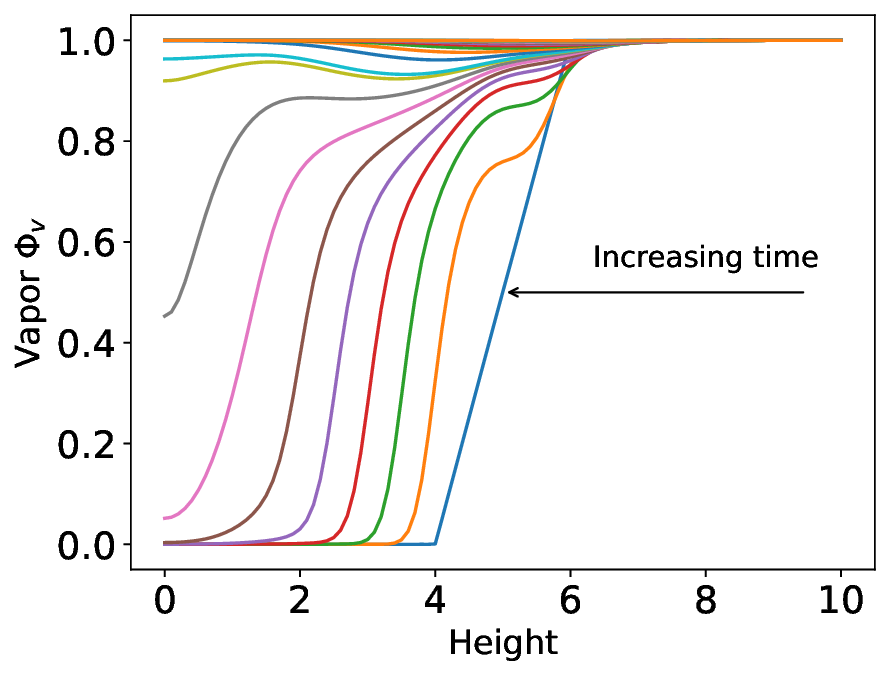}
\end{minipage}
\caption{Volume fraction fields for the solvent, air and vapor from left to right. The spatial discretization of a 1D domain is chosen $n_y =100$. Moreover, the parameters $\alpha_i = 10^{-8}$, $\beta_i = 10^{-1}$ for $ i \in \{ s, a \}$,  $\beta_v = 1$, $\gamma_i = 1$ for $ i \in \{ s, a, v \}$, and $\delta_v = 1$, the molar size of the fluid $N_s=N_s =1$, the final time $t_{\mathrm{max}} = 2.5$, the Flory–Huggins interaction parameters $\chi_{s,a} = 0$, and $\tau = 10^{-4}$ are used. }
\label{fig:sa_vol_t}
\end{figure}

Figure~\ref{fig:sa_vol_t} illustrates the phase fields for solvent, air, and vapor as functions of domain height within the 1D setting. Each subfigure represents the evolution of one variable over time, as indicated by the arrows. The solvent begins to evaporate, as we expect, and keeps doing so until it is almost completely eliminated from the domain. Additional results for the analogous 2D setting are presented in \ref{sec:appendix}. Observe that the results agree with the behavior seen in 1D results.


\begin{figure}[htp!]
\begin{minipage}[t]{0.32\textwidth}
\includegraphics[width=1.0\textwidth]{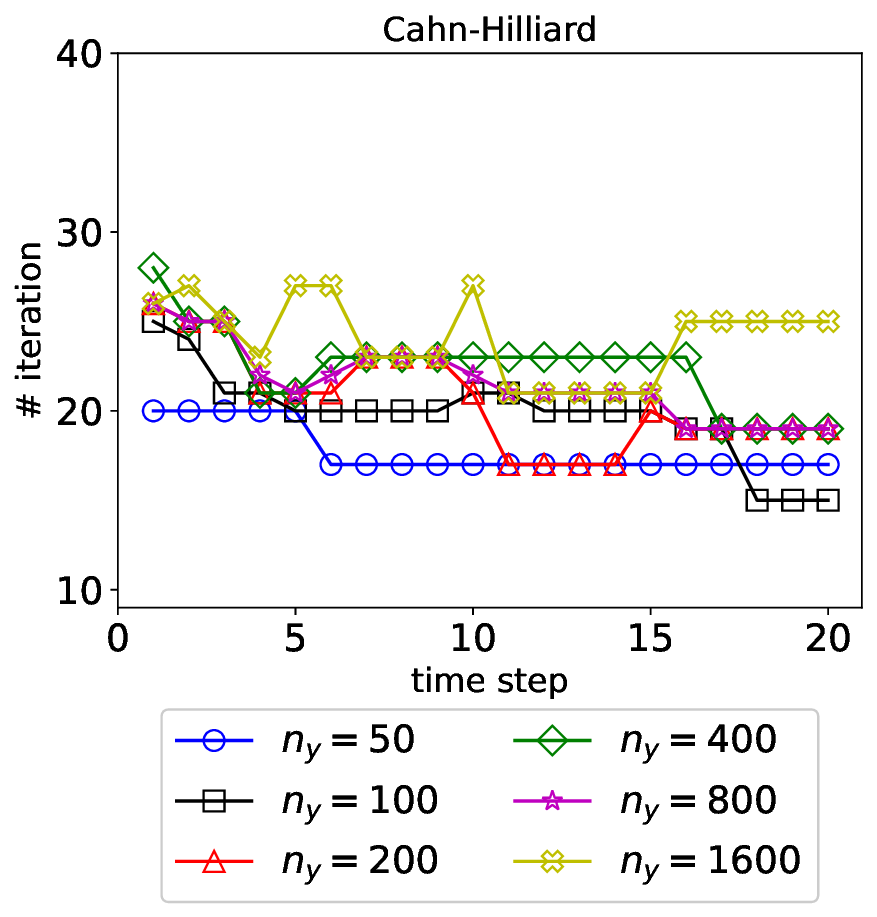}
\end{minipage}
\begin{minipage}[t]{0.32\textwidth}
\includegraphics[width=1.0\textwidth]{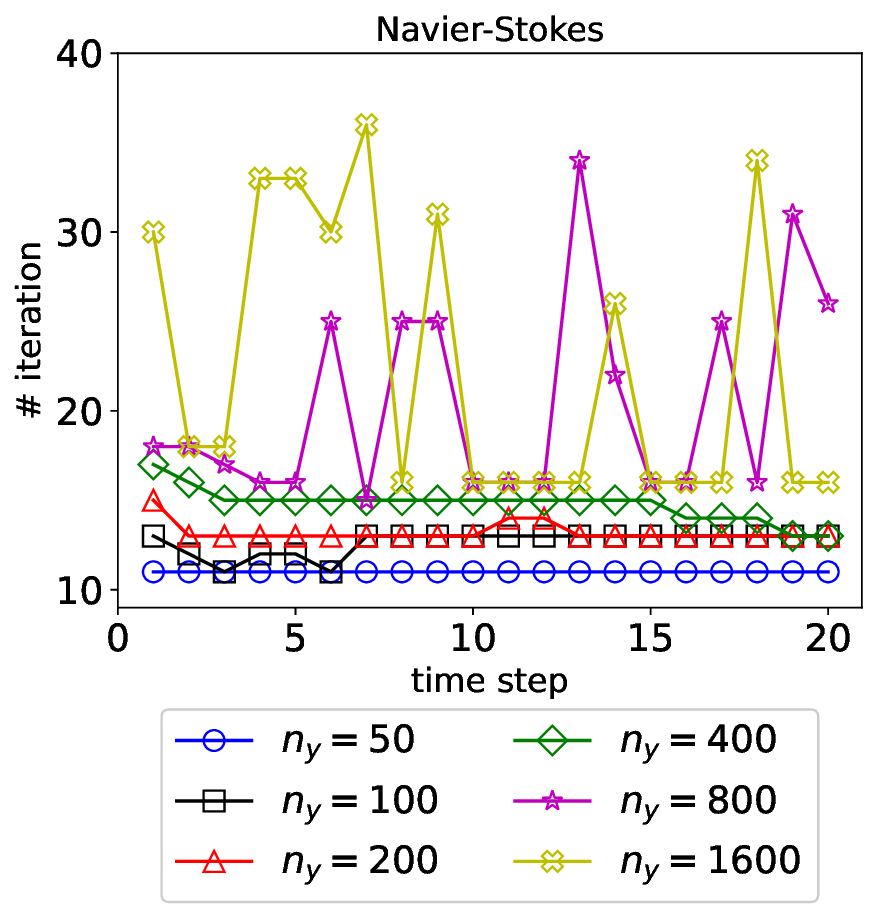}
\end{minipage}
\begin{minipage}[t]{0.32\textwidth}
\raisebox{6mm}{\includegraphics[width=1.1\textwidth]{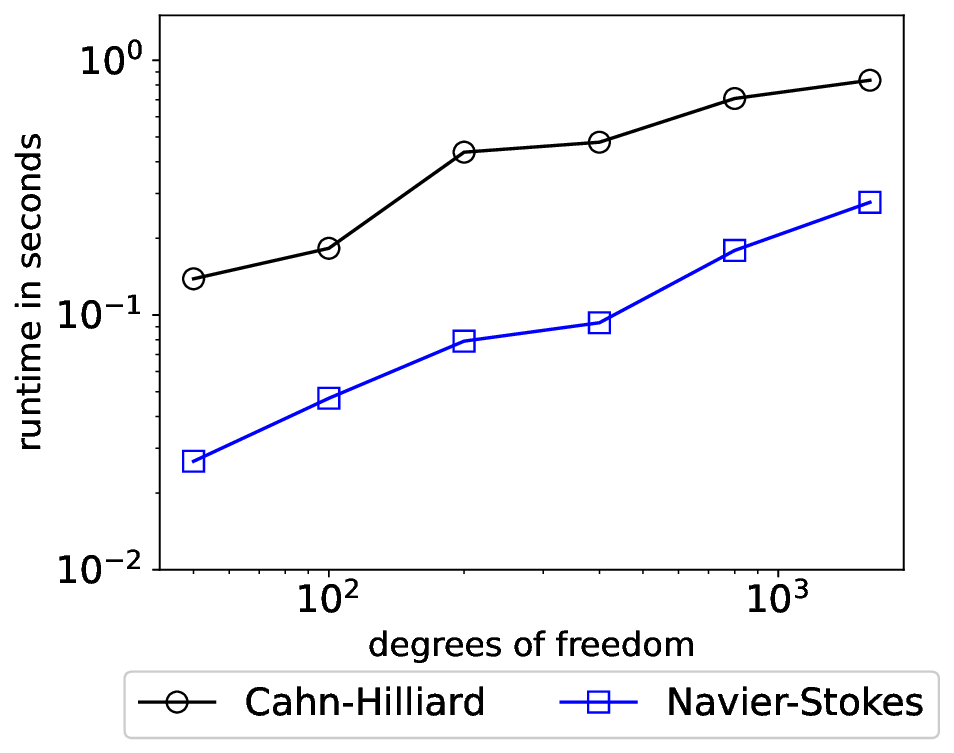}}
\end{minipage}
\caption{\texttt{GMRES} iteration numbers and runtimes (right) per time step for different spatial discretizations of a 1D domain of Cahn-Hilliard (left) and Navier-Stokes (middle).  The longest observed runtime of both equations throughout the first 20 time steps is shown. The other parameters and setup are the same as in Figure~\ref{fig:sa_vol_t}. }
\label{fig:sa_1D_gmres_disc}
\end{figure}

\begin{figure}[h!]
\begin{minipage}[t]{0.32\textwidth}
\includegraphics[width=1.0\textwidth]{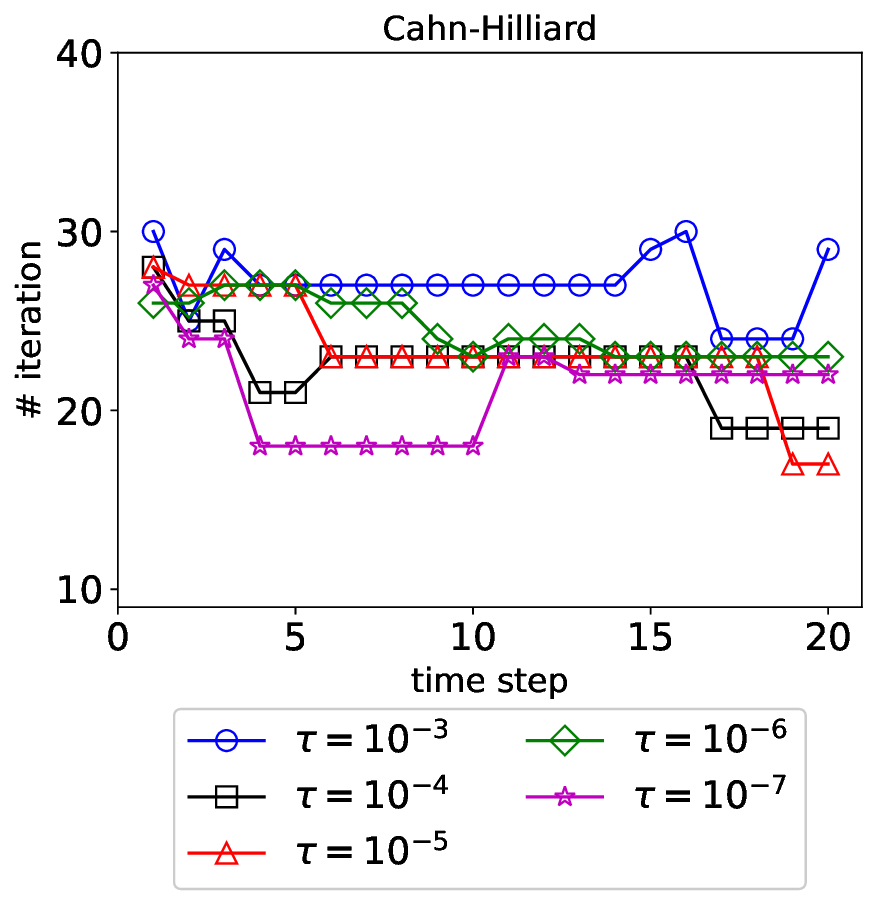}
\end{minipage}
\begin{minipage}[t]{0.32\textwidth}
\includegraphics[width=1.0\textwidth]{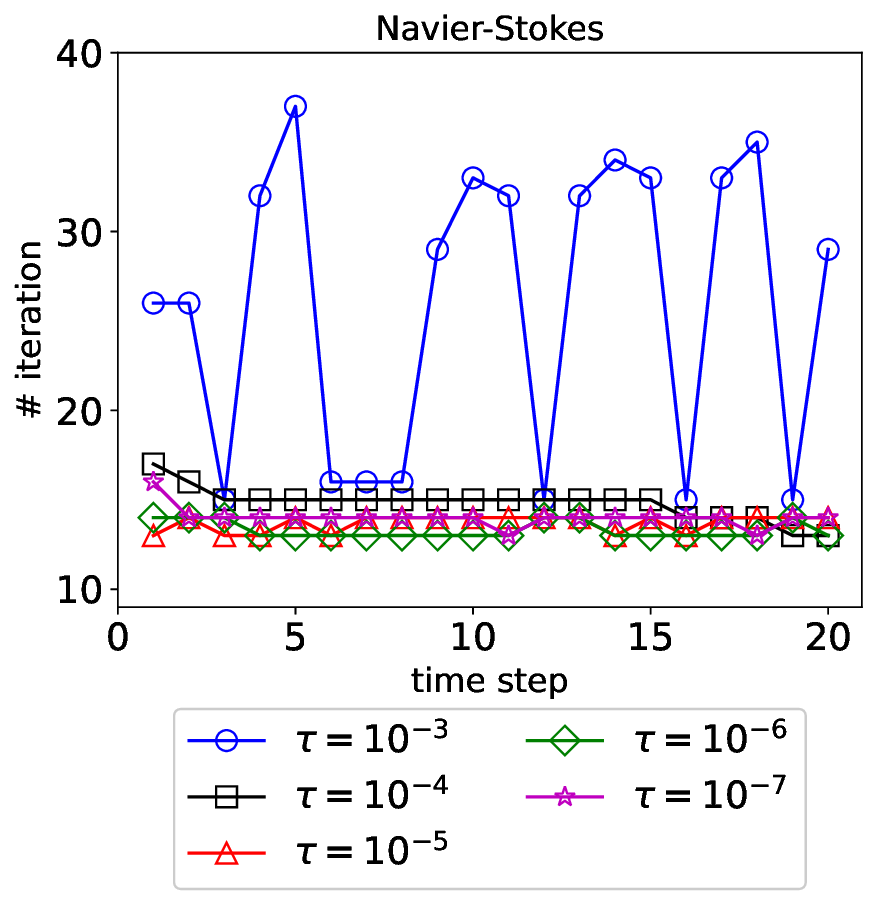}
\end{minipage}
\begin{minipage}[t]{0.32\textwidth}
\raisebox{5.0mm}{\includegraphics[width=1.1\textwidth]{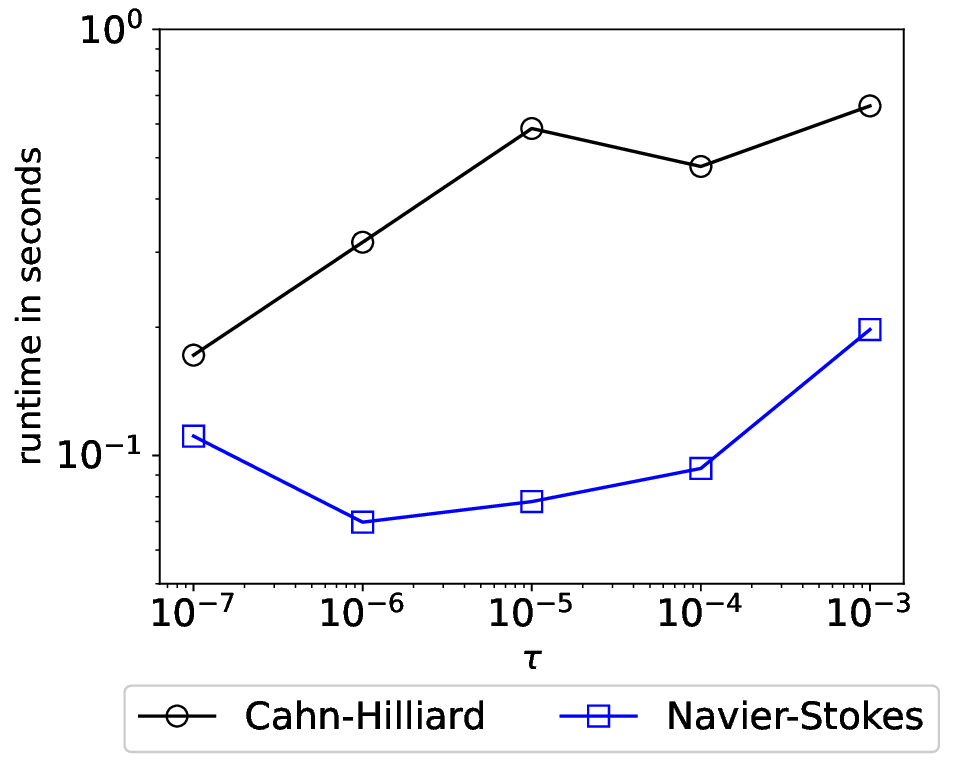}}
\end{minipage}
\caption{\texttt{GMRES} iteration numbers and runtimes (right) per time step for different time step $\tau$ of Cahn-Hilliard (left) and Navier-Stokes (middle) with $n_y = 400$ in a 1D domain. The longest observed runtime of both equations throughout the first 20 time steps is shown. The other parameters and setup are the same as in Figure~\ref{fig:sa_vol_t}. }
\label{fig:sa_1D_gmres_time}
\end{figure}

\begin{figure}[htp!]
\begin{minipage}[t]{0.32\textwidth}
\includegraphics[width=1.0\textwidth]{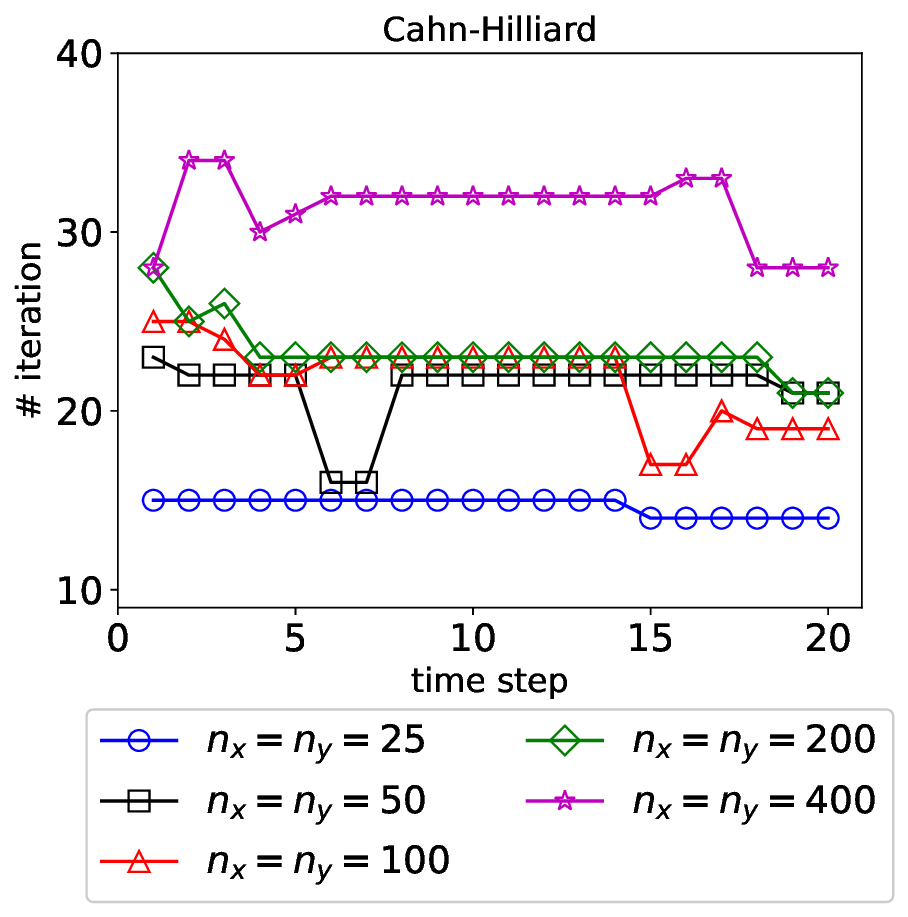}
\end{minipage}
\begin{minipage}[t]{0.32\textwidth}
\includegraphics[width=1.0\textwidth]{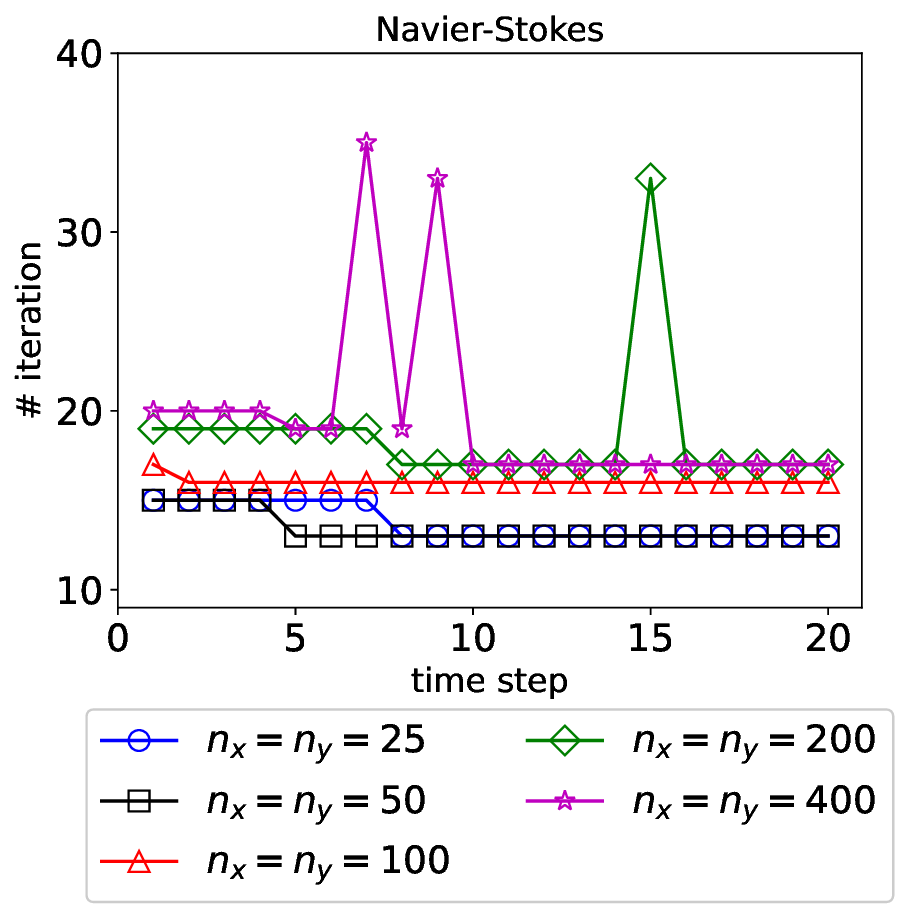}
\end{minipage}
\begin{minipage}[t]{0.32\textwidth}
\raisebox{5.8mm}{\includegraphics[width=1.1\textwidth]{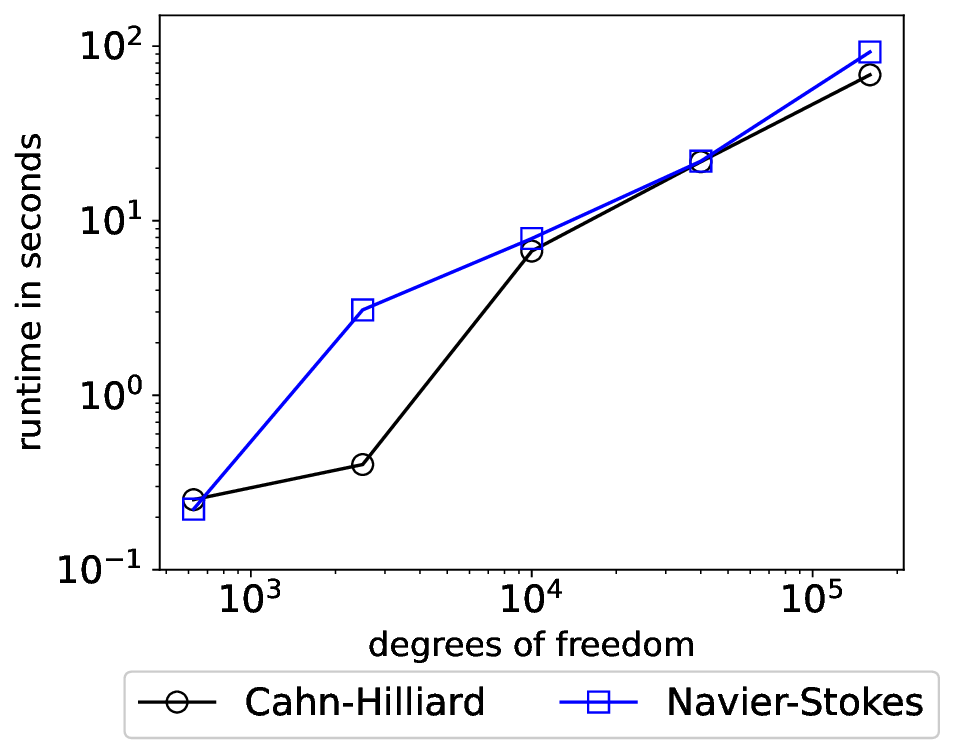}}
\end{minipage}
\caption{\texttt{GMRES} iteration numbers and runtimes (right) per time step for different spatial discretizations of a 2D domain of Cahn-Hilliard (left) and Navier-Stokes (middle).  The longest observed runtime of both equations throughout the first 20 time steps is shown.  The time step size $10^{-4}$ is used and the results shown  correspond to $t = 0$, $t =0.05$, $t = 0.15$, $t = 2.5$ (final time). The spatial discretization of the 2D domain is chosen $n_x = n_y =100$. Moreover, the parameters $\alpha_i = 10^{-8}$, $\beta_i = 10^{-1}$ for $ i \in \{ s, a \}$,  $\beta_v = 1$, $\gamma_i = 1$ for $ i \in \{ s, a, v \}$, and $\delta_v = 1$, and the molar size of the fluid $N_s=N_s =1$ are used. }
\label{fig:2D_sa_gmres_disc}
\end{figure}

\begin{figure}[h!]
\begin{minipage}[t]{0.32\textwidth}
\includegraphics[width=1.0\textwidth]{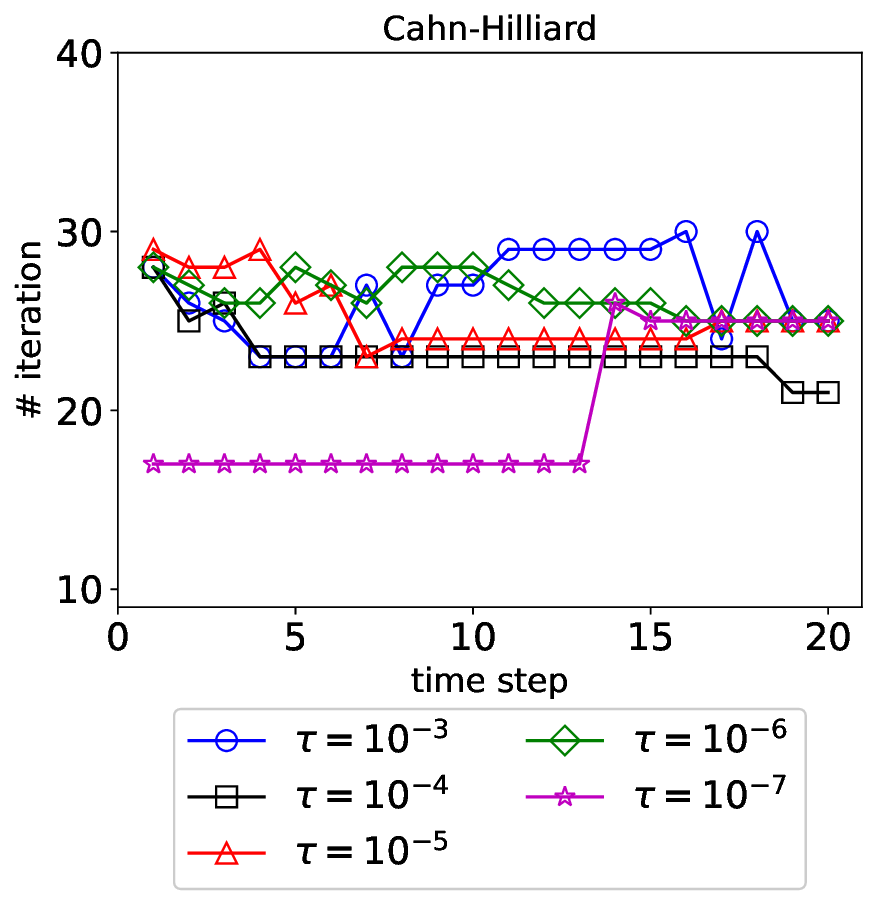}
\end{minipage}
\begin{minipage}[t]{0.32\textwidth}
\includegraphics[width=1.0\textwidth]{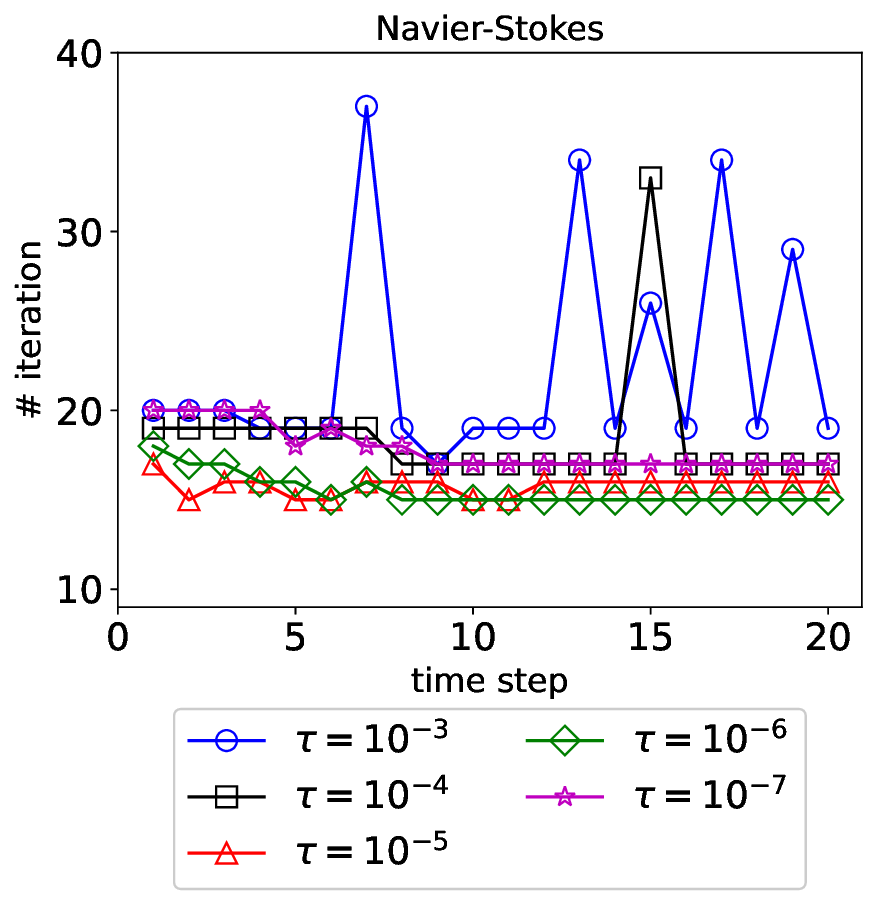}
\end{minipage}
\begin{minipage}[t]{0.32\textwidth}
\raisebox{6mm}{\includegraphics[width=1.05\textwidth]{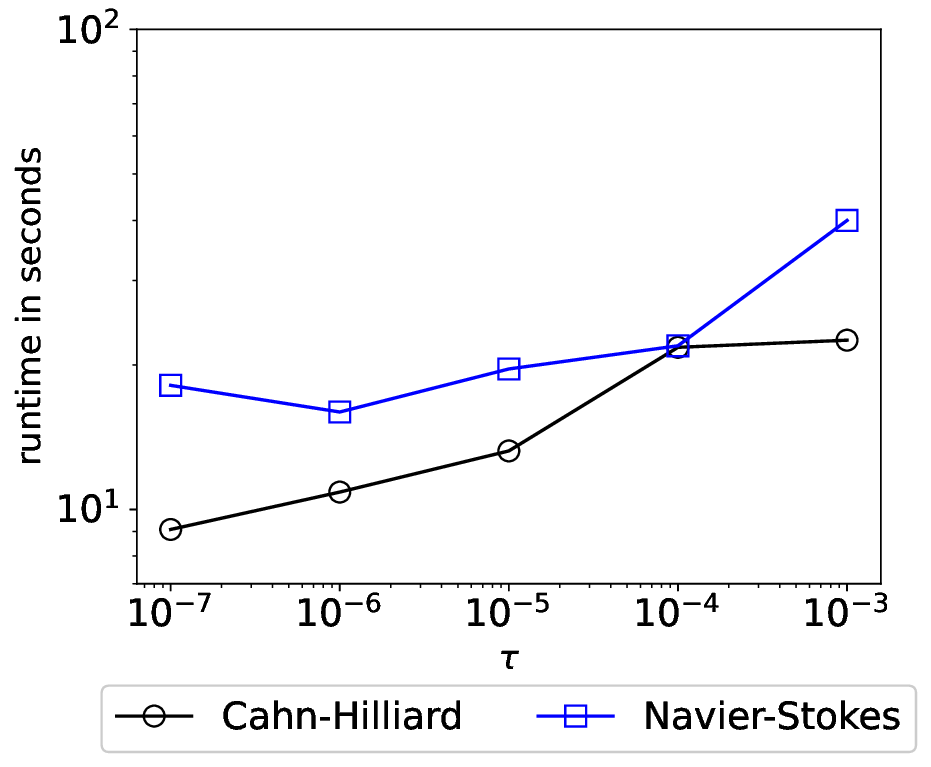}}
\end{minipage}
\caption{\texttt{GMRES} iteration numbers and runtimes (right) per time step for different time step $\tau$ of Cahn-Hilliard (left) and Navier-Stokes (middle) with $n_x=n_y = 200$ in a 2D domain. The longest observed runtime of both equations throughout the first 20 time steps is shown. The other parameters and setup are the same as in Figure~\ref{fig:2D_sa_gmres_disc}.}
\label{fig:2D_sa_gmres_t}
\end{figure}

The number of \texttt{GMRES} iterations and runtime needed to solve a single preconditioned linear system of the form \eqref{eqn:CH_linear} and \eqref{eqn:linearnav}  within a given tolerance are shown in Figures~\ref{fig:sa_1D_gmres_disc} and~\ref{fig:sa_1D_gmres_time} for the 1D domain and in Figures~\ref{fig:2D_sa_gmres_disc} and~\ref{fig:2D_sa_gmres_t} for the 2D domain. In the first $20$ time steps, the number of \texttt{GMRES} iterations stays small when varying the number of grid points, while the run time shows the linear incremental behavior for both 1D and 2D problems. Additionally, Figures~\ref{fig:sa_1D_gmres_time} and~\ref{fig:2D_sa_gmres_t} shows the stable numerical solution for the different time steps $\tau$ while the number of iterations remains small under all conditions.

\subsection{Phase separation and evaporation}
In this subsection, we expand our setting to include polymer and NFA in our system. Governed by the model equations presented in the previous sections, the aim is to provide insights into the complex morphologies and transport phenomena that arise in the production of organic solar cells. 
Both 2D and 3D results will be presented to provide a comprehensive parameter study.

\begin{figure}[h!]
\begin{minipage}[t]{0.32\textwidth}
\includegraphics[width=1.05\textwidth]{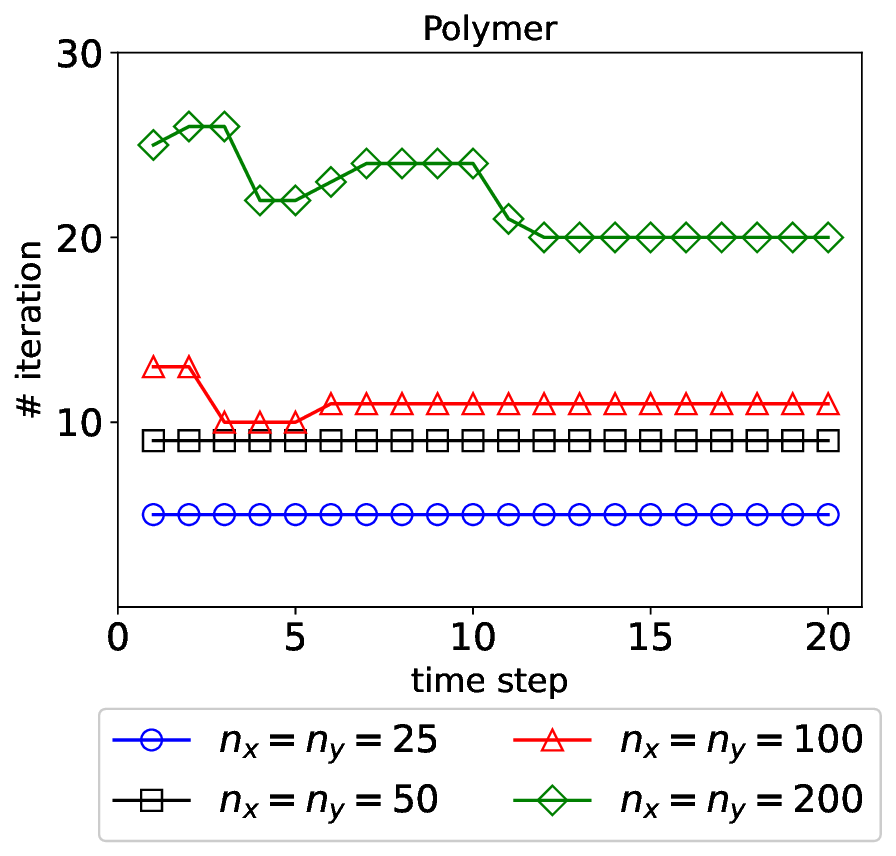}
\end{minipage}
\begin{minipage}[t]{0.32\textwidth}
\includegraphics[width=1.05\textwidth]{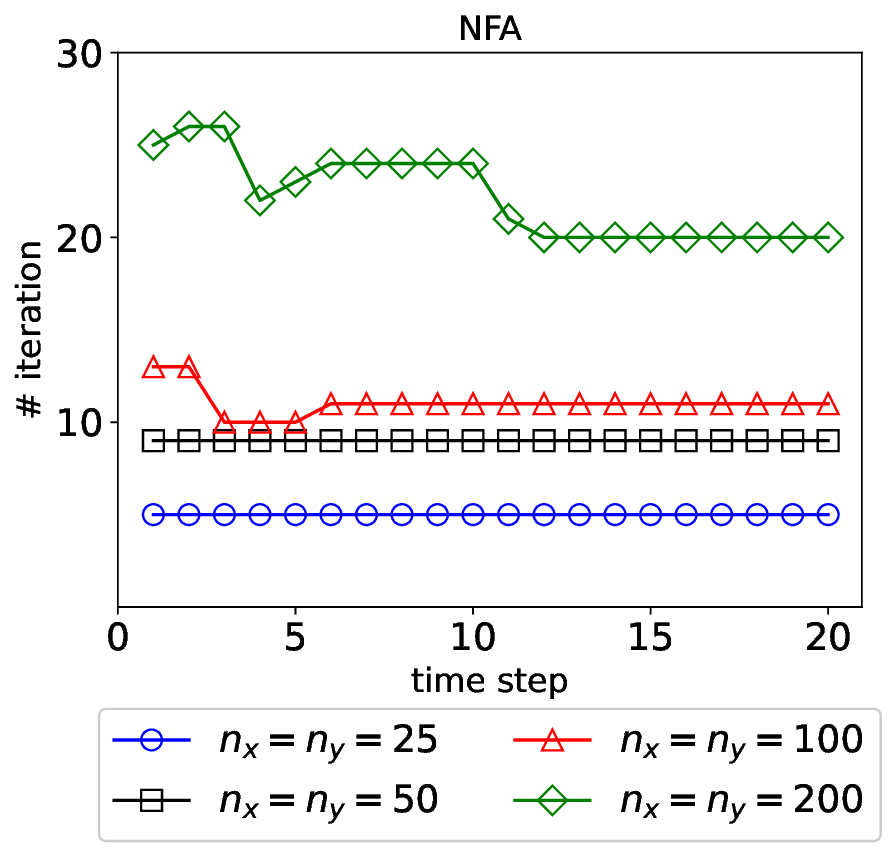}
\end{minipage}
\begin{minipage}[t]{0.32\textwidth}
\includegraphics[width=1.05\textwidth]{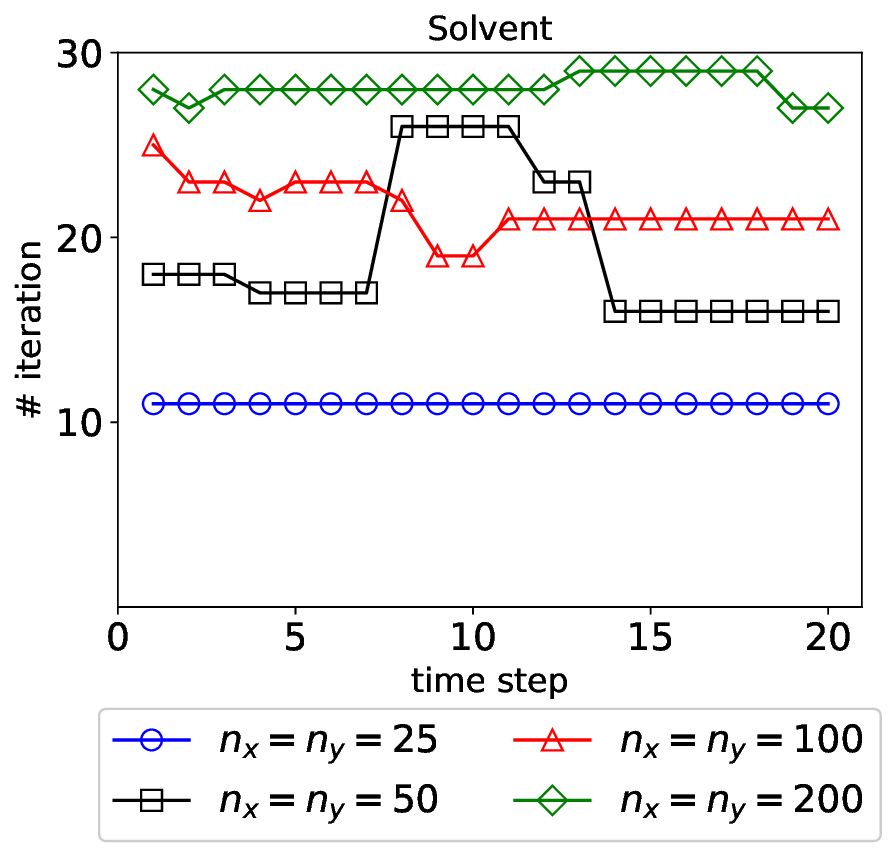}
\end{minipage}

\begin{minipage}[t]{0.32\textwidth}
\includegraphics[width=1.05\textwidth]{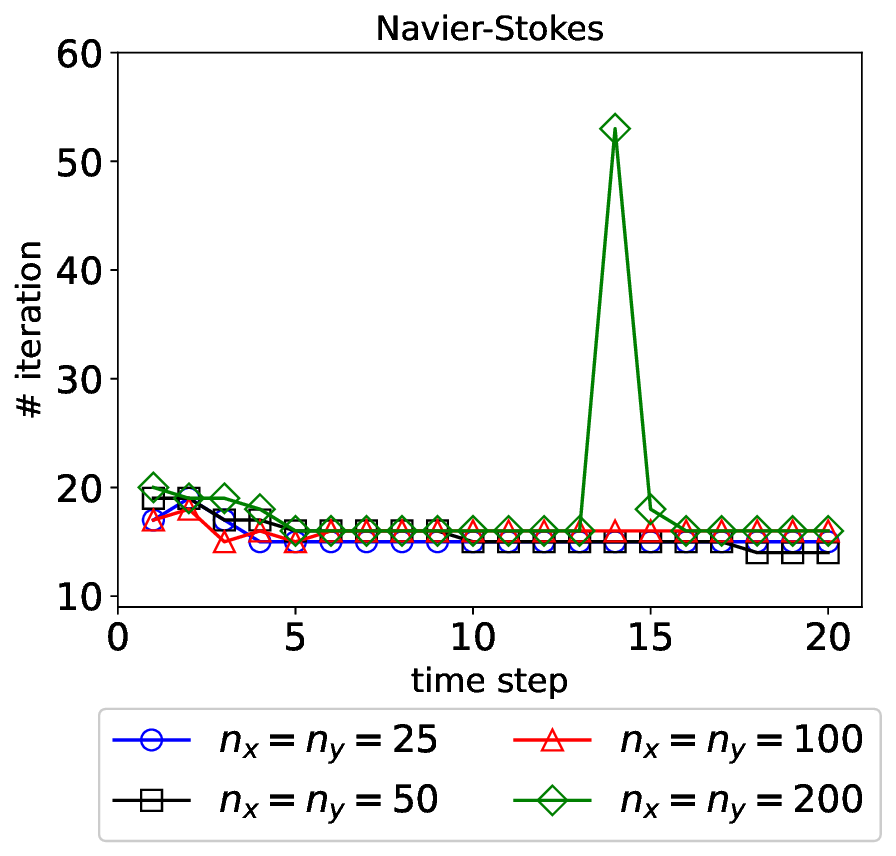}
\end{minipage}
\begin{minipage}[t]{0.32\textwidth}
\includegraphics[width=1.1\textwidth]{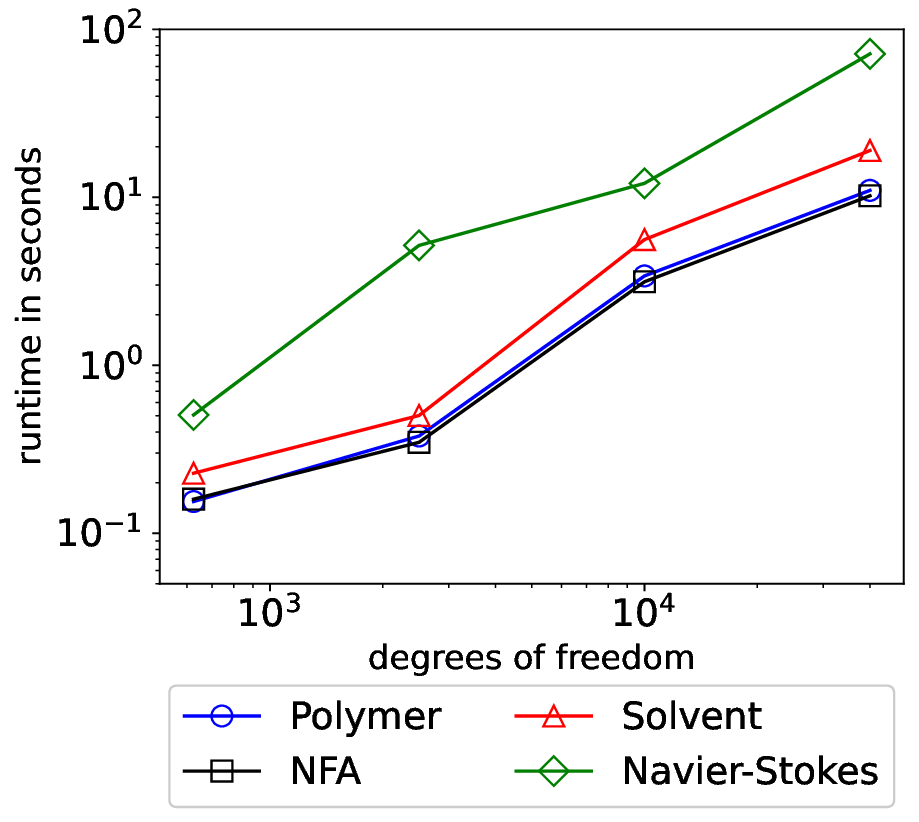}
\end{minipage}

\caption{\texttt{GMRES} iteration numbers and runtimes per time step for different spatial discretizations of Cahn-Hilliard (polymer, NFA, solvent) and Navier-Stokes with $\tau = 10^{-4}$ in a 2D domain. The spatial discretization of a 2D domain is chosen $n_x = n_y =50$. Moreover, the parameters $\alpha_i = 10^{-8}$,  $\beta_p= \beta_{\nfa} = 10^{-3}$, $\beta_s = 10^{-1}$, $\beta_a = 10^{0}$,  $\beta_v = 1$, $\gamma_i = 1$ for $ i \in \{ p, \nfa, s, a, v \}$, and $\delta_v = 1$, and the molar size of the fluid  $N_p=N_{\nfa} =5$  and $N_s=N_s =1$,  initial concentrations $a = b = 0.33$, the Flory–Huggins interaction parameters  $\chi_{p,\nfa} = 4$, $\chi_{p,s} =\chi_{\nfa,s}=0.3$,   $\chi_{p,a} = \chi_{\nfa,a} =\chi_{s,a} = 0$ are used.}
\label{fig:2D_all_gmres_iter_spatial}
\end{figure}

\begin{figure}[h!]
\begin{minipage}[t]{0.32\textwidth}
\includegraphics[width=1.05\textwidth]{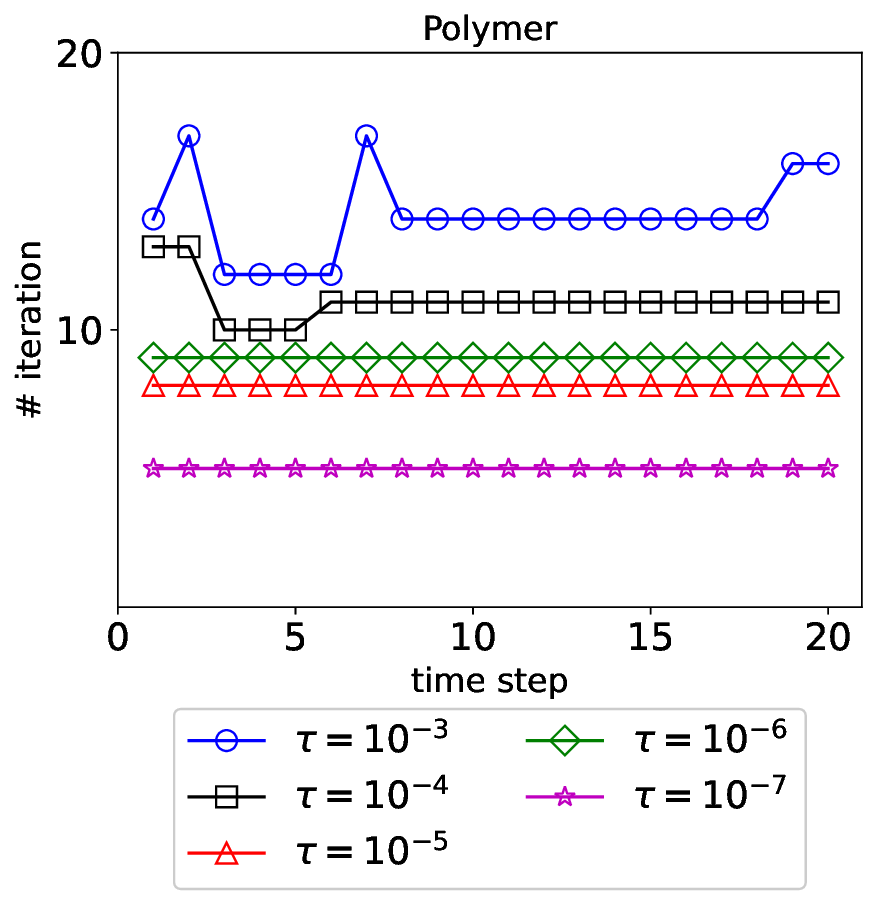}
\end{minipage}
\begin{minipage}[t]{0.32\textwidth}
\includegraphics[width=1.05\textwidth]{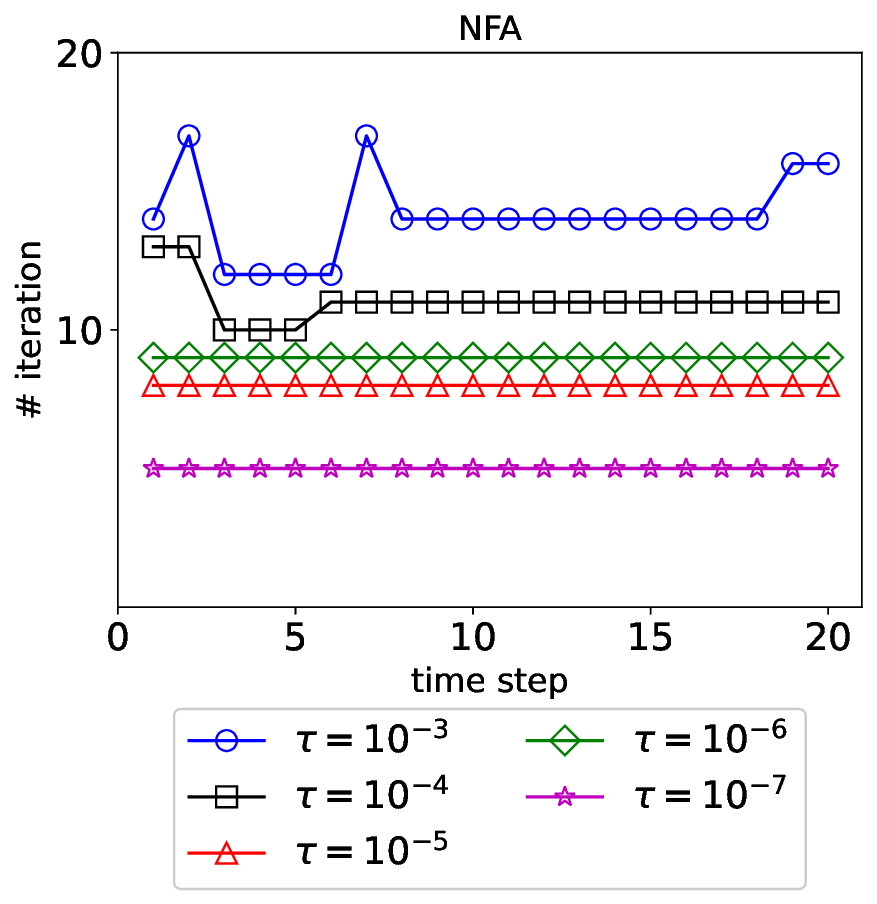}
\end{minipage}
\begin{minipage}[t]{0.32\textwidth}
\includegraphics[width=1.05\textwidth]{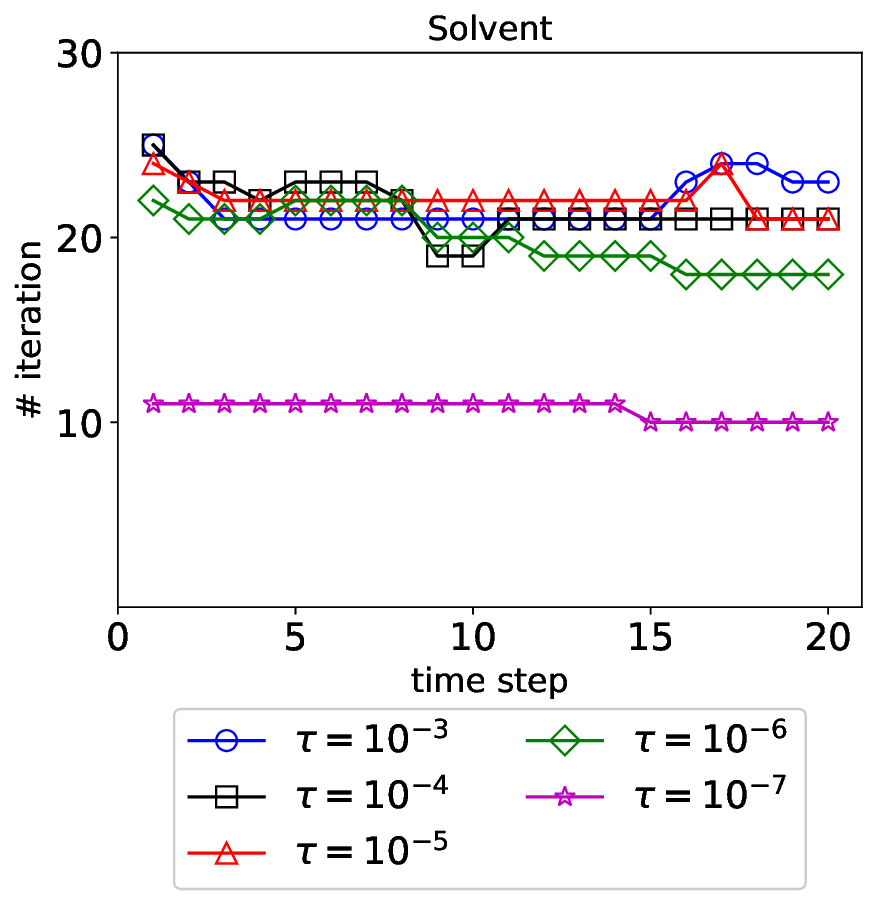}
\end{minipage}

\begin{minipage}[t]{0.32\textwidth}
\includegraphics[width=1.05\textwidth]{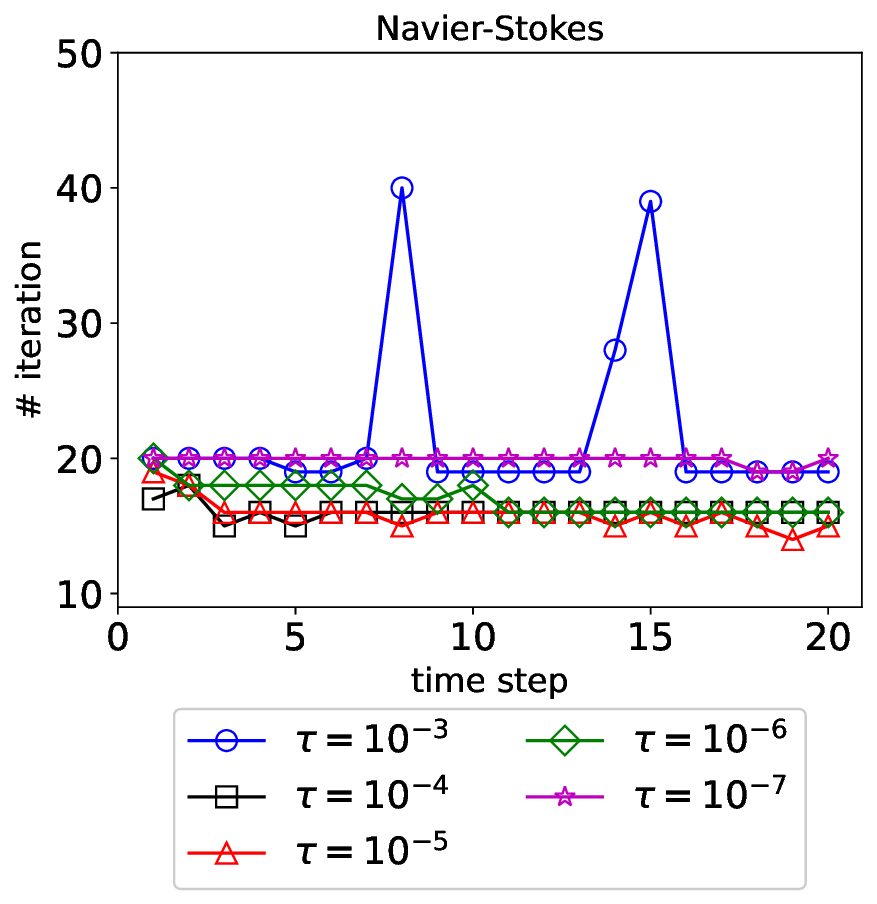}
\end{minipage}
\begin{minipage}[t]{0.32\textwidth}
\raisebox{2.5mm}{\includegraphics[width=1.1\textwidth]{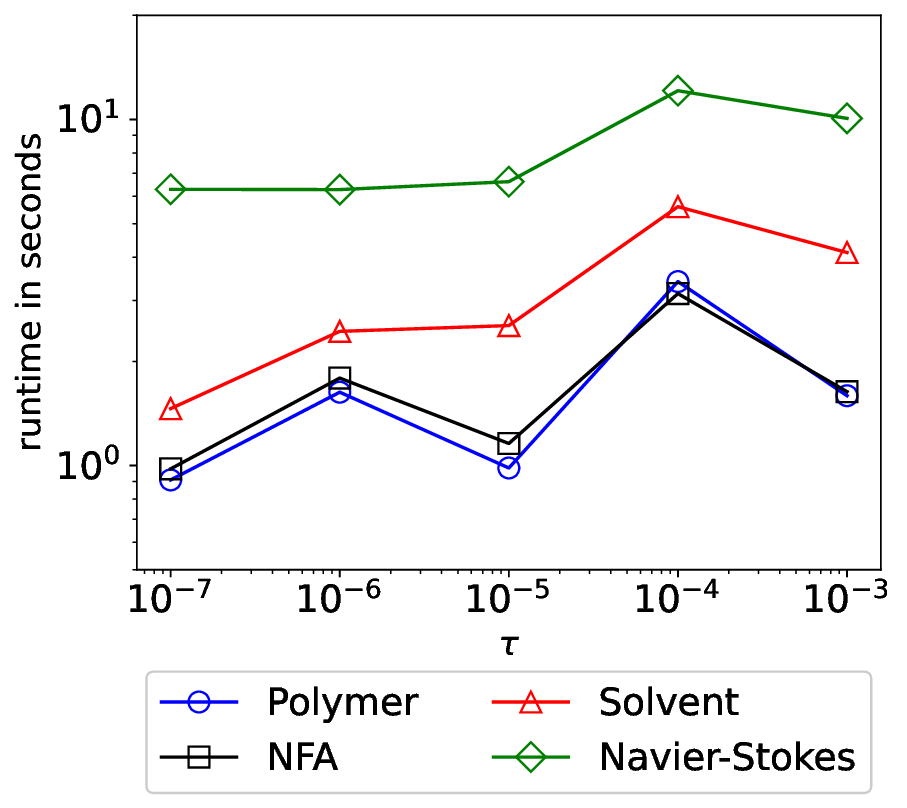}}
\end{minipage}

\caption{\texttt{GMRES} iteration numbers and runtimes per time step for different time step $\tau$ of Cahn-Hilliard (polymer, NFA, solvent) and Navier-Stokes with $n_x=n_y=100$ in a 2D domain. The other parameters are the same as in Figure~\ref{fig:2D_all_gmres_iter_spatial}.}
\label{fig:2D_all_gmres_iter_tau}
\end{figure}

We evaluate the performance of the preconditioned \texttt{GMRES} algorithm. As demonstrated in Section~\ref{subsec:sa}, the algorithm exhibits robust behavior despite varying parameters, see Figures~\ref{fig:2D_all_gmres_iter_spatial} and~\ref{fig:2D_all_gmres_iter_tau}. The runtime scales linearly, while the required \texttt{GMRES} iterations remain within a reasonable range. The obtained morphologies at different times are depicted of Figure~\ref{fig:2D_all_vol_t}. We observe that, although the results for the solvent and evaporation fields show promising results comparable to the findings in Section~\ref{subsec:sa}, the morphologies for polymer and NFA seem highly structured and influenced by meshing parameters. Numerical inaccuracies are evident in these simulations. Since we are dealing with a complex model, there may be many reasons for these errors. One aspect prone to produce numerical instabilities is the complexity of the free energy functional. For this reason, the logarithmic free energy functional of two phase separationg species is often substituted by the much simpler double-well potential of the form $ f^{loc}_{liq/gas}(\phi_1, \phi_2) =  \phi_1^2 \phi_2^2$ , see, e.g., \cite{Elliott1989}. While the function is much simpler, the characteristic minima at the pure phases are retained. Using this double-well potential results in a numerical model less prone to instabilities that still captures the desired phase separation process. In the case of two mixing species, a functional of the form $ f^{loc}_{liq/gas}(\phi_1, \phi_2) =  \phi_1^2 + \phi_2^2$ shows the characteristic minimum at the mixed state. The above reasons promoted us to revisit the approximation of the free energy function. Inspired by the double-well potential and the free energy functionals presented in \cite{bergermann_modeling_2019, Bergermann2023} we searched for a free energy density of the form 
\begin{align*}
f^{loc}_{liq/gas}(\boldsymbol{\phi}) = \sum_i \sum_{j<i} a_{ij} \phi_i^2 \phi_j^2 + \sum_i b_i \phi_i^2 + \sum_i \sum_{j<i} c_{ij} \phi_i \phi_j .
\end{align*}
Performing a numerical fit while keeping the parameters in Figure~\ref{fig:2D_all_gmres_iter_spatial}, 
the functional 
\begin{align}\label{eqn:appPotential}
       f^{loc}(\boldsymbol{\phi}, \Phi_v) &= \left( 1- p(\Phi_v)\right) \left( 3.5 \phi_p^2  \phi_{\nfa}^2 + 0.3\phi_s^2 + 0.3\phi_a^2 - 10\phi_s\phi_a \right)  \nonumber \\ 
       & + p(\Phi_v)  \left( 3.5 \phi_p^2 + 3.5 \phi_{\nfa}^2 + 2\phi_s^2 + 0.1\phi_a^2\right) 
\end{align}
was found to be a good representation of the local free energy density change, while still mimicking all desired characteristics.  

\begin{figure}[htp!]
\begin{minipage}[t]{0.23\textwidth}
\includegraphics[width=1.2\textwidth]{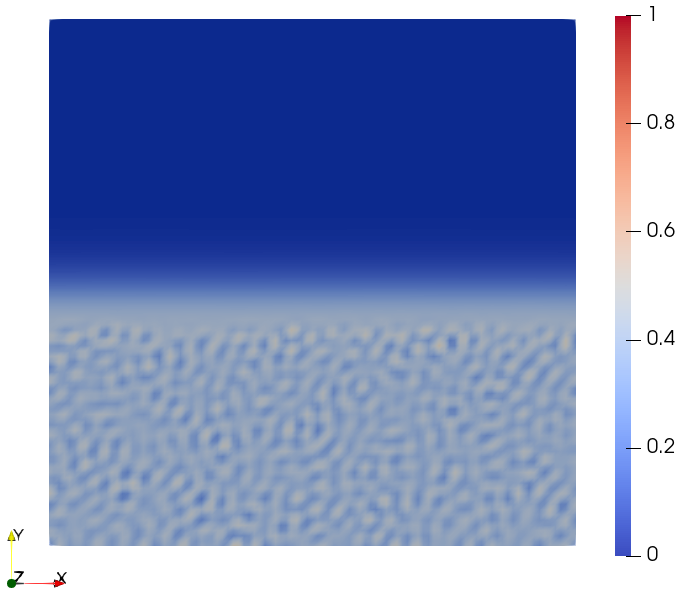}
\end{minipage}
\hfill
\begin{minipage}[t]{0.23\textwidth}
\includegraphics[width=1.2\textwidth]{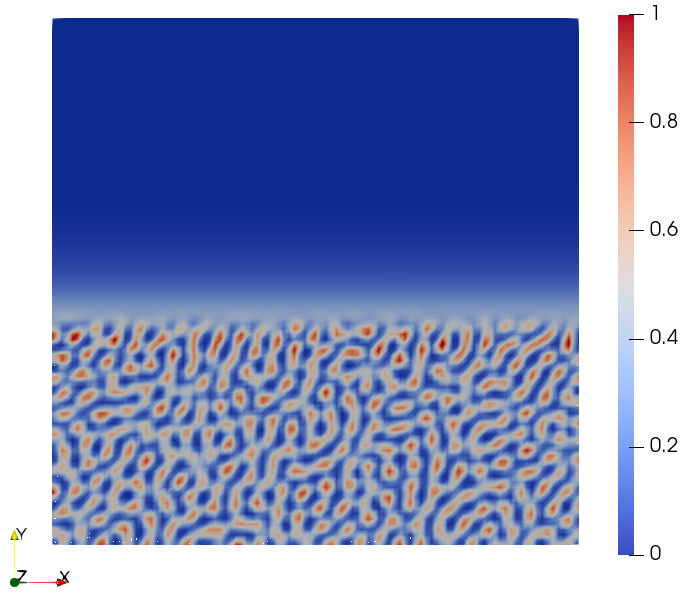}
\end{minipage}
\hfill
\begin{minipage}[t]{0.23\textwidth}
\includegraphics[width=1.2\textwidth]{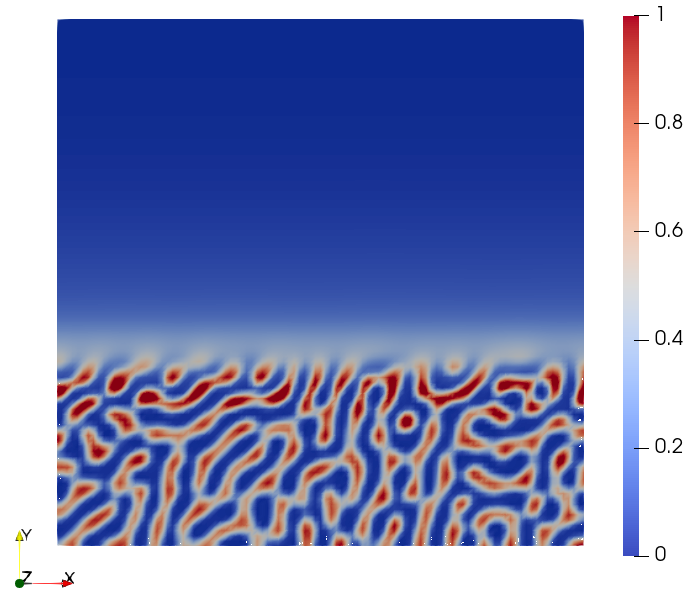}
\end{minipage}
\hfill
\begin{minipage}[t]{0.23\textwidth}
\includegraphics[width=1.2\textwidth]{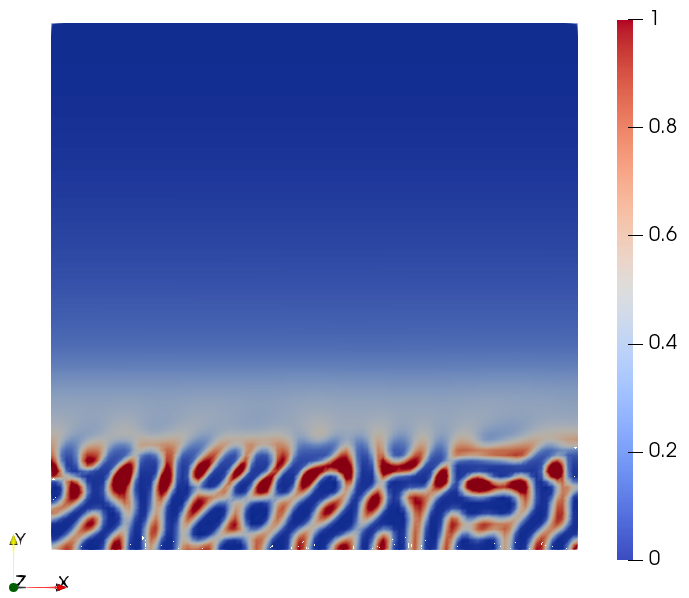}
\end{minipage}

\begin{minipage}[t]{0.23\textwidth}
\includegraphics[width=1.2\textwidth]{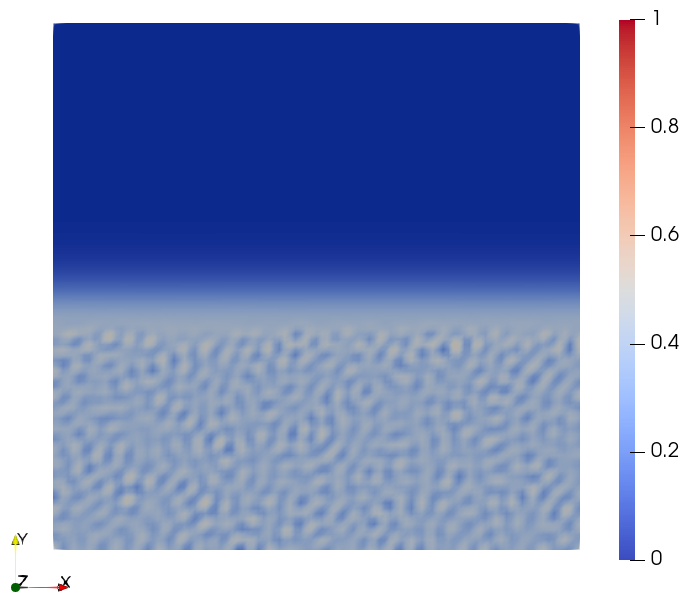}
\end{minipage}
\hfill
\begin{minipage}[t]{0.23\textwidth}
\includegraphics[width=1.2\textwidth]{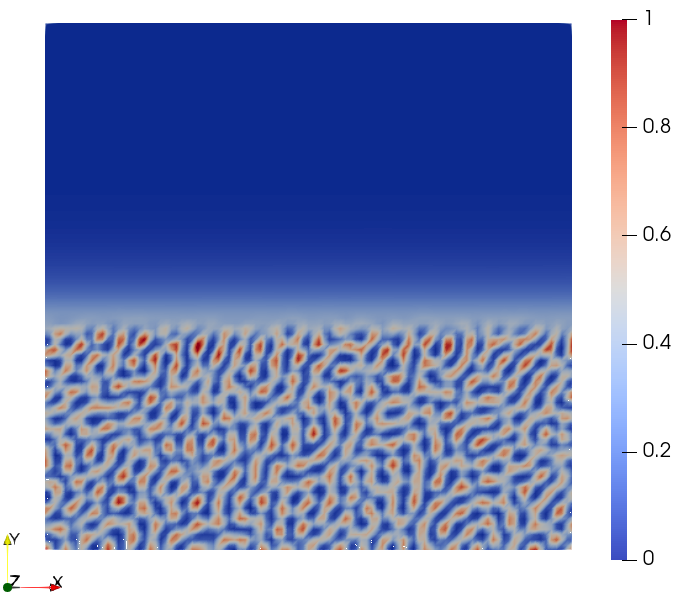}
\end{minipage}
\hfill
\begin{minipage}[t]{0.23\textwidth}
\includegraphics[width=1.2\textwidth]{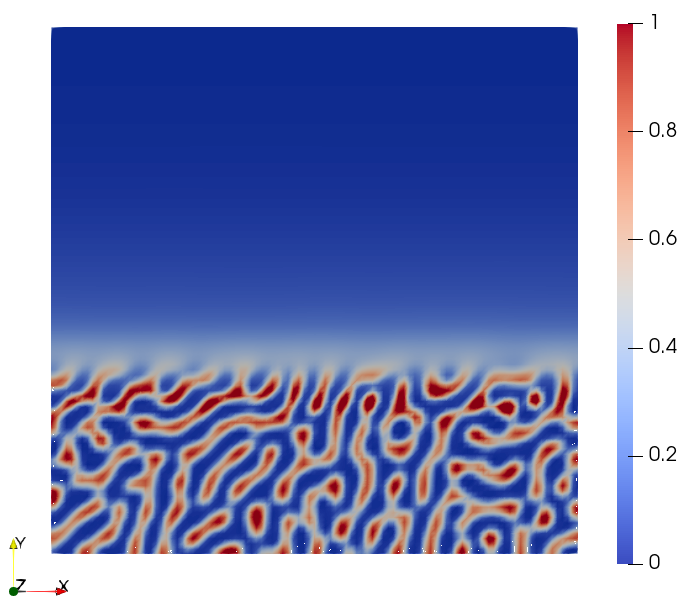}
\end{minipage}
\hfill
\begin{minipage}[t]{0.23\textwidth}
\includegraphics[width=1.2\textwidth]{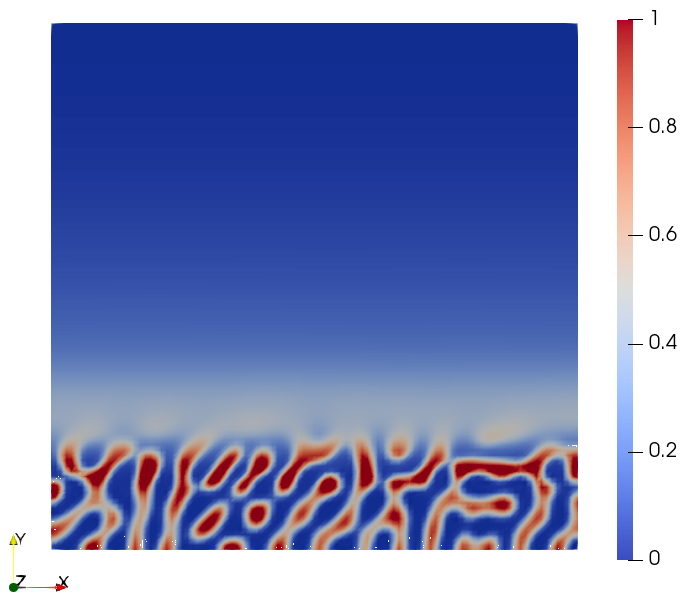}
\end{minipage}

\begin{minipage}[t]{0.23\textwidth}
\includegraphics[width=1.2\textwidth]{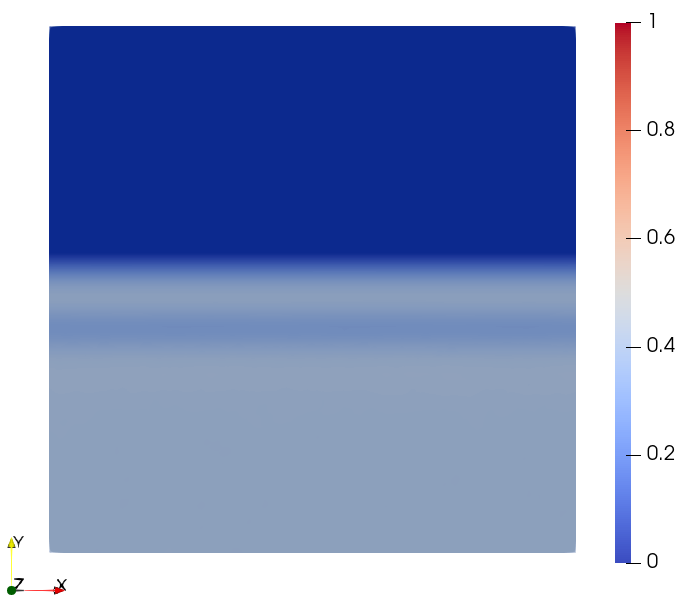}
\end{minipage}
\hfill
\begin{minipage}[t]{0.23\textwidth}
\includegraphics[width=1.2\textwidth]{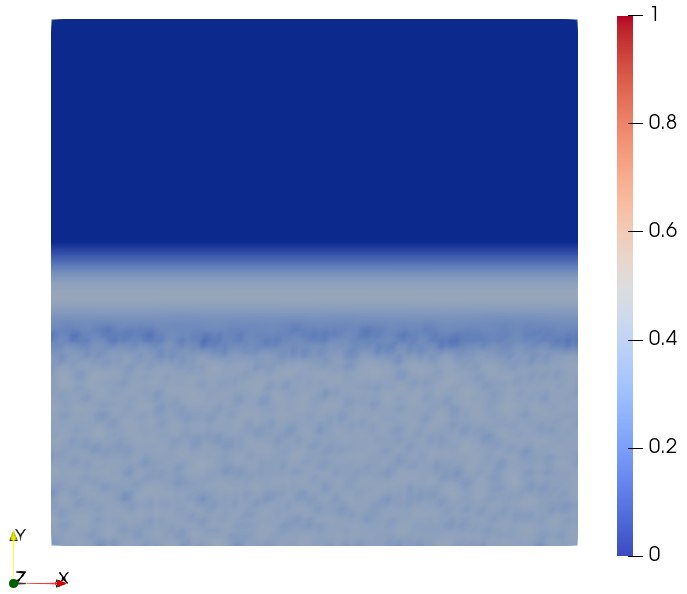}
\end{minipage}
\hfill
\begin{minipage}[t]{0.23\textwidth}
\includegraphics[width=1.2\textwidth]{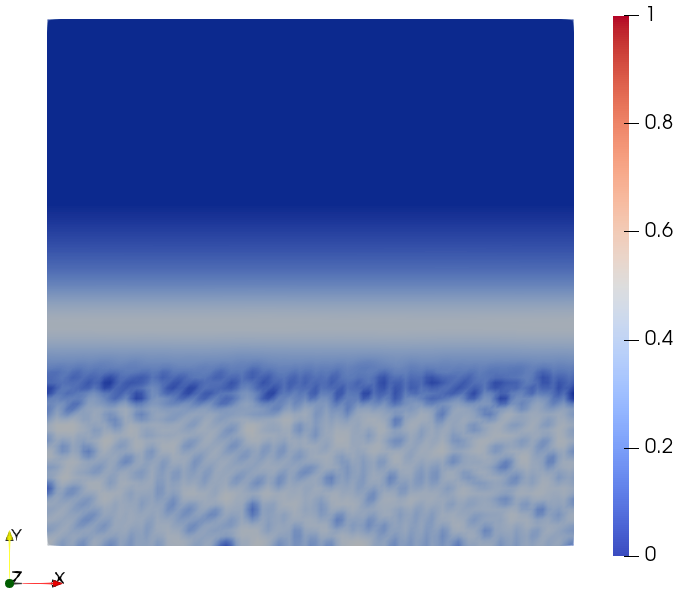}
\end{minipage}
\hfill
\begin{minipage}[t]{0.23\textwidth}
\includegraphics[width=1.2\textwidth]{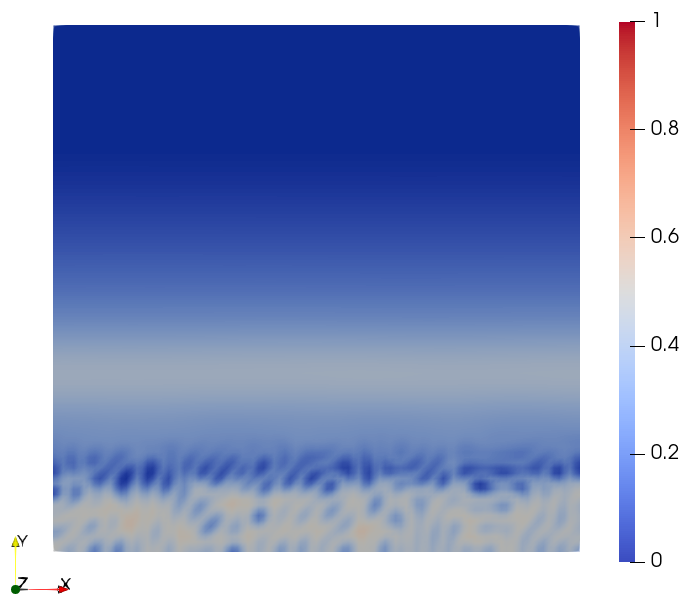}
\end{minipage}

\begin{minipage}[t]{0.23\textwidth}
\includegraphics[width=1.2\textwidth]{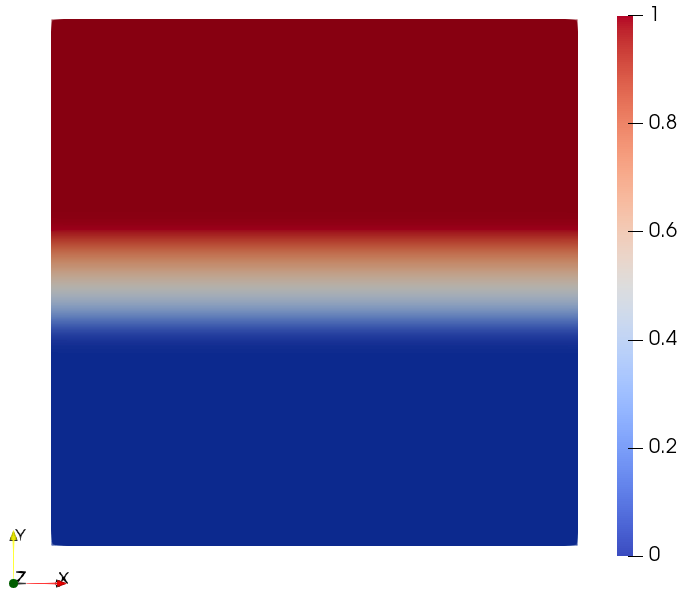}
\end{minipage}
\hfill
\begin{minipage}[t]{0.23\textwidth}
\includegraphics[width=1.2\textwidth]{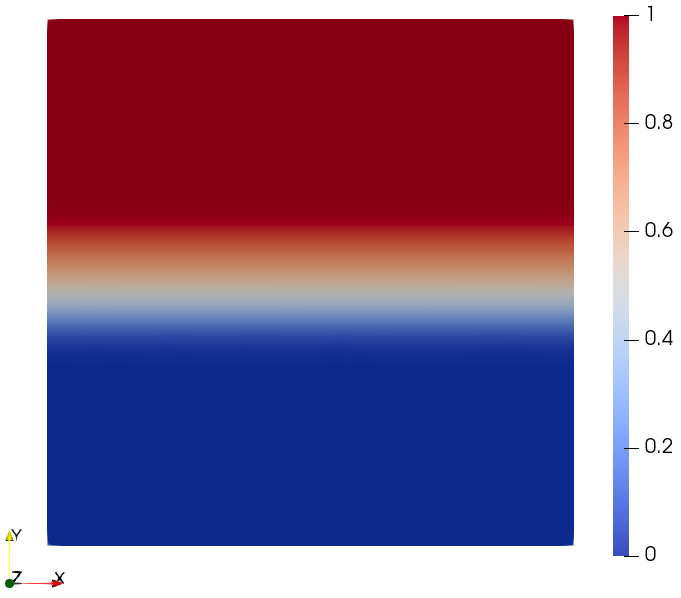}
\end{minipage}
\hfill
\begin{minipage}[t]{0.23\textwidth}
\includegraphics[width=1.2\textwidth]{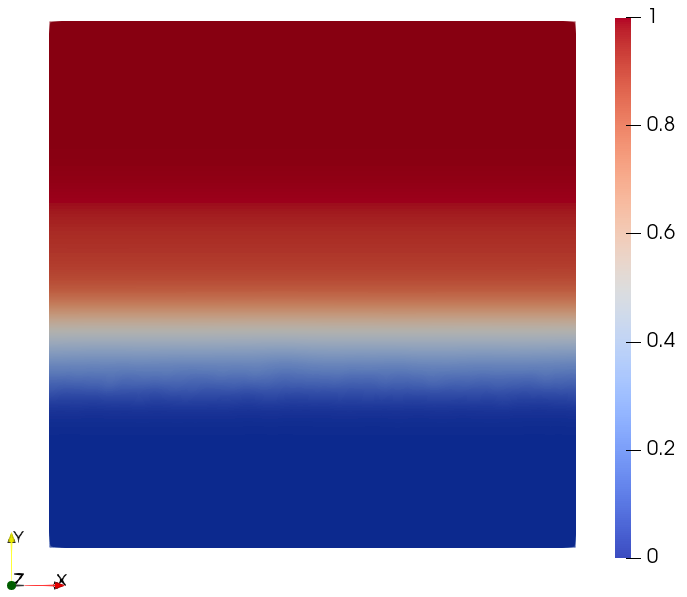}
\end{minipage}
\hfill
\begin{minipage}[t]{0.23\textwidth}
\includegraphics[width=1.2\textwidth]{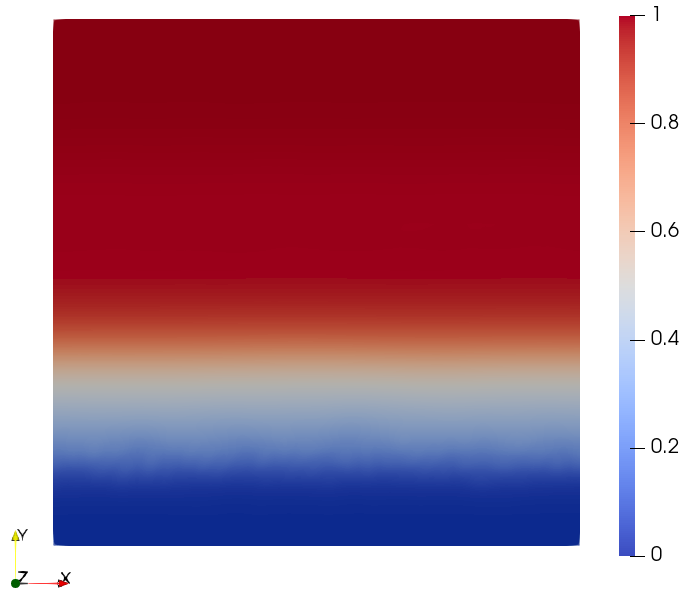}
\end{minipage}
\caption{Typical volume fraction field polymer, NFA, solvent, and order parameter field vapor (from top to bottom) during evaporation. The time step size $10^{-4}$ is used and the results shown  correspond to $t = 0.025$, $t =0.05$, $t = 0.25$, $t = 0.5$. The spatial discretization of a 2D domain is chosen $n_x = n_y =100$. The parameters and setup are the same as in Figure~\ref{fig:2D_all_gmres_iter_spatial}, and the polynomial approximation \eqref{eqn:appPotential} of the logarithmic potential are used. }
\label{fig:2D_all_vol_newf}
\end{figure}

By utilizing this new approximation, the morphologies in Figure~\ref{fig:2D_all_vol_newf} are obtained. 
The obtained numerical results show the expected phase separation of polymer and NFA, whilst the phase fields of solvent and vapor order parameter still show the behavior already presented in Section~\ref{subsec:sa}. To highlight the efficiency of the derived preconditioning system, we extend the numerical simulations to a three-dimensional setting. 

In the 3D domain, we use the Smoothed Aggregation solver for \texttt{AMG} to efficiently precondition the Navier-Stokes equations for both symmetric and non-symmetric matrices, thereby reducing runtime. All 3D results are obtained by using the polynomial approximation \eqref{eqn:appPotential} of the logarithmic potential. 

\begin{figure}[htp!]
\begin{minipage}[t]{0.48\textwidth}
\includegraphics[width=0.85\textwidth]{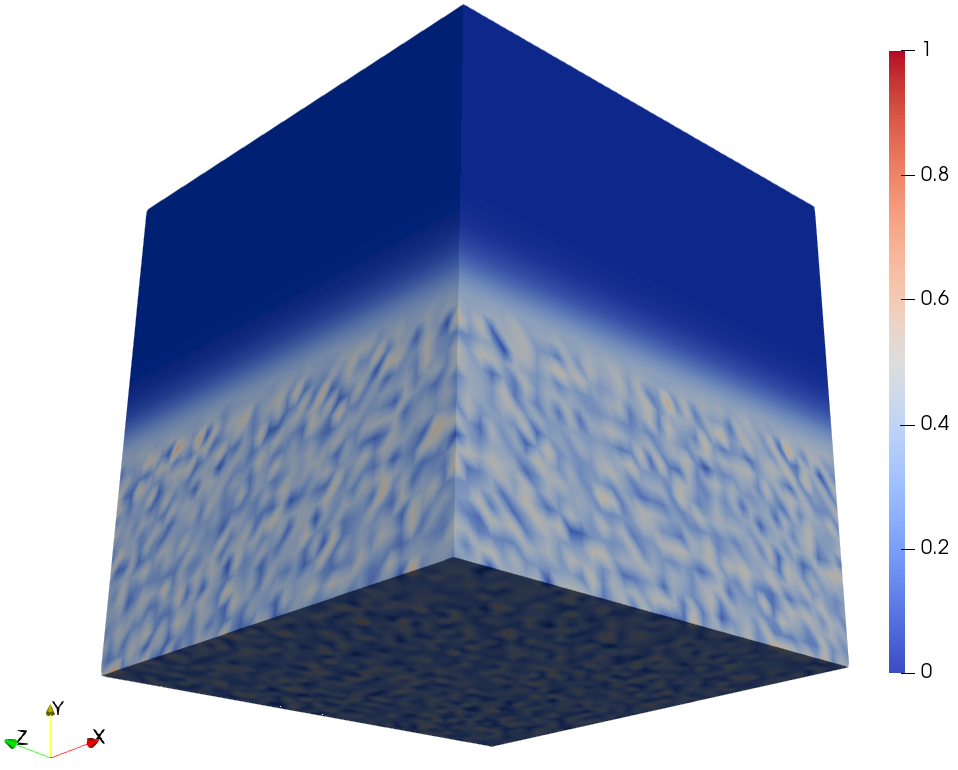}
\end{minipage}
\hfill
\begin{minipage}[t]{0.48\textwidth}
\includegraphics[width=0.85\textwidth]{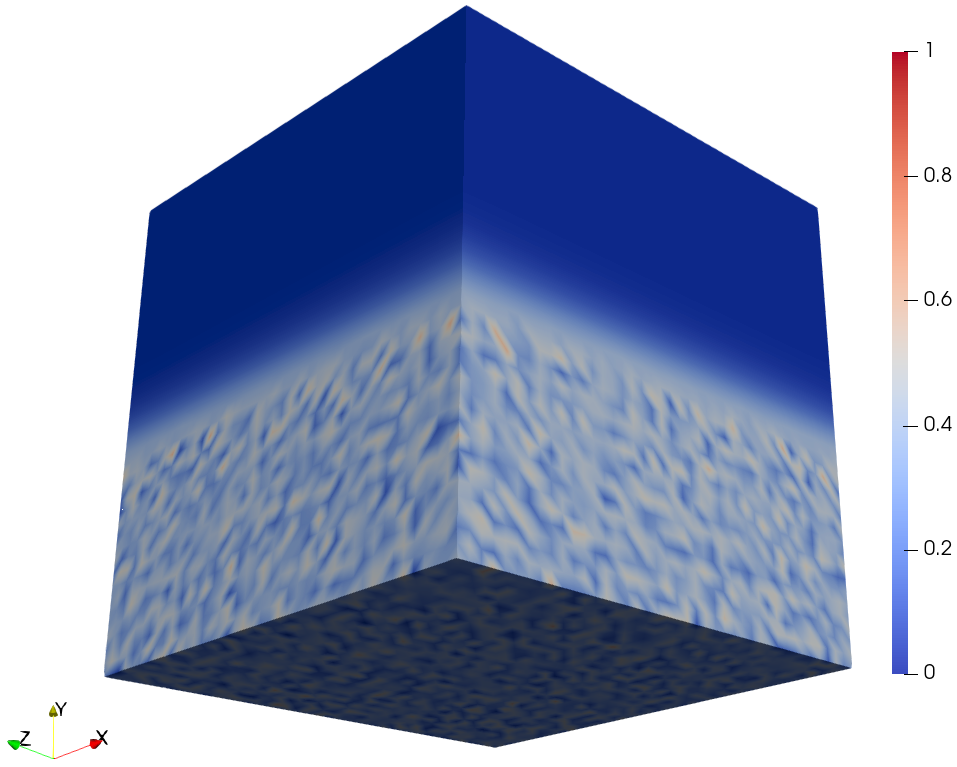}
\end{minipage}

\begin{minipage}[t]{0.48\textwidth}
\includegraphics[width=0.85\textwidth]{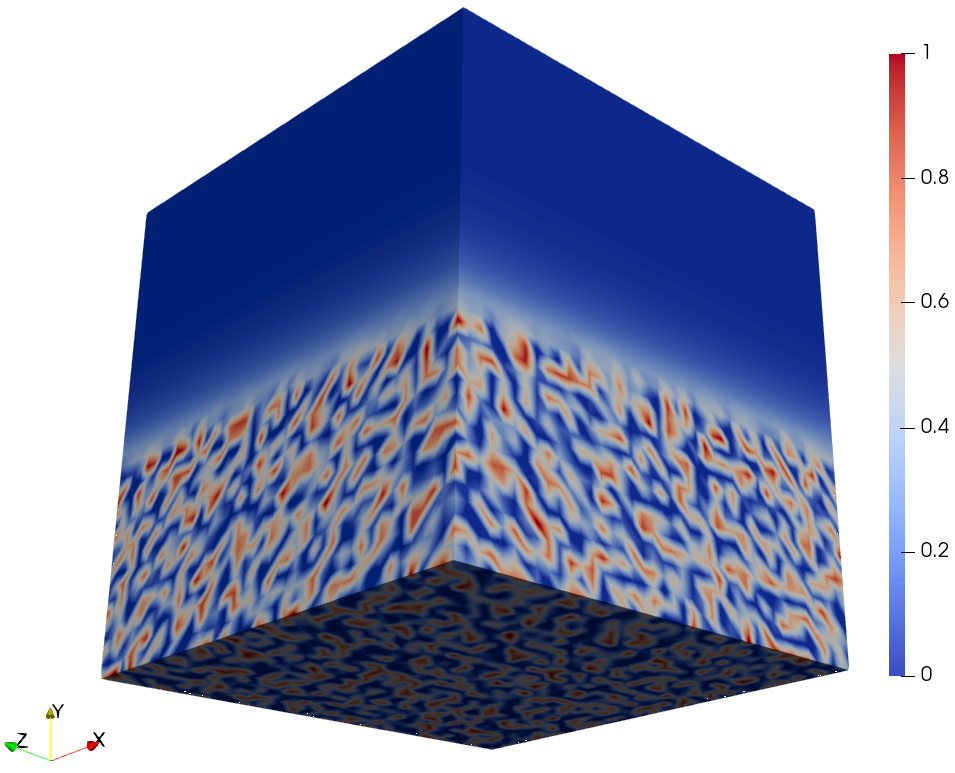}
\end{minipage}
\hfill
\begin{minipage}[t]{0.48\textwidth}
\includegraphics[width=0.85\textwidth]{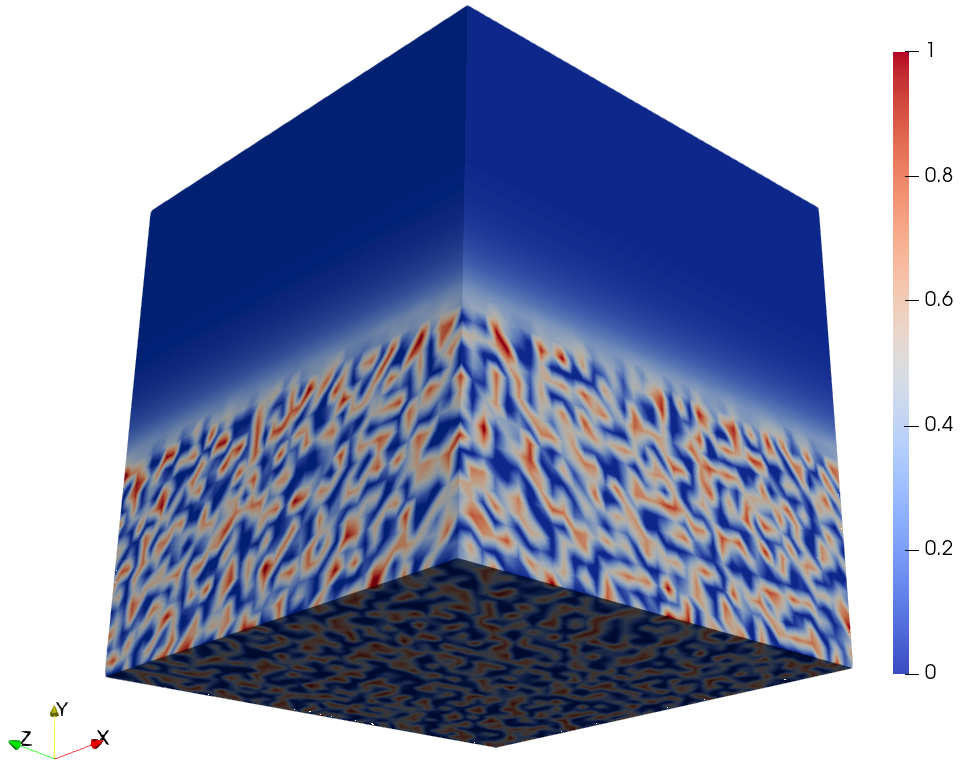}
\end{minipage}

\begin{minipage}[t]{0.48\textwidth}
\includegraphics[width=0.85\textwidth]{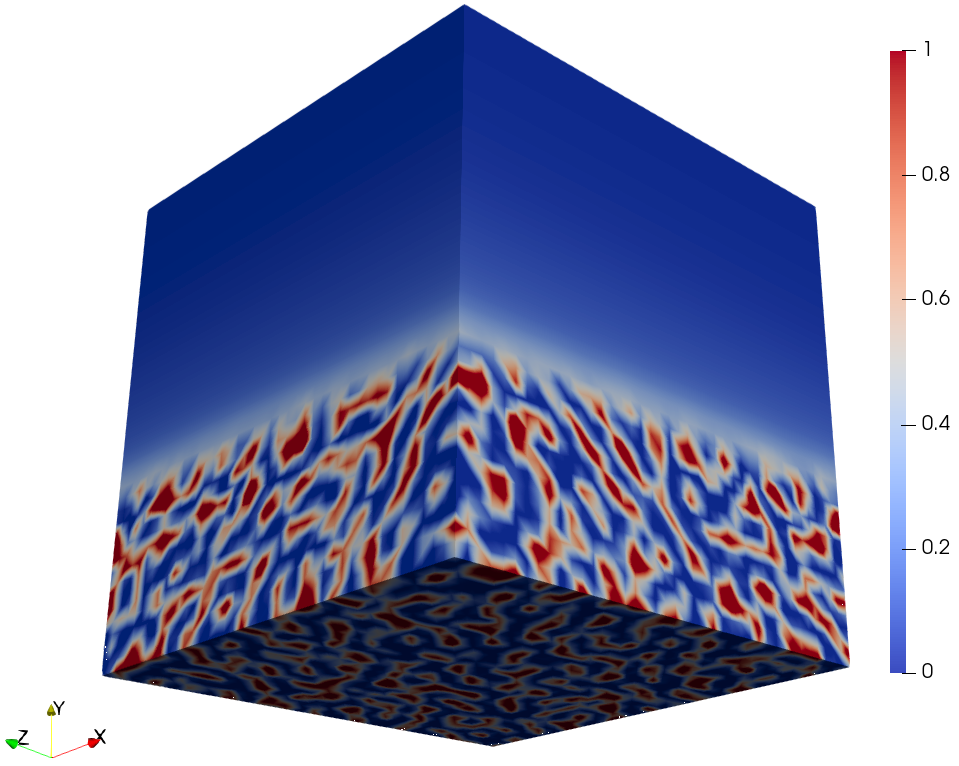}
\end{minipage}
\hfill
\begin{minipage}[t]{0.48\textwidth}
\includegraphics[width=0.85\textwidth]{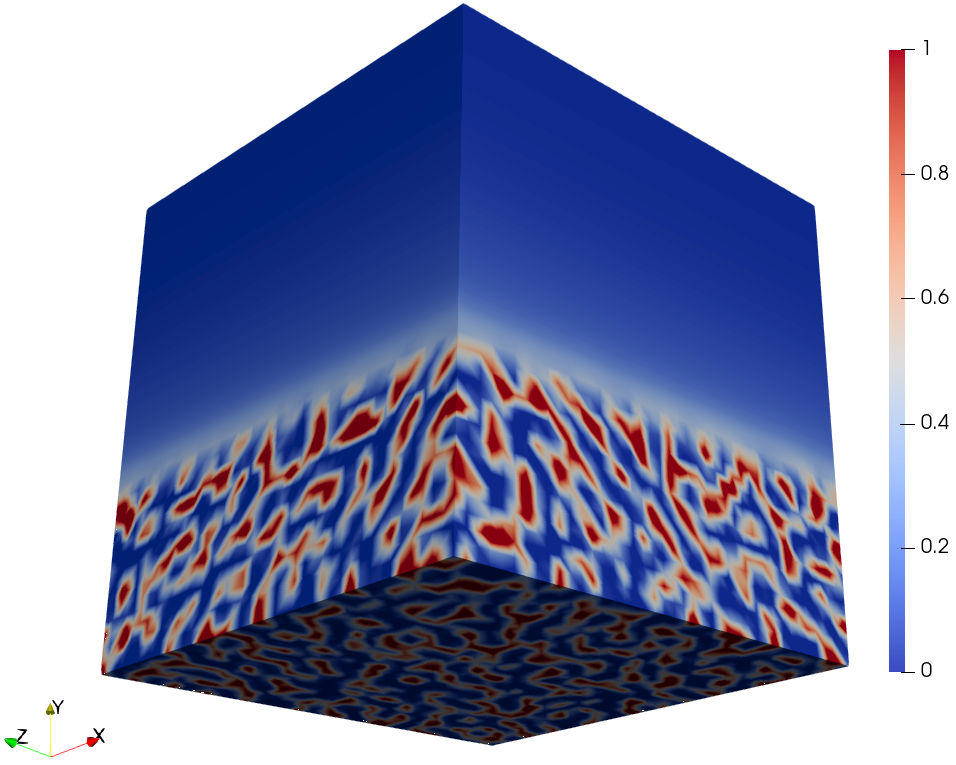}
\end{minipage}

\begin{minipage}[t]{0.48\textwidth}
\includegraphics[width=0.85\textwidth]{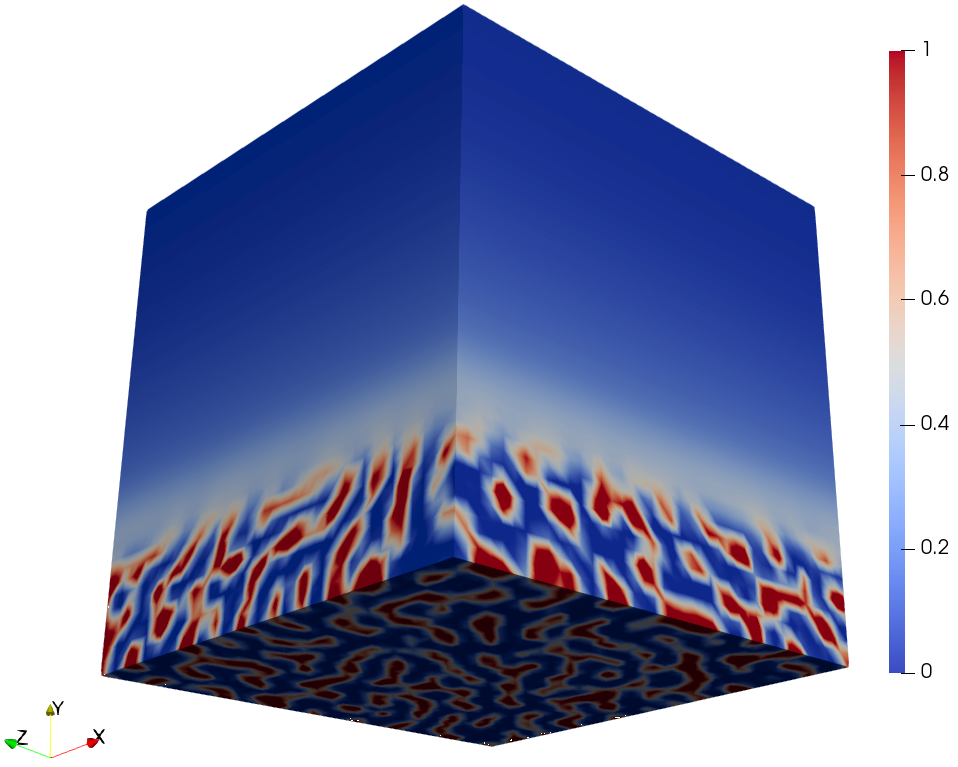}
\end{minipage}
\hfill
\begin{minipage}[t]{0.48\textwidth}
\includegraphics[width=0.85\textwidth]{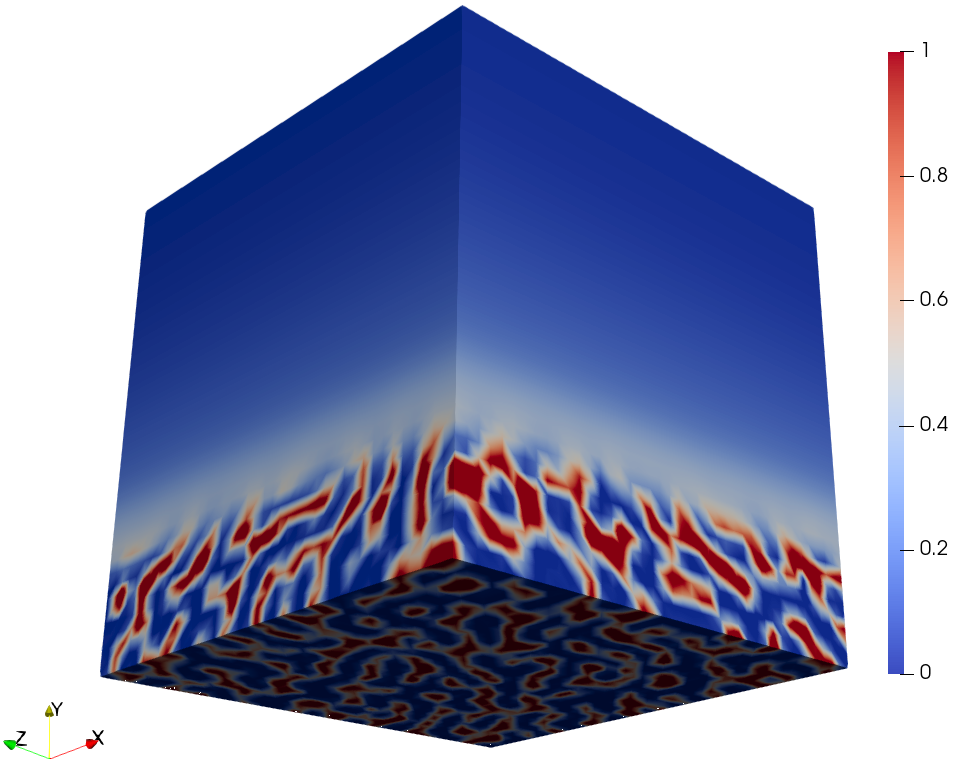}
\end{minipage}

\caption{Typical volume fraction field polymer (left) and NFA (right) during evaporation. The time step size $10^{-4}$ is used and the results shown  correspond to $t = 0.025$, $t =0.05$, $t = 0.25$, $t = 0.5$ (from top to bottom). The spatial discretization of a 3D domain is chosen $n_x = n_y = n_z = 40$. Moreover, the other parameters are the same as in Figure~\ref{fig:2D_all_gmres_iter_spatial}. }
\label{fig:3D_all_vol_t}
\end{figure}

\begin{figure}[t!]
\begin{minipage}[t]{0.32\textwidth}
\includegraphics[width=1.05\textwidth]{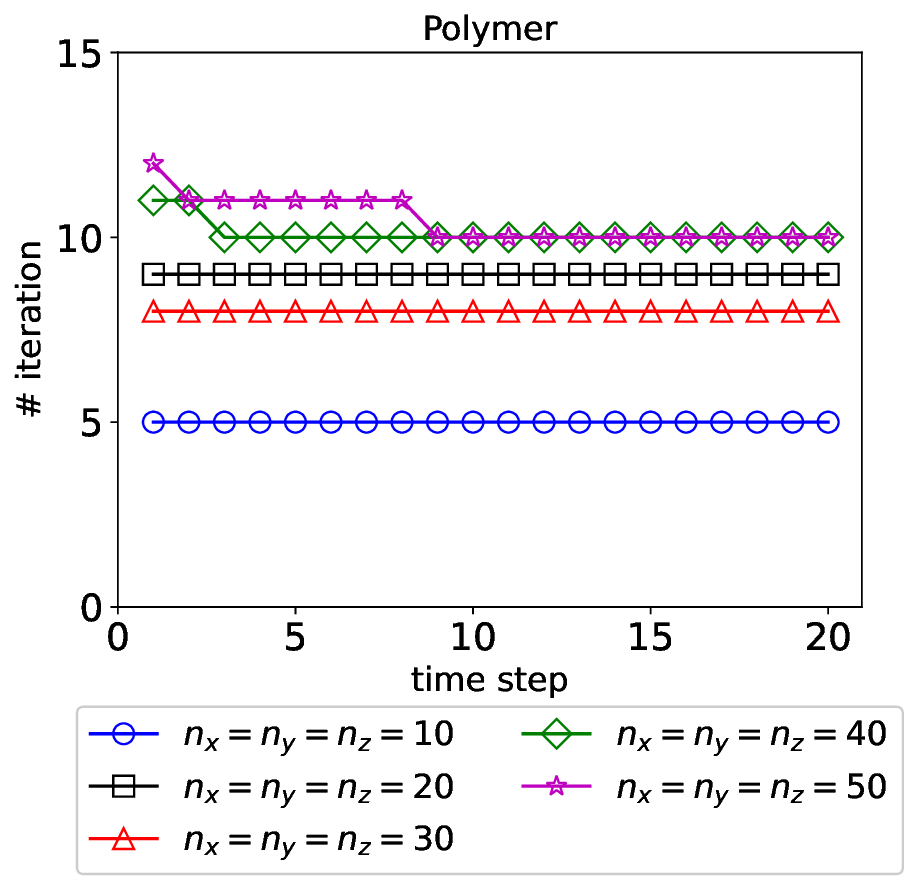}
\end{minipage}
\begin{minipage}[t]{0.32\textwidth}
\includegraphics[width=1.05\textwidth]{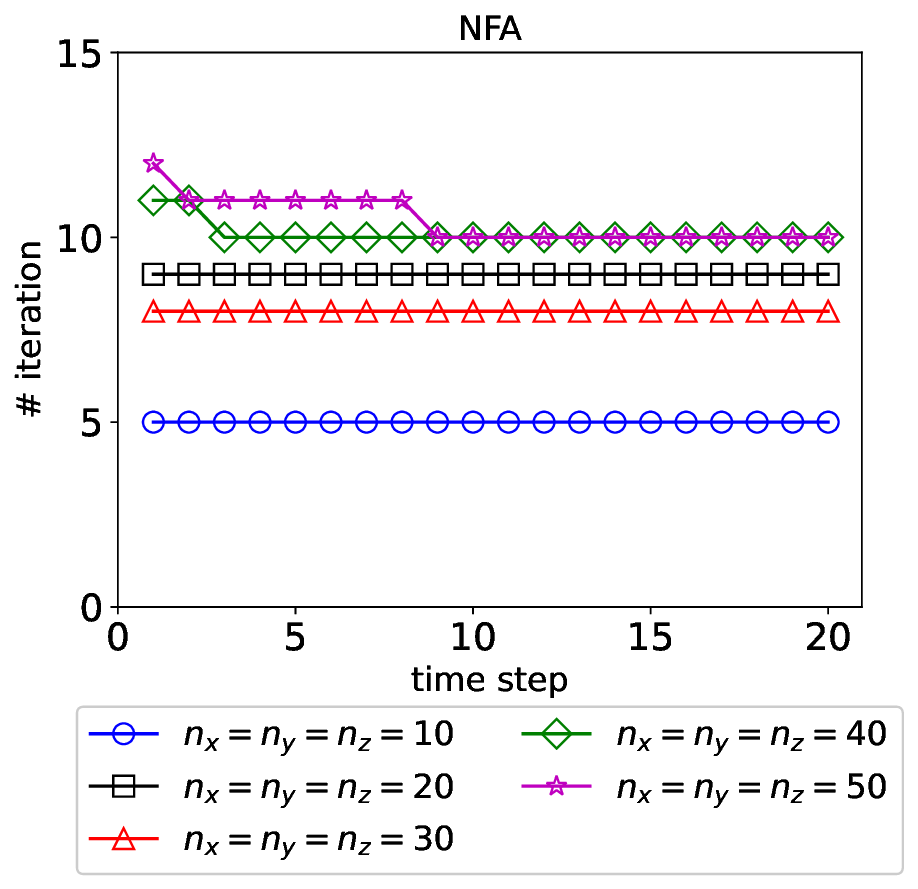}
\end{minipage}
\begin{minipage}[t]{0.32\textwidth}
\includegraphics[width=1.05\textwidth]{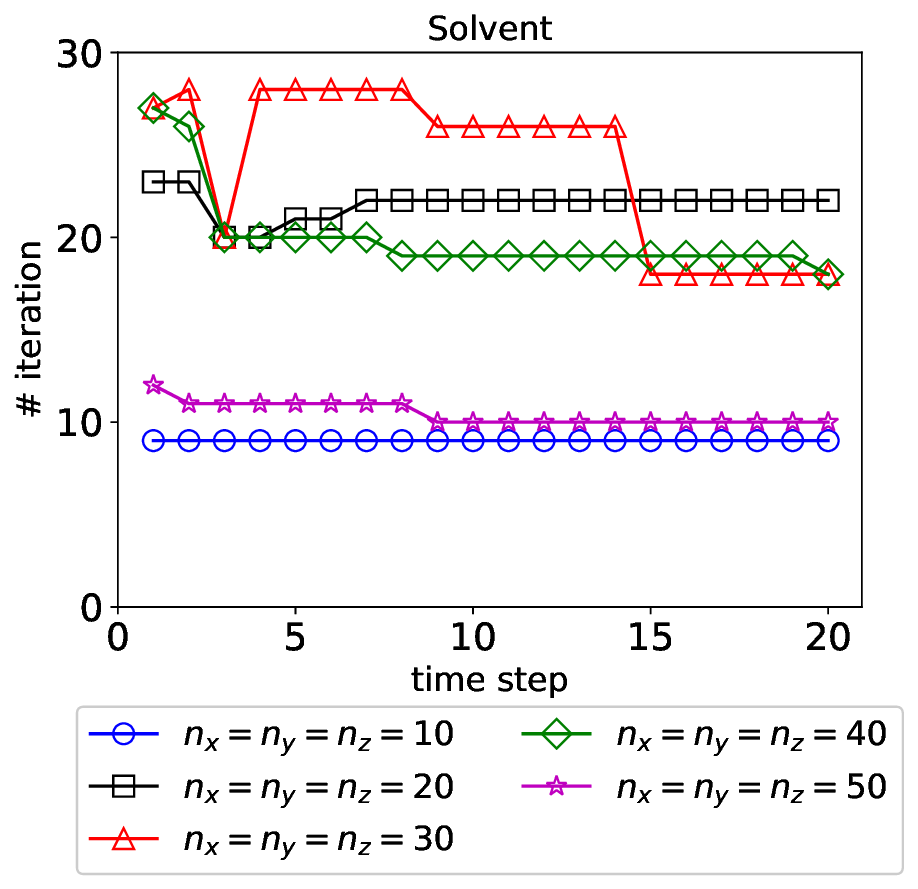}
\end{minipage}

\begin{minipage}[t]{0.32\textwidth}
\includegraphics[width=1.05\textwidth]{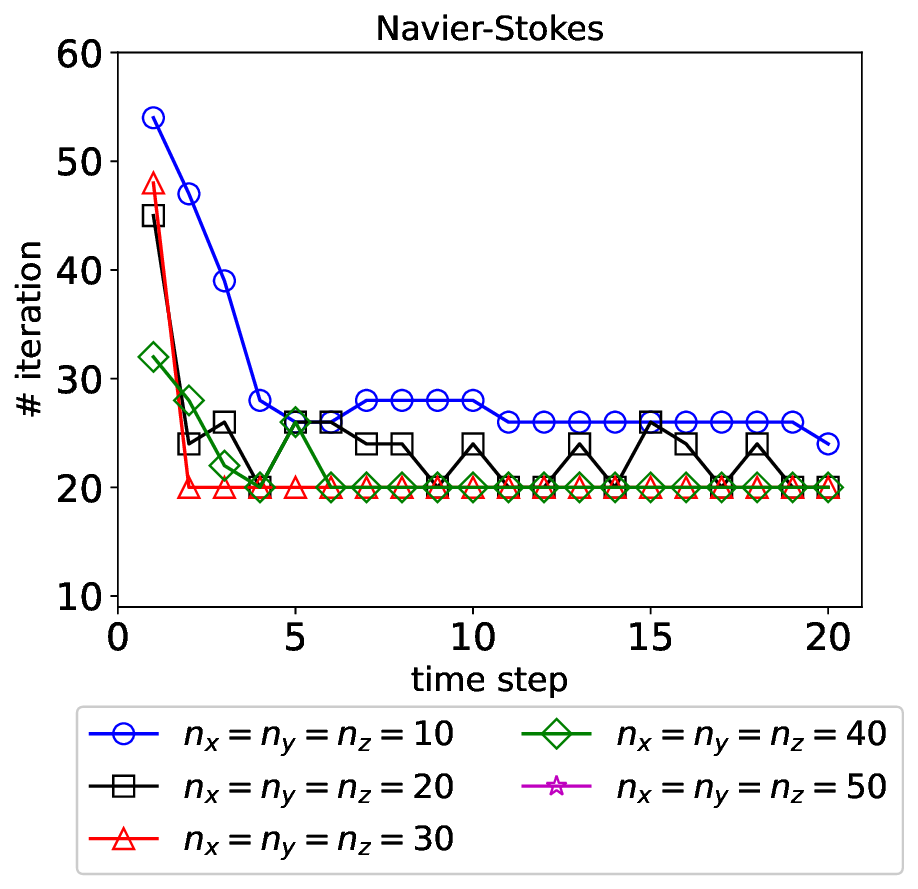}
\end{minipage}
\begin{minipage}[t]{0.32\textwidth}
\includegraphics[width=1.1\textwidth]{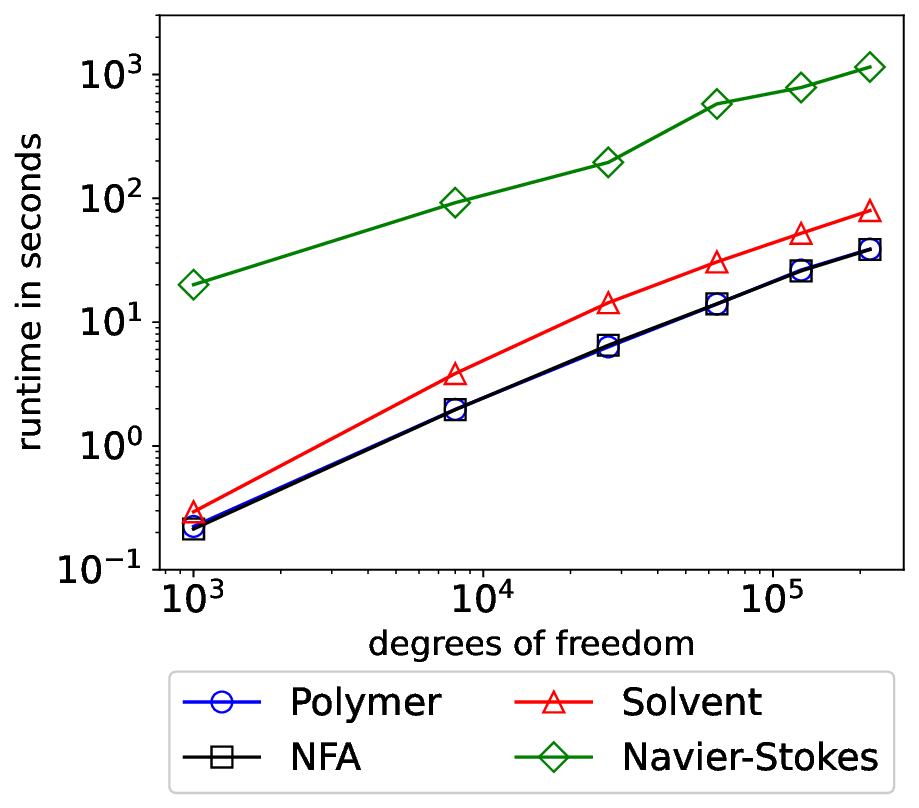}
\end{minipage}
\caption{\texttt{GMRES} iteration numbers and runtimes per time step for different spatial discretizations of Cahn-Hilliard (polymer, NFA, solvent) and Navier-Stokes with $\tau = 10^{-4}$ in a 3D domain. The other parameters and setup are the same as in Figure~\ref{fig:3D_all_vol_t}.}
\label{fig:3D_all_gmres_iter_spatial}
\end{figure}

Figure~\ref{fig:3D_all_vol_t} illustrates the 3D morphology evolution at four different times, with parameter choices detailed in the caption. To demonstrate the robustness of our solver in the 3D domain, we report the number of \texttt{GMRES} iterations and the runtime required to meet the tolerance on a given mesh in Figure~\ref{fig:3D_all_gmres_iter_spatial}. These results align with observations from previous experiments.

\section{Conclusion and Outlook}
In this paper, we presented a phase-field model describing the morphology evolution of organic solar cells. The resulting model consists of coupled Navier-Stokes, Cahn-Hilliard, and Allen-Cahn equations. To deal with the large-scale, poorly conditioned linear systems, a preconditioning technique was introduced, handling the system ensuring robust handling of the system with respect to variations in the discretization parameters. We have developed and analyzed block-diagonal preconditioners that utilize a Schur complement approximation, which is efficient to implement. This approximation is combined with multilevel methods, such as \texttt{AMG}. Numerous numerical results have been presented for the evaporation of solvent in the simple model, and for phase separation and evaporation in the complex model, across 1D, 2D, and 3D domains. Being able to model three-dimensional morphology formation is an important tool in finding optimized production conditions and can possibly improve OPV performances. 

To express significant conclusions on the effects of production or material properties to the production process, as a next step, physical relevant parameters need to be introduced. Yielding possible large parameter intervals, this provides new difficulties to the derived numerical system that need to be treated with care. Often these range over several orders of magnitude, which could likely require the application and design of multiscale numerical methods regarding discretization and iterative solution. 



\section*{Declaration of competing interest}
The authors declare that they have no known competing financial interests or personal relationships that could have appeared to influence the work reported in this paper.

\section*{Acknowledgements}
We thank the Deutsche Forschungsgemeinschaft (DFG) for funding this work (Research
Unit FOR 5387 POPULAR, Project No. 461909888).

\section*{Data availability}
No data was used for the research described in the article.


\appendix
\section{ } \label{sec:appendix}
\begin{figure}[htp!]
\begin{minipage}[t]{0.32\textwidth}
\includegraphics[width=1.05\textwidth]{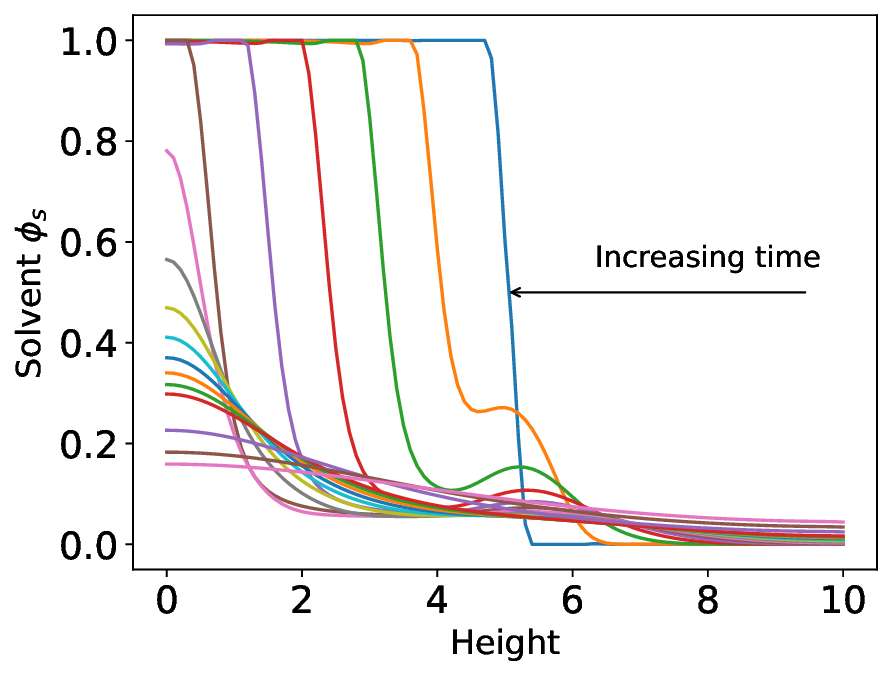}
\end{minipage}
\hfill
\begin{minipage}[t]{0.32\textwidth}
\includegraphics[width=1.05\textwidth]{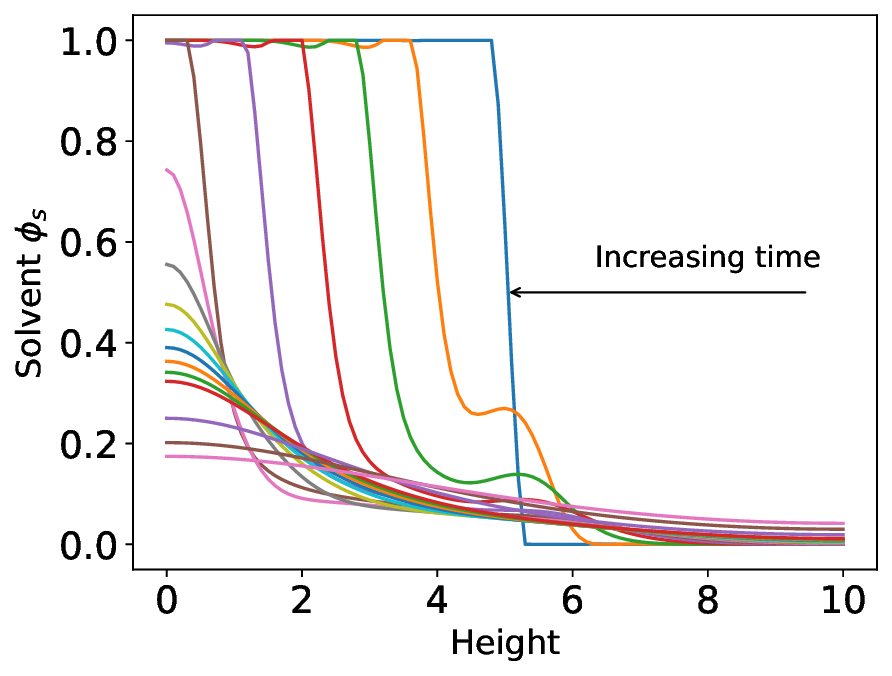}
\end{minipage}
\hfill
\begin{minipage}[t]{0.32\textwidth}
\includegraphics[width=1.05\textwidth]{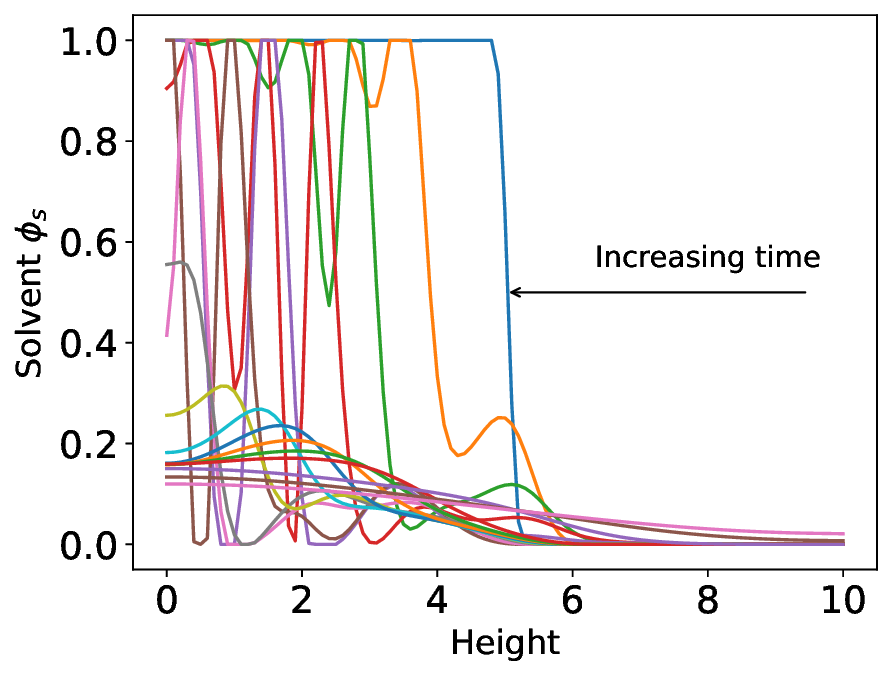}
\end{minipage}

\begin{minipage}[t]{0.32\textwidth}
\includegraphics[width=1.05\textwidth]{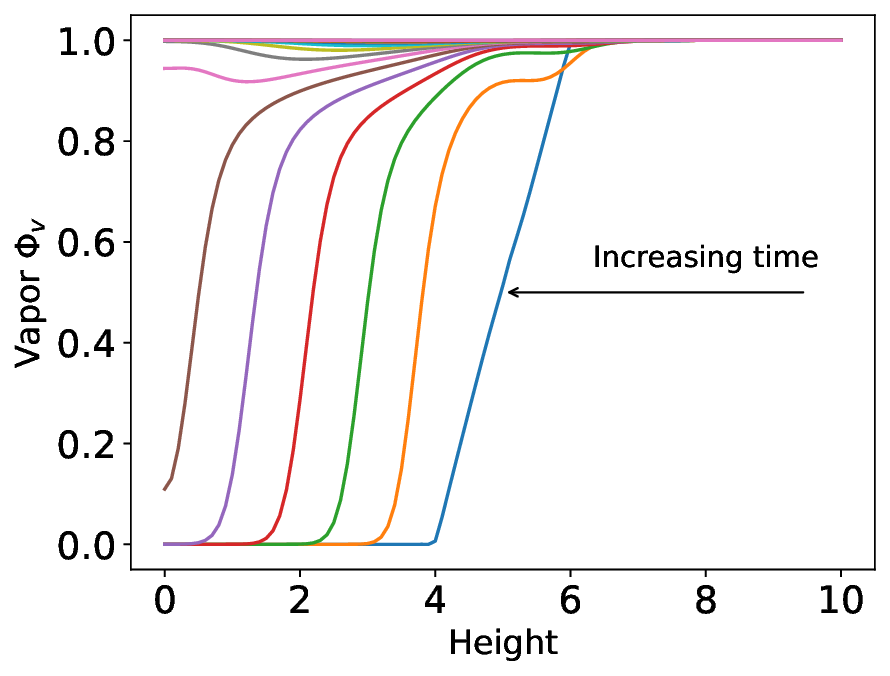}
\end{minipage}
\hfill
\begin{minipage}[t]{0.32\textwidth}
\includegraphics[width=1.05\textwidth]{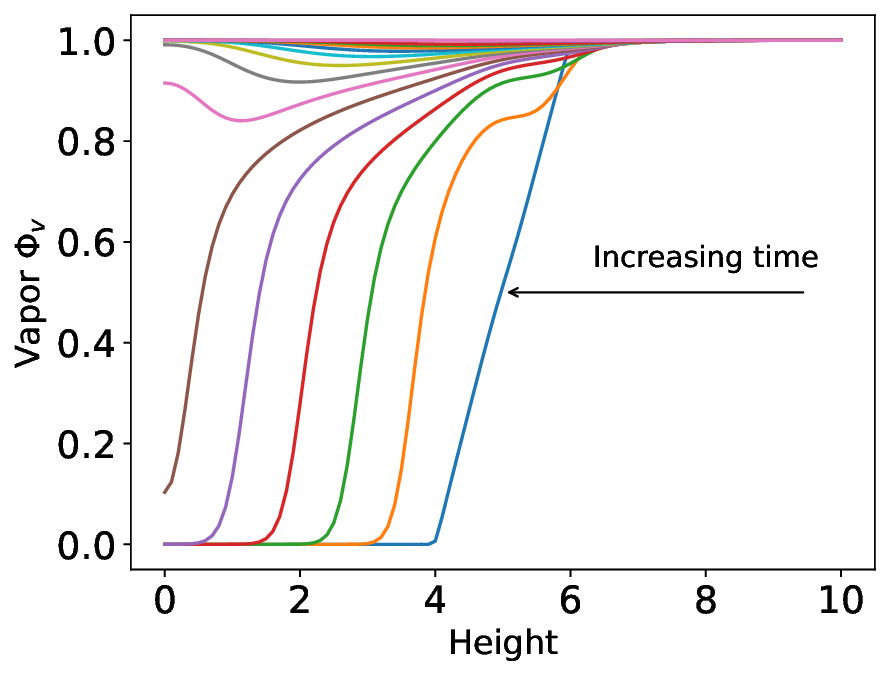}
\end{minipage}
\hfill
\begin{minipage}[t]{0.32\textwidth}
\includegraphics[width=1.05\textwidth]{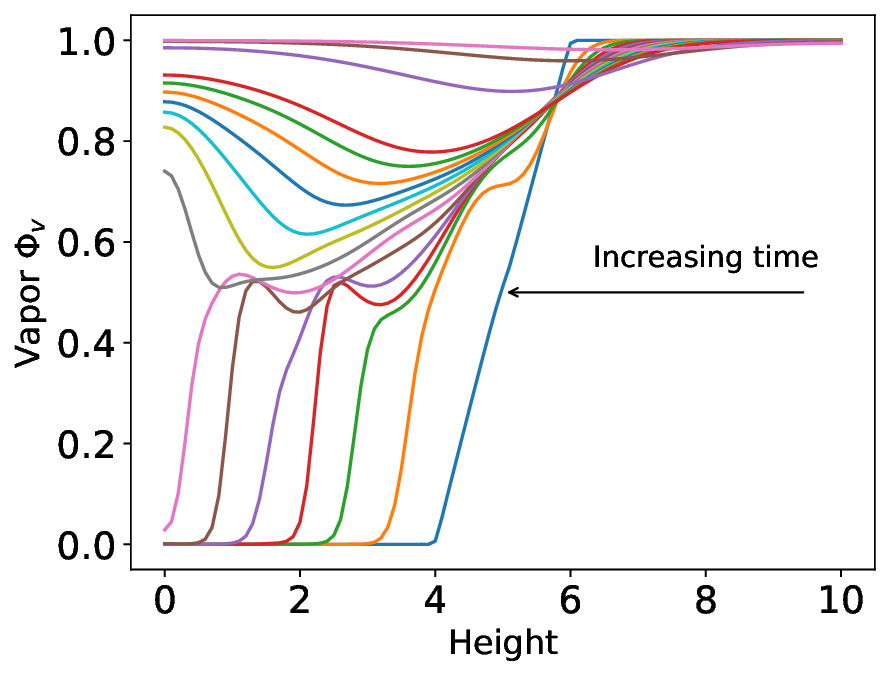}
\end{minipage}
\caption{Volume fraction fields for the solvent (top) and the vapor (bottom)  with $\tau = 10^{-3}$ and $\chi_{s,a} = -10, 0 , 10$  from left to right. The spatial discretization a 1D domain is chosen $n_y =100$. Moreover, the parameters $\alpha_i = 10^{-8}$, $\beta_i = 10^{-1}$ for $ i \in \{ s, a \}$,  $\beta_v = 1$, $\gamma_i = 1$ for $ i \in \{ s, a, v \}$, and $\delta_v = 1$, the final time $t_{\mathrm{max}} = 2.5$, and the molar size of the fluid $N_s=N_s =1$ are used. }
\label{fig:1D_inter}
\end{figure}

In Figure~\ref{fig:1D_inter}, we show the evolution of solvent and vapor volume fractions for three different interaction parameters $\chi_{s,a} : - 10, \, 0, \, 10$. We observe that increasing the interaction parameter causes the fluctuations in the solution.  

\begin{figure}[h!]
\begin{minipage}[t]{0.23\textwidth}
\includegraphics[width=1.15\textwidth]{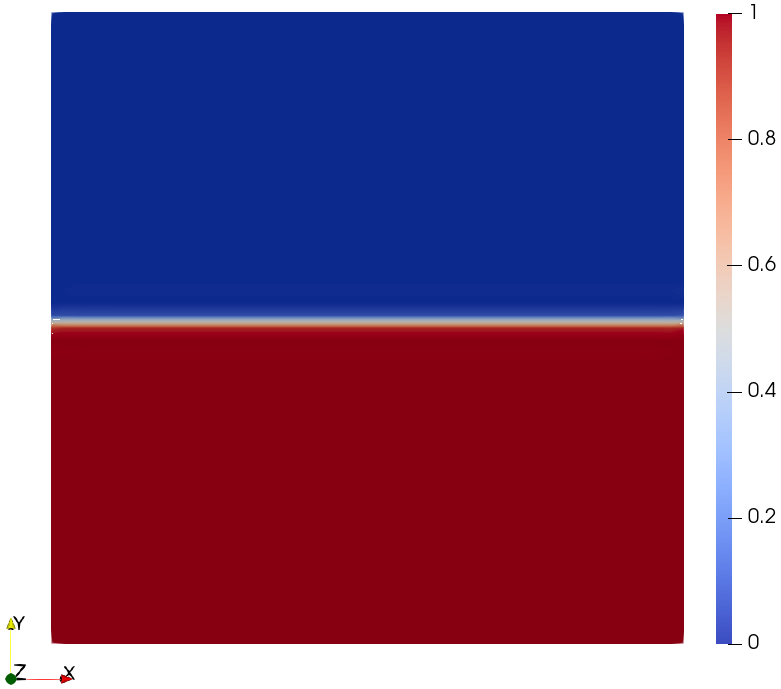}
\end{minipage}
\hfill
\begin{minipage}[t]{0.23\textwidth}
\includegraphics[width=1.15\textwidth]{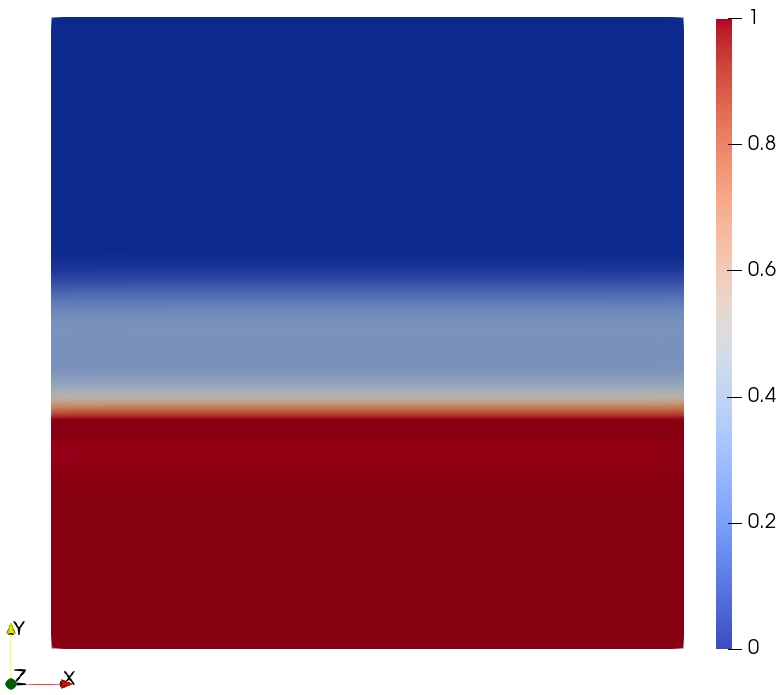}
\end{minipage}
\hfill
\begin{minipage}[t]{0.23\textwidth}
\includegraphics[width=1.15\textwidth]{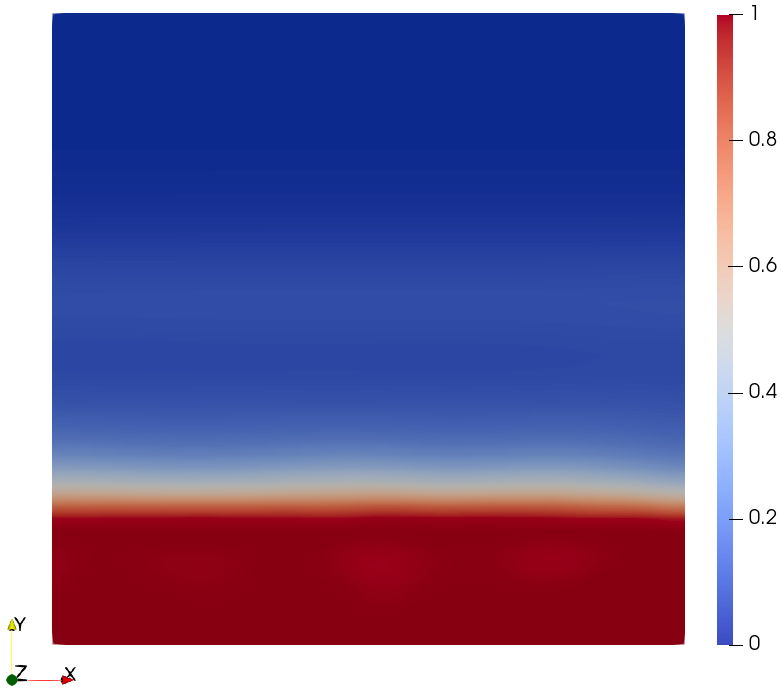}
\end{minipage}
\hfill
\begin{minipage}[t]{0.23\textwidth}
\includegraphics[width=1.15\textwidth]{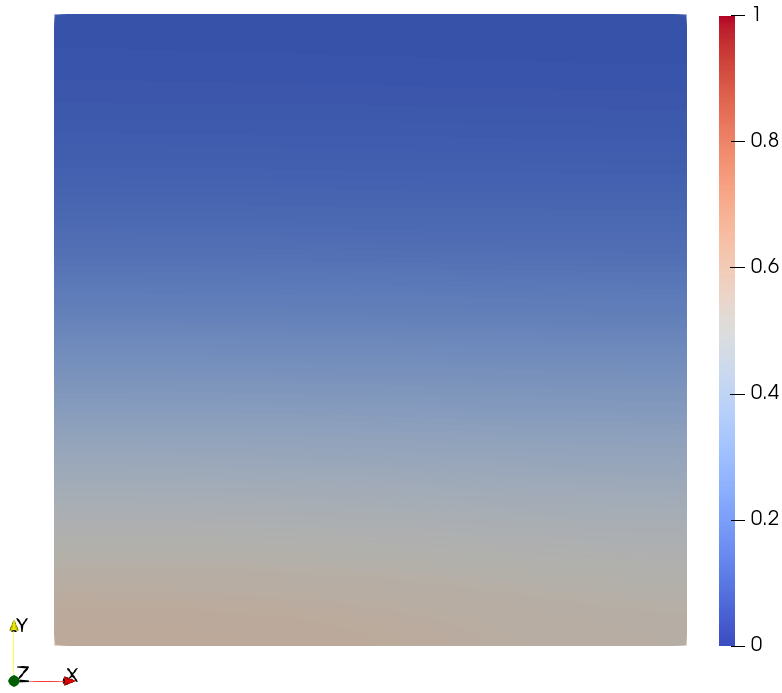}
\end{minipage}

\begin{minipage}[t]{0.23\textwidth}
\includegraphics[width=1.15\textwidth]{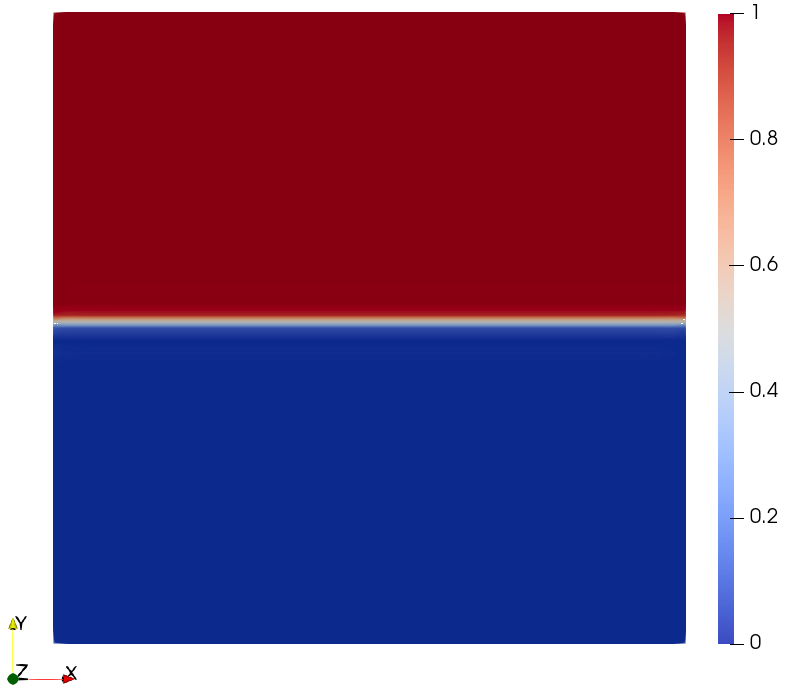}
\end{minipage}
\hfill
\begin{minipage}[t]{0.23\textwidth}
\includegraphics[width=1.15\textwidth]{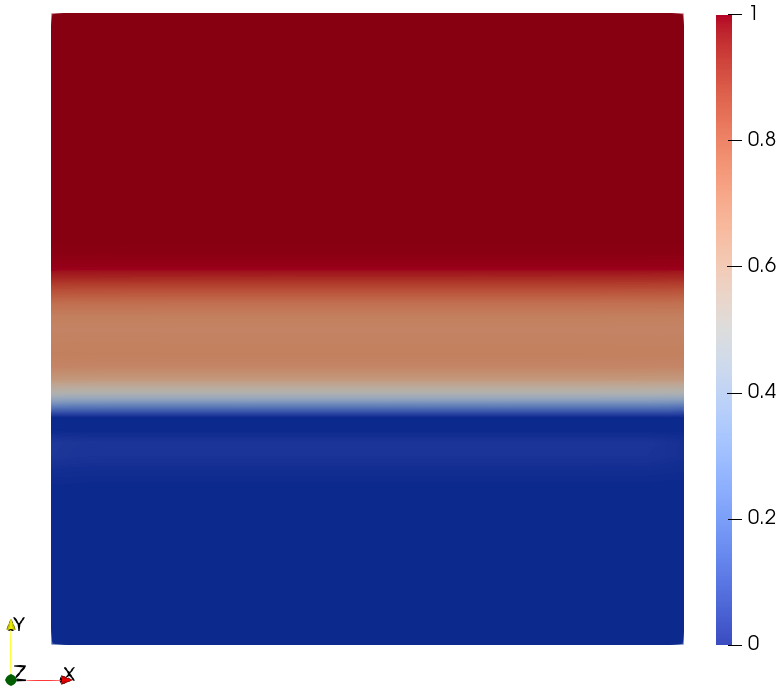}
\end{minipage}
\hfill
\begin{minipage}[t]{0.23\textwidth}
\includegraphics[width=1.15\textwidth]{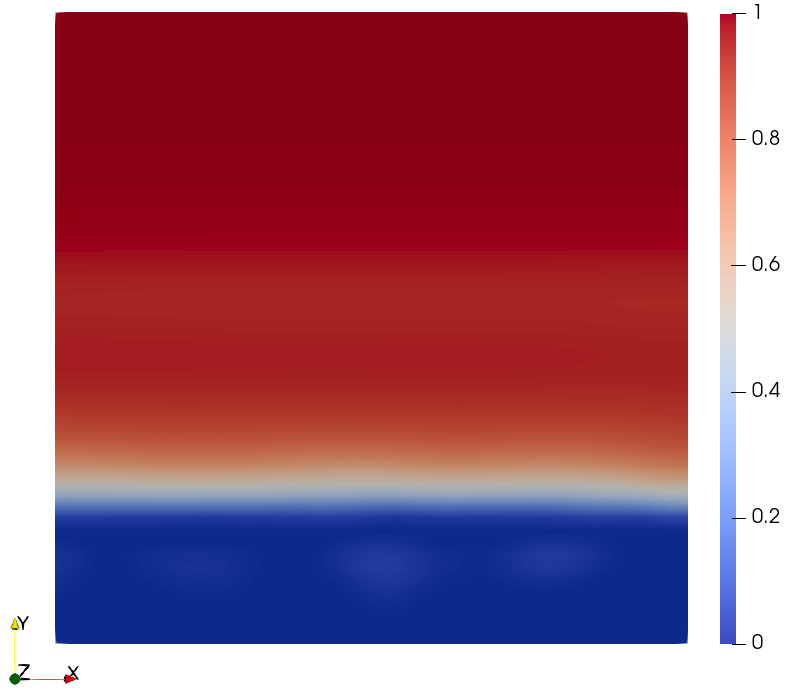}
\end{minipage}
\hfill
\begin{minipage}[t]{0.23\textwidth}
\includegraphics[width=1.15\textwidth]{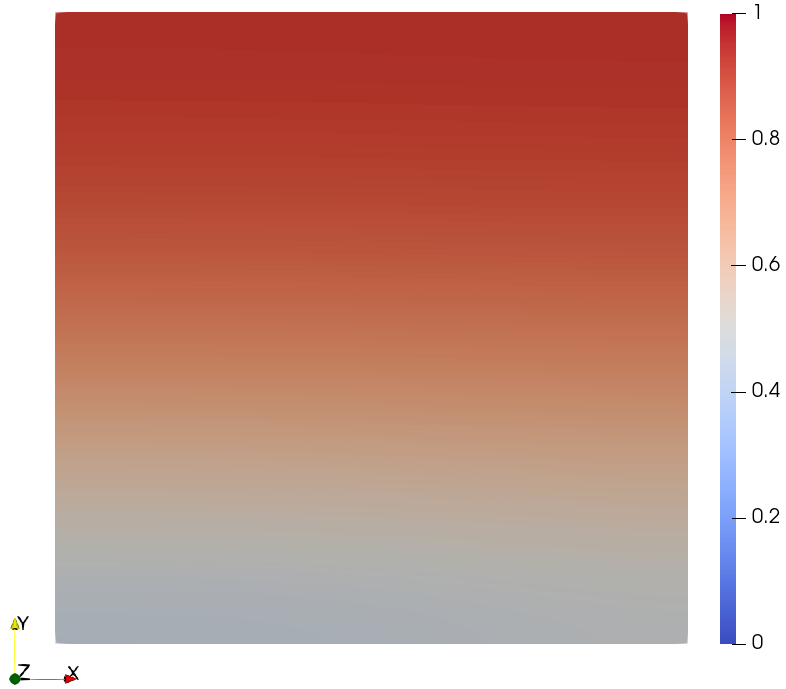}
\end{minipage}

\begin{minipage}[t]{0.23\textwidth}
\includegraphics[width=1.15\textwidth]{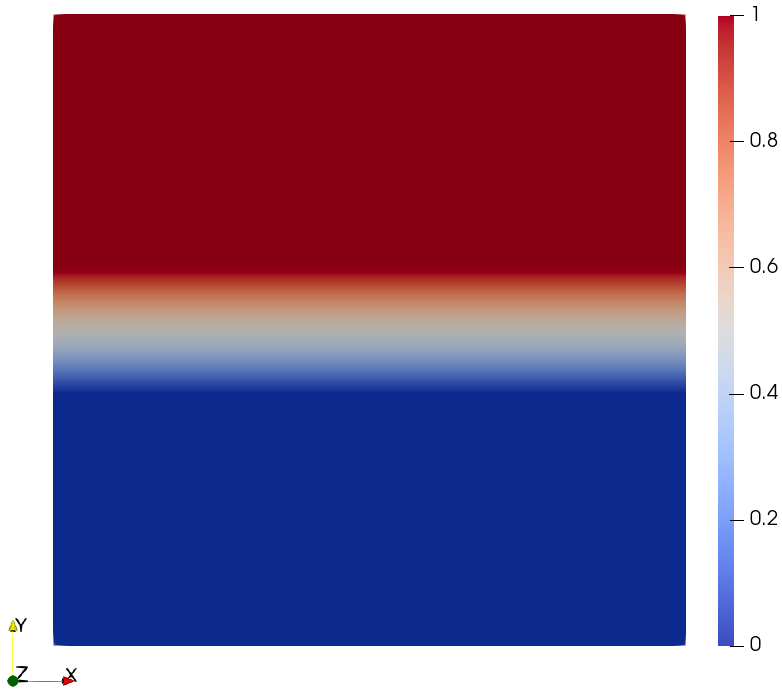}
\end{minipage}
\hfill
\begin{minipage}[t]{0.23\textwidth}
\includegraphics[width=1.15\textwidth]{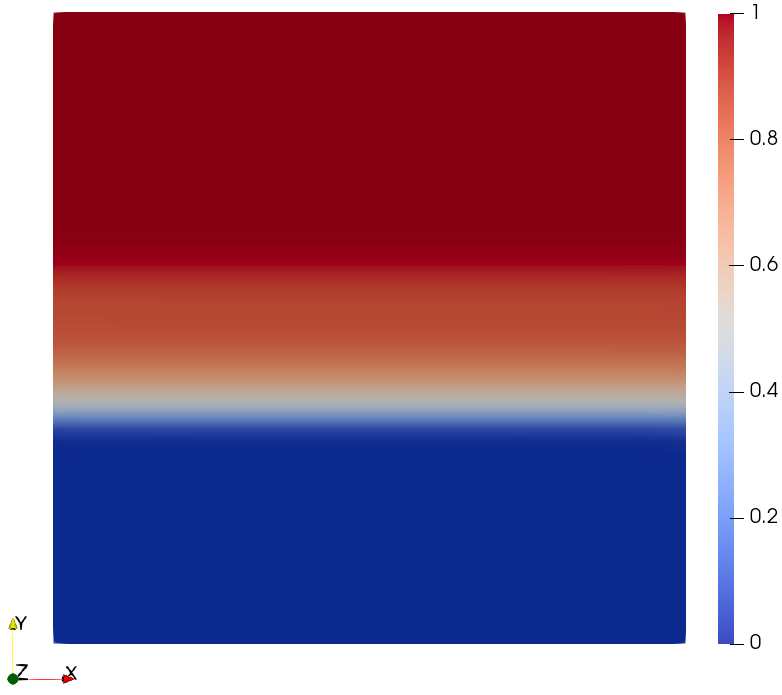}
\end{minipage}
\hfill
\begin{minipage}[t]{0.23\textwidth}
\includegraphics[width=1.15\textwidth]{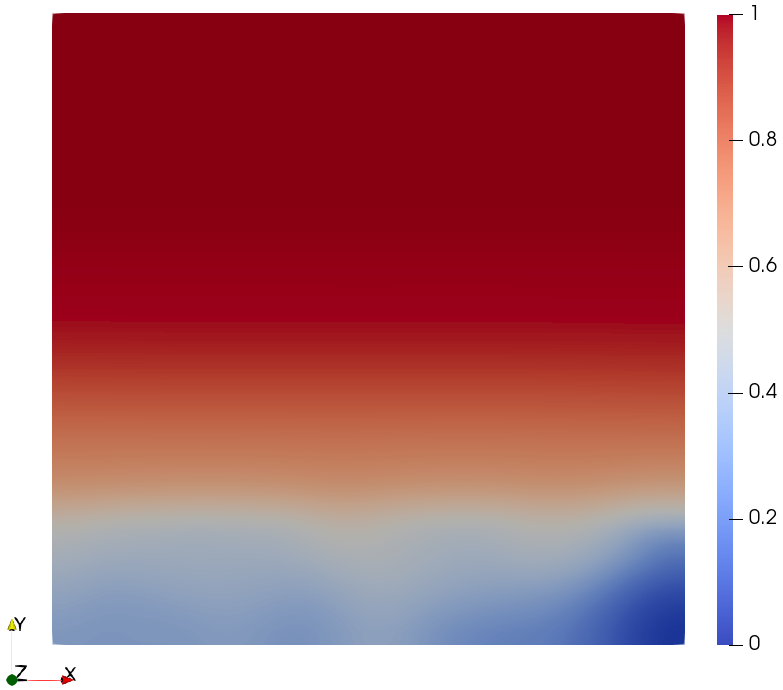}
\end{minipage} 
\hfill
\begin{minipage}[t]{0.23\textwidth}
\includegraphics[width=1.15\textwidth]{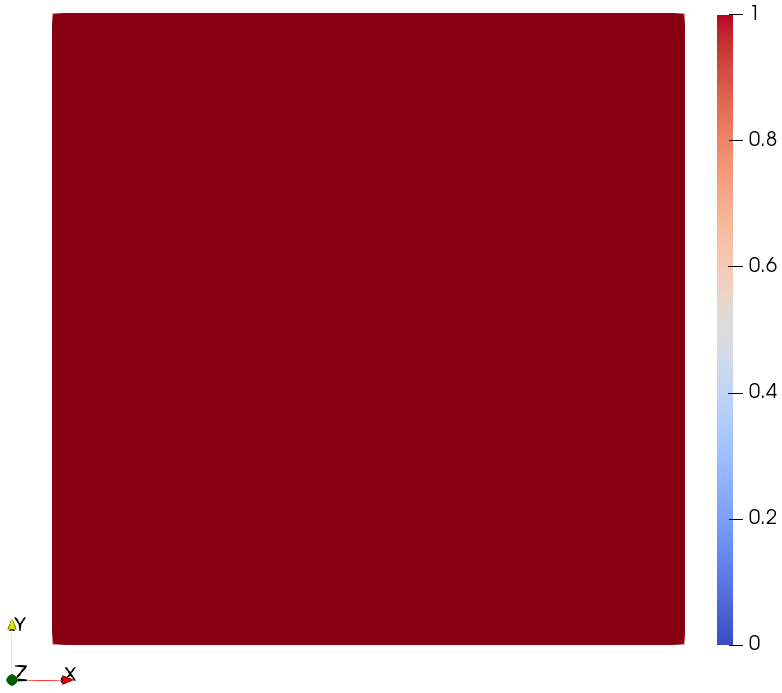}
\end{minipage}
\caption{Typical volume fraction field solvent (top), air (middle) and order parameter field vapor (bottom) during evaporation. The time step size $10^{-4}$ is used and the results shown  correspond to $t = 0$, $t =0.05$, $t = 0.15$, $t = 2.5$ (final time). The spatial discretization of a 2D domain is chosen $n_x = n_y =100$. Moreover, the parameters $\alpha_i = 10^{-8}$, $\beta_i = 10^{-1}$ for $ i \in \{ s, a \}$,  $\beta_v = 1$, $\gamma_i = 1$ for $ i \in \{ s, a, v \}$, and $\delta_v = 1$, and the molar size of the fluid $N_s=N_s =1$ are used. }
\label{fig:2D_sa_vol_t}
\end{figure}

As we observed in 1D simulations, over time the solvent evaporates and the height of the film drops and continues until a certain amount remains in the simulation box in 2D domain, see Figure~\ref{fig:2D_sa_vol_t}.

\begin{figure}[h!]
\begin{minipage}[t]{0.23\textwidth}
\includegraphics[width=1.17\textwidth]{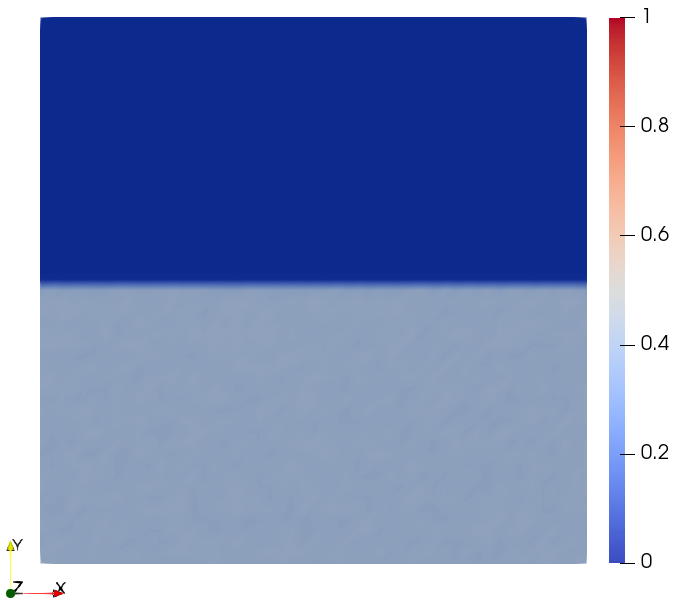}
\end{minipage}
\hfill
\begin{minipage}[t]{0.23\textwidth}
\includegraphics[width=1.17\textwidth]{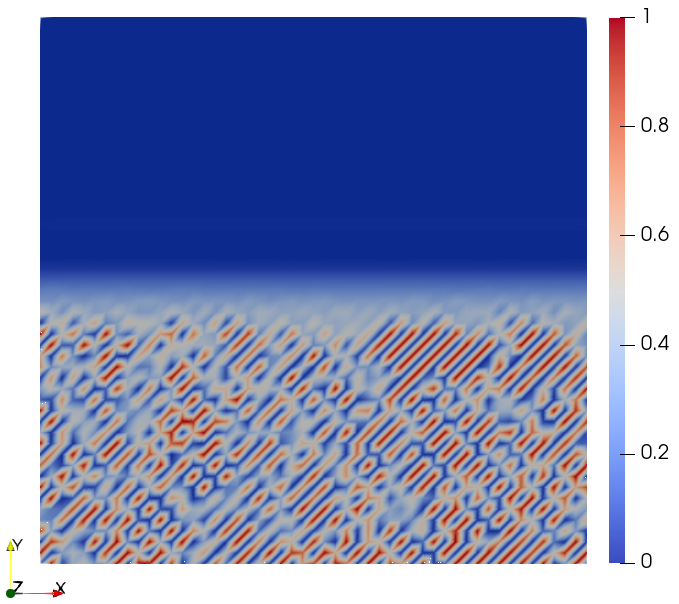}
\end{minipage}
\hfill
\begin{minipage}[t]{0.23\textwidth}
\includegraphics[width=1.17\textwidth]{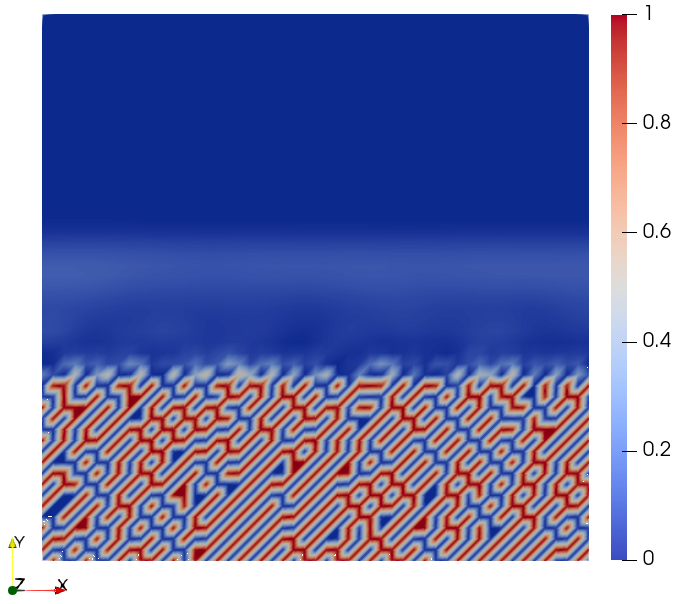}
\end{minipage}
\hfill
\begin{minipage}[t]{0.23\textwidth}
\includegraphics[width=1.17\textwidth]{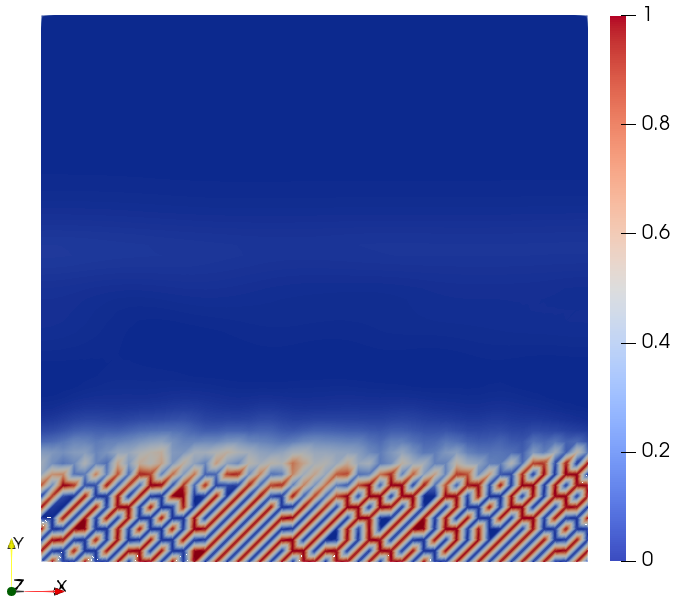}
\end{minipage}

\begin{minipage}[t]{0.23\textwidth}
\includegraphics[width=1.17\textwidth]{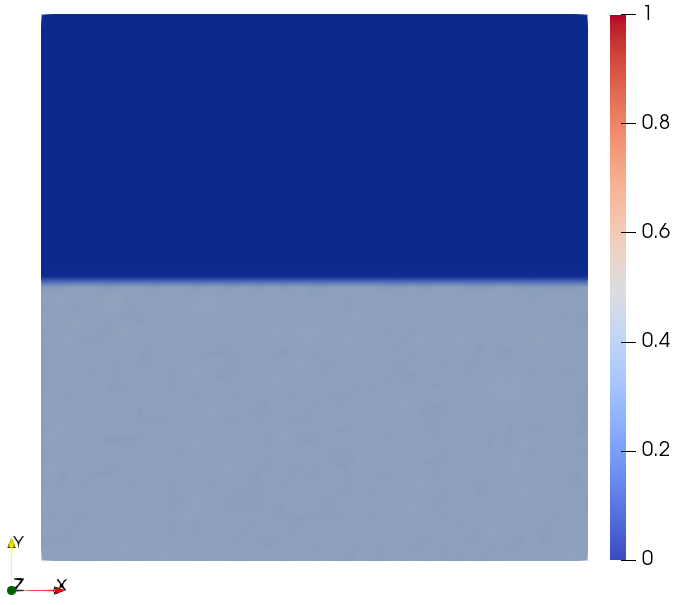}
\end{minipage}
\hfill
\begin{minipage}[t]{0.23\textwidth}
\includegraphics[width=1.17\textwidth]{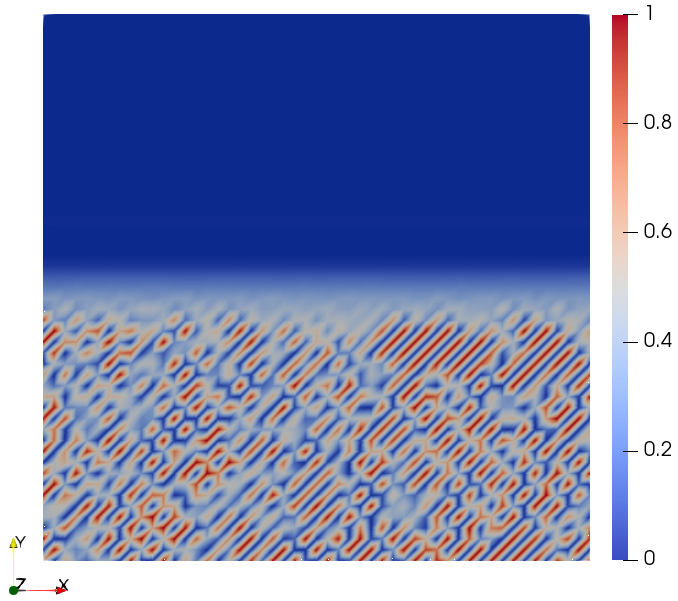}
\end{minipage}
\hfill
\begin{minipage}[t]{0.23\textwidth}
\includegraphics[width=1.17\textwidth]{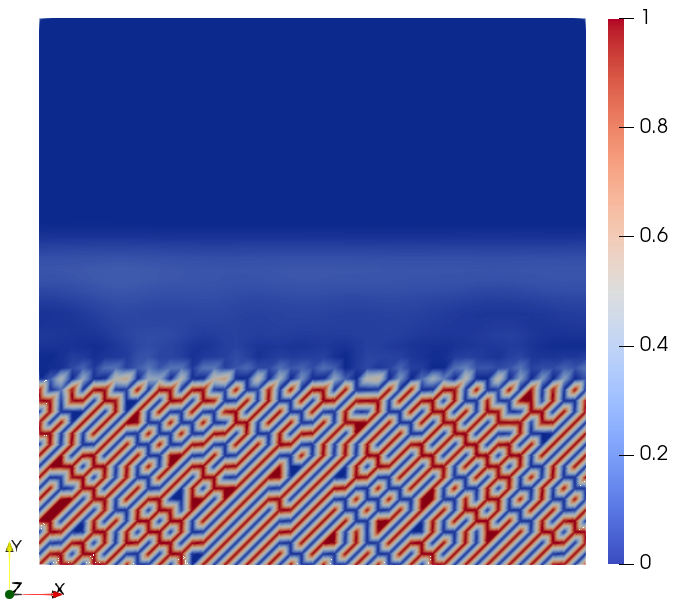}
\end{minipage}
\hfill
\begin{minipage}[t]{0.23\textwidth}
\includegraphics[width=1.2\textwidth]{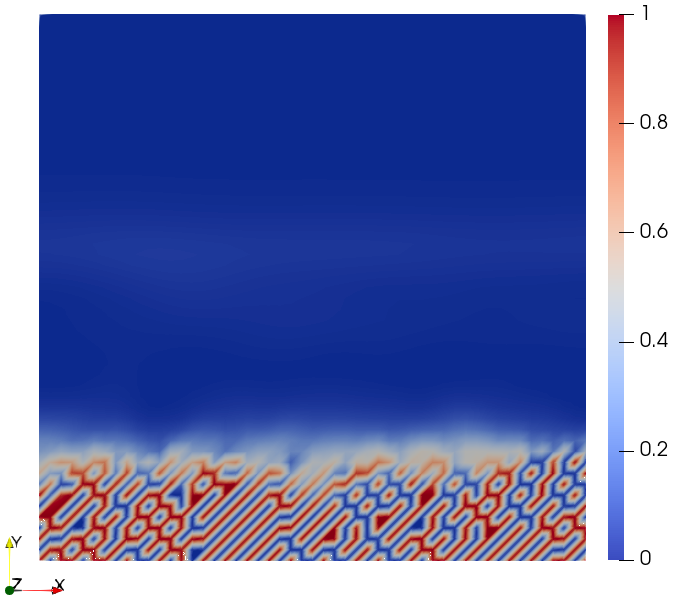}
\end{minipage}

\begin{minipage}[t]{0.23\textwidth}
\includegraphics[width=1.17\textwidth]{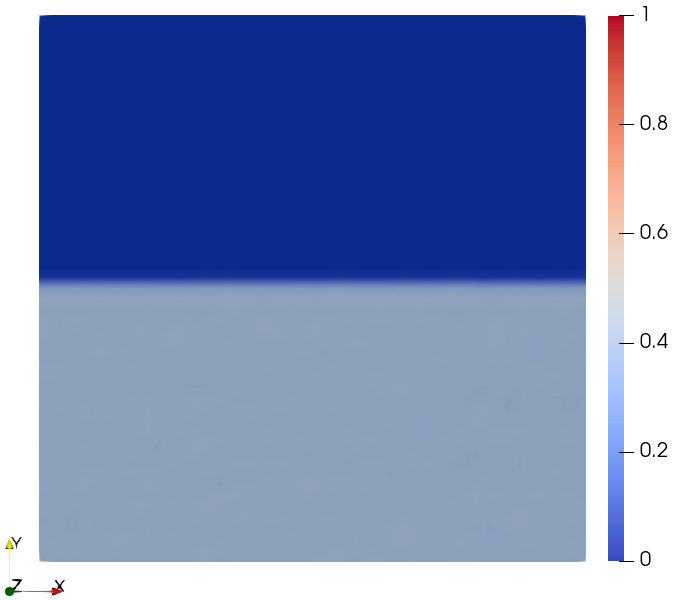}
\end{minipage}
\hfill
\begin{minipage}[t]{0.23\textwidth}
\includegraphics[width=1.17\textwidth]{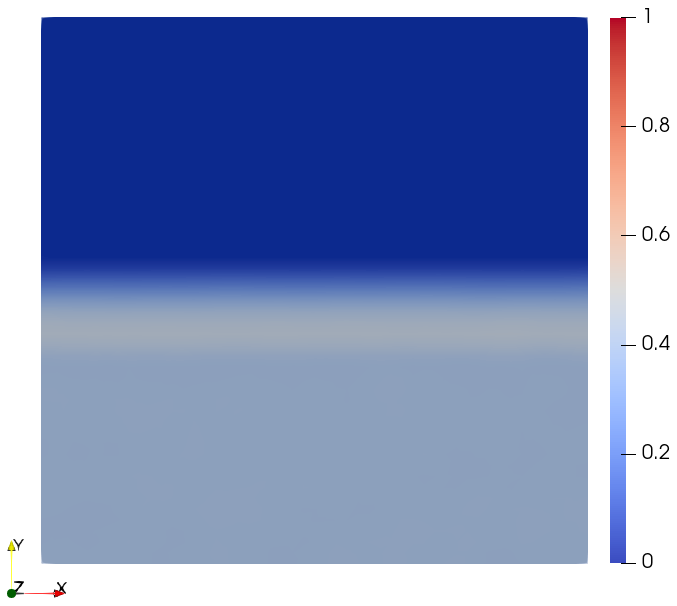}
\end{minipage}
\hfill
\begin{minipage}[t]{0.23\textwidth}
\includegraphics[width=1.17\textwidth]{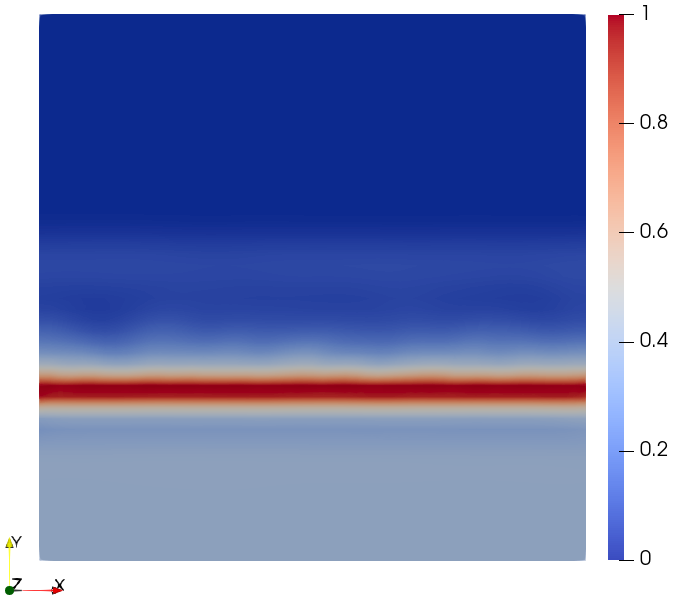}
\end{minipage}
\hfill
\begin{minipage}[t]{0.23\textwidth}
\includegraphics[width=1.17\textwidth]{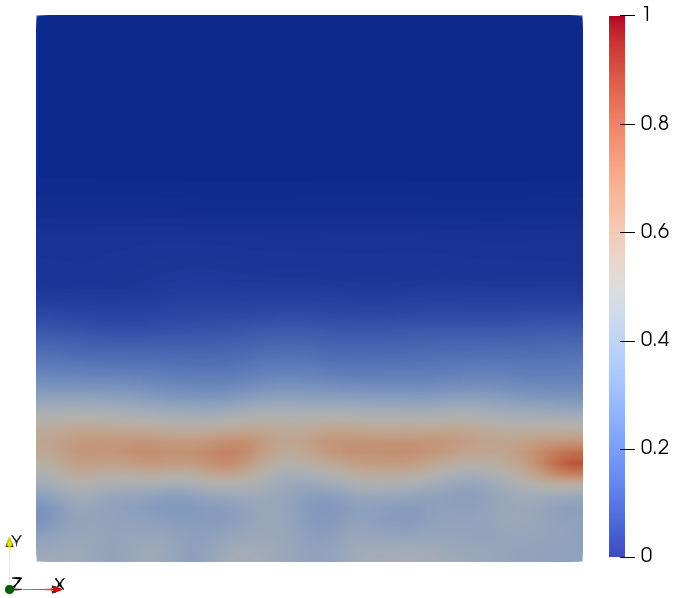}
\end{minipage}

\begin{minipage}[t]{0.23\textwidth}
\includegraphics[width=1.17\textwidth]{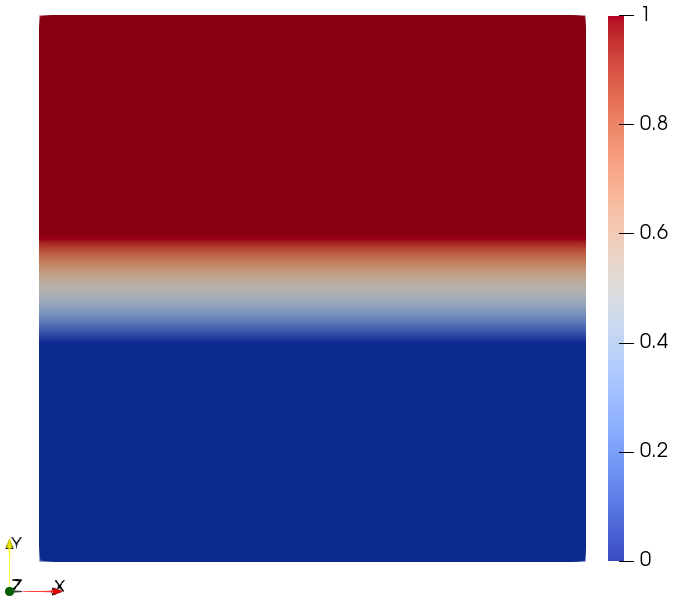}
\end{minipage}
\hfill
\begin{minipage}[t]{0.23\textwidth}
\includegraphics[width=1.17\textwidth]{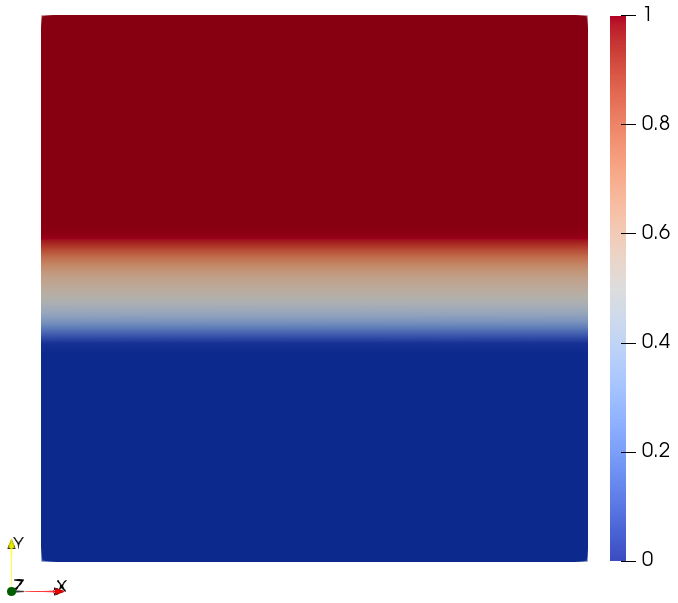}
\end{minipage}
\hfill
\begin{minipage}[t]{0.23\textwidth}
\includegraphics[width=1.17\textwidth]{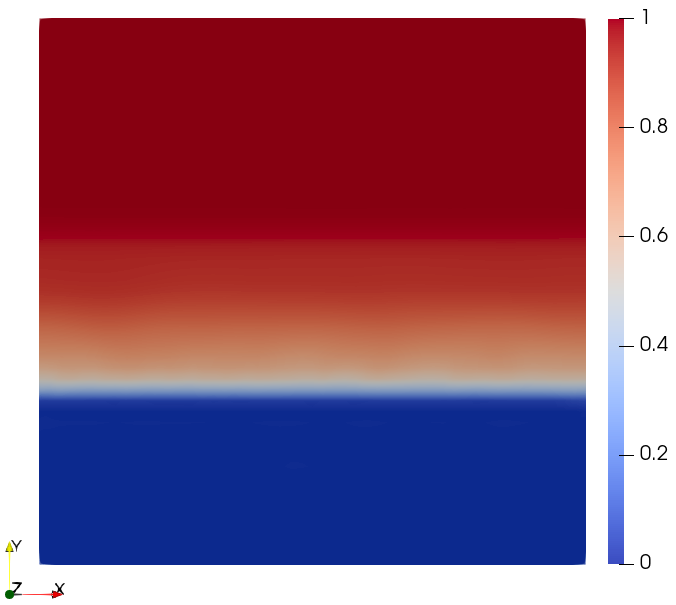}
\end{minipage}
\hfill
\begin{minipage}[t]{0.23\textwidth}
\includegraphics[width=1.17\textwidth]{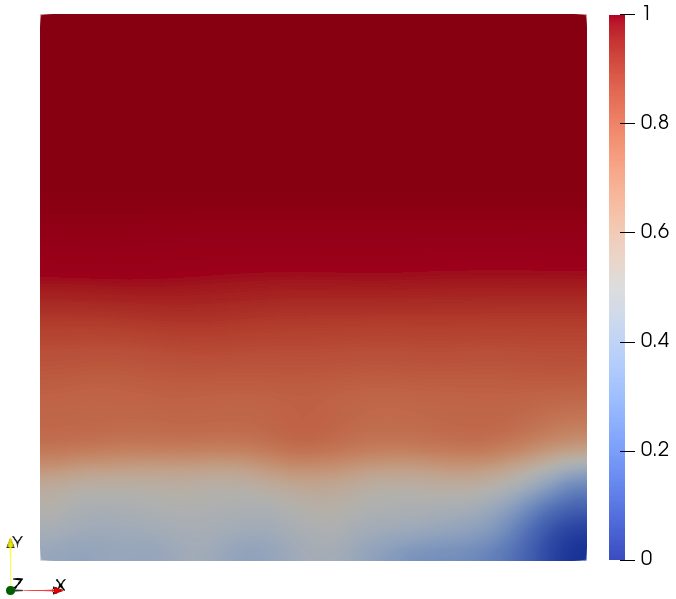}
\end{minipage}
\caption{Typical volume fraction field polymer, NFA, solvent, and order parameter field vapor (from top to bottom) during evaporation. The time step size $10^{-4}$ is used and the results shown  correspond to $t = 0$, $t =0.004$, $t = 0.05$, $t = 0.1$. The spatial discretization of a 2D domain is chosen $n_x = n_y =50$. The other parameters and setup are the same as in Figure~\ref{fig:2D_all_gmres_iter_spatial}. }
\label{fig:2D_all_vol_t}
\end{figure}







\begin{thebibliography}{10}
	
	\bibitem{abels2013incompressible}
	Helmut Abels, Daniel Depner, and Harald Garcke.
	\newblock On an incompressible {Navier}--{Stokes}/{Cahn}--{Hilliard} system
	with degenerate mobility.
	\newblock {\em Annales de l'Institut Henri Poincar{\'e} C}, 30(6):1175--1190,
	2013.
	
	\bibitem{Garcke2011_CHNS}
	Helmut Abels, Harald Garcke, and Günther Grün.
	\newblock Thermodynamically consistent, frame indifferent diffuse interface
	models for incompressible two-phase flows with different densities.
	\newblock {\em Math. Models Meth. Appl. Sci.}, 22(03):1150013, 2012.
	
	\bibitem{NAdam_FFranke_SAland_2020}
	Nadja Adam, Florian Franke, and Sebastian Aland.
	\newblock A simple parallel solution method for the navier–stokes
	cahn–hilliard equations.
	\newblock {\em Mathematics}, 8(8), 2020.
	
	\bibitem{fenics}
	Martin~S. Alnaes, Jan Blechta, Johan Hake, August Johansson, Benjamin Kehlet,
	Anders Logg, Chris~N. Richardson, Johannes Ring, Marie~E. Rognes, and
	Garth~N. Wells.
	\newblock The {FEniCS} project version 1.5.
	\newblock {\em Archive of Numerical Software}, 3, 2015.
	
	\bibitem{ufl}
	Martin~S. Aln\ae{}s, Anders Logg, Kristian~B. \O{}lgaard, Marie~E. Rognes, and
	Garth~N. Wells.
	\newblock Unified form language: A domain-specific language for weak
	formulations of partial differential equations.
	\newblock {\em ACM Trans. Math. Softw.}, 40(2), mar 2014.
	
	\bibitem{dolfinx}
	Igor~A. Baratta, Joseph~P. Dean, Jørgen~S. Dokken, Michal Habera, Jack~S.
	Hale, Chris~N. Richardson, Marie~E. Rognes, Matthew~W. Scroggs, Nathan Sime,
	and Garth~N. Wells.
	\newblock {{DOLFINx}: The next generation {FEniCS} problem solving
		environment}, December 2023.
	
	\bibitem{pyamg2023}
	Nathan Bell, Luke~N. Olson, Jacob Schroder, and Ben Southworth.
	\newblock {PyAMG}: Algebraic multigrid solvers in python.
	\newblock {\em Journal of Open Source Software}, 8(87):5495, 2023.
	
	\bibitem{benzi2005numerical}
	Michele Benzi, Gene~H Golub, and J{\"o}rg Liesen.
	\newblock Numerical solution of saddle point problems.
	\newblock {\em Acta numerica}, 14:1--137, 2005.
	
	\bibitem{bergermann_modeling_2019}
	Kai Bergermann.
	\newblock Modeling the morphology evolution of organic solar cells.
	\newblock {\em GAMM Archive for Students}, 1(1):18--27, October 2019.
	
	\bibitem{Bergermann2023}
	Kai Bergermann, Carsten Deibel, Roland Herzog, Roderick C.~I. MacKenzie,
	Jan-Frederik Pietschmann, and Martin Stoll.
	\newblock Preconditioning for a phase-field model with application to
	morphology evolution in organic semiconductors.
	\newblock {\em Communications in Computational Physics}, 34(1):1--17, 2023.
	
	\bibitem{JBosch_CKahle_MStoll_2018}
	Jessica Bosch, Christian Kahle, and Martin Stoll.
	\newblock Preconditioning of a coupled {Cahn}-{Hilliard} {Navier}-{Stokes}
	system.
	\newblock {\em Communications in Computational Physics}, 23(2):603--628, 2018.
	
	\bibitem{bosch2015preconditioning}
	Jessica Bosch and Martin Stoll.
	\newblock Preconditioning for vector-valued {Cahn}-{Hilliard} equations.
	\newblock {\em SIAM Journal on Scientific Computing}, 37(5):S216--S243, 2015.
	
	\bibitem{bosch2014fast}
	Jessica Bosch, Martin Stoll, and Peter Benner.
	\newblock Fast solution of {Cahn}-{Hilliard} variational inequalities using
	implicit time discretization and finite elements.
	\newblock {\em Journal of Computational Physics}, 262:38--57, 2014.
	
	\bibitem{boyanova2012efficient}
	Petia Boyanova, Minh Do-Quang, and Maya Neytcheva.
	\newblock Efficient preconditioners for large scale binary {Cahn}-{Hilliard}
	models.
	\newblock {\em Computational Methods in Applied Mathematics}, 12(1):1--22,
	2012.
	
	\bibitem{boyanova2014efficient}
	Petia Boyanova and Maya Neytcheva.
	\newblock Efficient numerical solution of discrete multi-component
	{Cahn}--{Hilliard} systems.
	\newblock {\em Computers \& Mathematics with Applications}, 67(1):106--121,
	2014.
	
	\bibitem{boyer2010cahn}
	Franck Boyer, C{\'e}line Lapuerta, Sebastian Minjeaud, Bruno Piar, and Michel
	Quintard.
	\newblock {Cahn}--{Hilliard}/{Navier}--{Stokes} model for the simulation of
	three-phase flows.
	\newblock {\em Transport in Porous Media}, 82:463--483, 2010.
	
	\bibitem{brenner2008mathematical}
	Susanne~C. Brenner and Larkin~R. Scott.
	\newblock {\em The mathematical theory of finite element methods}.
	\newblock Springer, 2008.
	
	\bibitem{Elliott1989}
	C.~M. Elliott.
	\newblock {\em The Cahn-Hilliard Model for the Kinetics of Phase Separation},
	pages 35--73.
	\newblock Birkh{\"a}user Basel, Basel, 1989.
	
	\bibitem{elman2014finite}
	Howard Elman, David Silvester, and Andy Wathen.
	\newblock {\em Finite elements and fast iterative solvers: with applications in
		incompressible fluid dynamics}.
	\newblock OUP Oxford, 2014.
	
	\bibitem{falgout2006introduction}
	Robert~D Falgout.
	\newblock An introduction to algebraic multigrid.
	\newblock Technical report, Lawrence Livermore National Lab.(LLNL), Livermore,
	CA (United States), 2006.
	
	\bibitem{HGarke_MHinze_CKahle_2016}
	Harald Garcke, Michael Hinze, and Christian Kahle.
	\newblock A stable and linear time discretization for a thermodynamically
	consistent model for two-phase incompressible flow.
	\newblock {\em Applied Numerical Mathematics}, 99:151--171, 2016.
	
	\bibitem{gurtin_generalized_1996}
	Morton~E. Gurtin.
	\newblock Generalized {Ginzburg}-{Landau} and {Cahn}-{Hilliard} equations based
	on a microforce balance.
	\newblock {\em Physica D: Nonlinear Phenomena}, 92(3-4):178--192, May 1996.
	
	\bibitem{janssen2007optimization}
	G~Janssen, Aranzazu Aguirre, E~Goovaerts, Peter Vanlaeke, Jef Poortmans, and
	Jean Manca.
	\newblock Optimization of morphology of p3ht/pcbm films for organic solar
	cells: effects of thermal treatments and spin coating solvents.
	\newblock {\em The European Physical Journal-Applied Physics}, 37(3):287--290,
	2007.
	
	\bibitem{kay2007efficient}
	David Kay and Richard Welford.
	\newblock Efficient numerical solution of
	{Cahn}--{Hilliard}--{Navier}--{Stokes} fluids in 2d.
	\newblock {\em SIAM Journal on Scientific Computing}, 29(6):2241--2257, 2007.
	
	\bibitem{autofenics}
	Anders Logg, Kent-Andre Mardal, Garth~N. Wells, et~al.
	\newblock {\em Automated Solution of Differential Equations by the Finite
		Element Method}.
	\newblock Springer, 2012.
	
	\bibitem{murphy2000note}
	Malcolm~F Murphy, Gene~H Golub, and Andrew~J Wathen.
	\newblock A note on preconditioning for indefinite linear systems.
	\newblock {\em SIAM Journal on Scientific Computing}, 21(6):1969--1972, 2000.
	
	\bibitem{pearson2012regularization}
	John~W Pearson, Martin Stoll, and Andrew~J Wathen.
	\newblock Regularization-robust preconditioners for time-dependent
	pde-constrained optimization problems.
	\newblock {\em SIAM Journal on Matrix Analysis and Applications},
	33(4):1126--1152, 2012.
	
	\bibitem{brabec2020_evaporation}
	Olivier J.~J. Ronsin, DongJu Jang, Hans-Joachim Egelhaaf, Christoph~J. Brabec,
	and Jens Harting.
	\newblock A phase-field model for the evaporation of thin film mixtures.
	\newblock {\em Phys. Chem. Chem. Phys.}, 22:6638--6652, 2020.
	
	\bibitem{brabec2021_evaporation}
	Olivier J.~J. Ronsin, DongJu Jang, Hans-Joachim Egelhaaf, Christoph~J. Brabec,
	and Jens Harting.
	\newblock Phase-field simulation of liquid–vapor equilibrium and evaporation
	of fluid mixtures.
	\newblock {\em ACS Applied Materials \& Interfaces}, 13(47):55988--56003, 2021.
	\newblock PMID: 34792348.
	
	\bibitem{ruge1987algebraic}
	John~W Ruge and Klaus St{\"u}ben.
	\newblock Algebraic multigrid.
	\newblock In {\em Multigrid methods}, pages 73--130. SIAM, 1987.
	
	\bibitem{saad2003iterative}
	Yousef Saad.
	\newblock {\em Iterative methods for sparse linear systems}.
	\newblock SIAM, 2003.
	
	\bibitem{basix}
	Matthew~W. Scroggs, Igor~A. Baratta, Chris~N. Richardson, and Garth~N. Wells.
	\newblock Basix: a runtime finite element basis evaluation library.
	\newblock {\em Journal of Open Source Software}, 7(73):3982, 2022.
	
	\bibitem{JShen_XYang_2010}
	Jie Shen and Xiaofeng Yang.
	\newblock A phase-field model and its numerical approximation for two-phase
	incompressible flows with different densities and viscosities.
	\newblock {\em SIAM Journal on Scientific Computing}, 32(3):1159--1179, 2010.
	
	\bibitem{strang2008analysis}
	Gilbert Strang and George Fix.
	\newblock {\em An Analysis of the Finite Element Methods, and Engineering}.
	\newblock SIAM, 2008.
	
	\bibitem{CTaylor_PHood_1973}
	C.~Taylor and P.~Hood.
	\newblock A numerical solution of the navier-stokes equations using the finite
	element technique.
	\newblock {\em Computers \& Fluids}, 1(1):73--100, 1973.
	
	\bibitem{MDijk_AWakker_1997}
	M.A. van Dijk and A.~Wakker.
	\newblock {\em Concepts of Polymer Thermodynamics}.
	\newblock Lancaster: Technomic Pub., 1997.
	
	\bibitem{SciPy}
	Pauli Virtanen, Ralf Gommers, Travis~E. Oliphant, Matt Haberland, Tyler Reddy,
	David Cournapeau, Evgeni Burovski, Pearu Peterson, Warren Weckesser, Jonathan
	Bright, St{\'e}fan~J. {van der Walt}, Matthew Brett, Joshua Wilson, K.~Jarrod
	Millman, Nikolay Mayorov, Andrew R.~J. Nelson, Eric Jones, Robert Kern, Eric
	Larson, C~J Carey, {\.I}lhan Polat, Yu~Feng, Eric~W. Moore, Jake
	{VanderPlas}, Denis Laxalde, Josef Perktold, Robert Cimrman, Ian Henriksen,
	E.~A. Quintero, Charles~R. Harris, Anne~M. Archibald, Ant{\^o}nio~H. Ribeiro,
	Fabian Pedregosa, Paul {van Mulbregt}, and {SciPy 1.0 Contributors}.
	\newblock {{SciPy} 1.0: Fundamental Algorithms for Scientific Computing in
		Python}.
	\newblock {\em Nature Methods}, 17:261--272, 2020.
	
	\bibitem{weitz_revealing_2023}
	Paul Weitz, Vincent~Marc Le~Corre, Xiaoyan Du, Karen Forberich, Carsten Deibel,
	Christoph~J. Brabec, and Thomas Heumüller.
	\newblock Revealing {Photodegradation} {Pathways} of {Organic} {Solar} {Cells}
	by {Spectrally} {Resolved} {Accelerated} {Lifetime} {Analysis}.
	\newblock {\em Advanced Energy Materials}, 13(2):2202564, 2023.
	
	\bibitem{Wodo2012_height}
	Olga Wodo and Baskar Ganapathysubramanian.
	\newblock Modeling morphology evolution during solvent-based fabrication of
	organic solar cells.
	\newblock {\em Computational Materials Science}, 55:113--126, 2012.
	
	\bibitem{wopke_traps_2022}
	Christopher Wöpke, Clemens Göhler, Maria Saladina, Xiaoyan Du, Li~Nian,
	Christopher Greve, Chenhui Zhu, Kaila~M. Yallum, Yvonne~J. Hofstetter, David
	Becker-Koch, Ning Li, Thomas Heumüller, Ilya Milekhin, Dietrich R.~T. Zahn,
	Christoph~J. Brabec, Natalie Banerji, Yana Vaynzof, Eva~M. Herzig, Roderick
	C.~I. MacKenzie, and Carsten Deibel.
	\newblock Traps and transport resistance are the next frontiers for stable
	non-fullerene acceptor solar cells.
	\newblock {\em Nature Communications}, 13(1):3786, July 2022.
	
	\bibitem{xu2017algebraic}
	Jinchao Xu and Ludmil Zikatanov.
	\newblock Algebraic multigrid methods.
	\newblock {\em Acta Numerica}, 26:591--721, 2017.
	
	\bibitem{zhang2022renewed}
	Guichuan Zhang, Francis~R Lin, Feng Qi, Thomas Heum\"uller, Andreas Distler,
	Hans-Joachim Egelhaaf, Ning Li, Philip~CY Chow, Christoph~J Brabec, Alex K-Y
	Jen, et~al.
	\newblock Renewed prospects for organic photovoltaics.
	\newblock {\em Chemical Reviews}, 122(18):14180--14274, 2022.
	
\end{thebibliography}

\clearpage

\end{document}